\documentstyle[12pt,a4wide]{article}

\input amssym.def
\input amssym

\font\teneufm=eufm10 scaled \magstep1
\font\seveneufm=eufm7 scaled \magstep1
\font\fiveeufm=eufm5 scaled \magstep1
\newfam\eufmfam
\textfont\eufmfam=\teneufm
\scriptfont\eufmfam=\seveneufm
\scriptscriptfont\eufmfam=\fiveeufm
\def\frak#1{{\fam\eufmfam\relax#1}}

\def\hfl#1#2{\smash{\mathop{\hbox to 12mm{\rightarrowfill}}
\limits^{\scriptstyle#1}_{\scriptstyle#2}}}
\def\vfl#1#2{\llap{$\scriptstyle #1$}\left\downarrow
\vbox to 6mm{}\right.\rlap{$\scriptstyle #2$}}
\def\hgfl#1#2{\smash{\mathop{\hbox to 12mm{\leftarrowfill}}
\limits^{\scriptstyle#1}_{\scriptstyle#2}}}
\def\vhfl#1#2{\llap{$\scriptstyle #1$}\left\uparrow
\vbox to 6mm{}\right.\rlap{$\scriptstyle #2$}}

\def\DEMONSTRATION{\smallskip{\noindent{\it D{\'e}monstration}.\ }}

\newtheorem{proposition}{{\sc Proposition}}[subsection]

\newtheorem{corollaire}[proposition]{{\sc Corollaire}}
\newtheorem{definition}[proposition]{{\sc D{\'e}finition}}
\newtheorem{lemme}[proposition]{{\sc Lemme}}
\newtheorem{theoreme}{{\sc Th{\'e}or{\`e}me}}

\def\carre{\square\bigskip}

\let\[\lbrack
\let\]\rbrack
\let\<\langle
\let\>\rangle
\let\rta\rightarrow

\def\AA{{\Bbb A}}
\def\NN{{\Bbb N}}
\def\CC{{\Bbb C}}
\def\FF{{\Bbb F}}
\def\ZZ{{\Bbb Z}}
\def\PP{{\Bbb P}}
\def\RR{{\rm R}}
\def\Ga{{\Bbb G}_a}
\def\Gm{{\Bbb G}_m}
\def\Ql{{\bar{\Bbb Q}}_\ell}
\def\A{{\cal A}}
\def\L{{\cal L}}
\def\F{{F}}

\def\OO{{\cal O}}

\def\R{{\cal R}}

\def\t{{\rm t}}
\def\Id{{\rm Id}}
\def\Fr{{\rm Fr}}

\def\Tr{{\rm Tr}}

\def\S{{\frak S}}

\def\Fix{{\rm Fix\,}}

\def\End{{\rm End\,}}
\def\Gal{{\rm Gal\,}}
\def\Lie{{\rm Lie\,}}
\def\Nil{{\rm Nil}}

\def\det{{\rm det}}
\def\pr{{\rm pr}}

\def\act{{\rm act}}
\def\res{{\rm res\,}}

\def\v{{\rm val}}
\def\d{{\rm d}}
\def\diag{{\rm diag\,}}
\def\GL{{\rm GL}}

\def\lg{{\rm lg}}
\def\vp{\varpi}

\def\mod{{\rm mod\,\,}}

\def\gl{{\frak {gl}}}
\def\g{{\frak g}}
\def\IC{{\rm IC}}
\def\und{\underline}

\def\Res{{\rm Res}}
\def\Ind{{\rm Ind}}
\def\pgcd{{\rm pgcd}\,}
\def\C{{\cal C}}
\def\H{{\cal H}}
\def\tA{{\tilde{\cal A}}}

\def\limpro{\mathop{\oalign{lim\cr\hidewidth$\longleftarrow$\hidewidth\cr}}}

\font\tenmsb=msbm10 scaled \magstep1
\font\sevenmsb=msbm7 scaled \magstep1
\font\fivemsb=msbm5 scaled \magstep1
\newfam\msbfam
\textfont\msbfam=\tenmsb
\scriptfont\msbfam=\sevenmsb
\scriptscriptfont\msbfam=\fivemsb
\def\Bbb#1{{\fam\msbfam\relax#1}}

\def\diagramme#1{\def\normalbaselines{\baselineskip=1pt
\lineskip=6pt\lineskiplimit=1pt} \matrix{#1}}

\def\nbc{Ng\^o Ba\hskip -6pt\raise 2pt\hbox{'}\hskip 2pt o Ch\^au}

\begin{document}
\title{Faisceaux pervers,\\ homomorphisme de changement de base\\ 
et lemme fondamental de Jacquet et Ye}
\author{\nbc}
\date{}
\maketitle

\begin{abstract}
We give a geometric interpretation of the base change homomorphism 
between the Hecke algebra of $\GL(n)$ for an unramified extension
of local fields of positive characteristic.
For this, we use some results of Ginzburg, Mirkovic and Vilonen related to
the geometric Satake isomorphism. We give new proof of these results
in the positive characteristic case.

By using that geometric interpretation of the base change homomorphism,
we prove the fundamental lemma of Jacquet and Ye for arbitrary Hecke function
in the equal characteristic case.
\end{abstract}

\section*{Introduction}

Soient $F$ un corps local de caract\'erisitique $p>0$, 
${\cal O}$ son anneau des entiers et $k=\FF_q$ son corps r\'esiduel.
Notons ${\cal H}^+$ l'alg{\`e}bre des fonctions complexes {\`a} support 
compact dans
$$\GL(n,F)^+=\GL(n,F)\cap\gl(n,{\cal O})$$ 
qui sont bi-$\GL(n,{\cal O})$-invariantes.
D'apr{\`e}s Satake (\cite{Sa}), on a un isomorphisme 
entre ${\cal H}^+$ et l'alg{\`e}bre des polyn{\^o}mes
sym{\'e}triques :
$${\cal H}^+{\tilde\rta}\,\CC[z_1,\ldots,z_n]^{{\frak S}_n}.$$

Soient $r$ un entier naturel, 
$F_r$ l'extension non ramifi{\'e}e de degr{\'e} $r$ de $F$,
${\cal O}_r$ son anneau des entiers et $k_r=\FF_{q^r}$ son corps r\'esiduel. 
Notons ${\cal H}_r^+$
l'alg{\`e}bre des fonctions complexes {\`a} support compact dans
$\GL(n,F_r)^+$ qui sont bi-$\GL(n,{\cal O}_r)$-invariantes.
On a comme pr\'ec\'edemment un isomorphisme de Satake
$${\cal H}^+_r{\tilde\rta}\,\CC[t_1,\ldots,t_n]^{{\frak S}_n}.$$

Compte tenu des isomorphismes de Satake pour $\H^+$ et pour $\H^+_r$,
l'homo\-morphisme de changement de base $b:{\cal H}^+_r\rta{\cal H}^+$
est d\'efini par l'homomor\-phisme
$$\CC[t_1,\ldots,t_n]^{{\frak S}_n}\rta\CC[z_1,\ldots,z_n]^{{\frak S}_n}$$
qui envoie $t_i$ sur $z_i^r$.

D'apr\`es la d\'ecomposition de Cartan
$$\GL(n,F_r)^+=\coprod_{\lambda=(\lambda_1\geq\cdots\geq\lambda_n\geq 0)}
\GL(n,\OO_r)\vp^\lambda\GL(n,\OO_r)$$
$\vp^\lambda$ \'etant la matrice diagonale 
$\diag(\vp^\lambda_1,\ldots,\vp^\lambda_n)$,
les fonctions caract\'eristiques $c_{r,\lambda}$ des doubles classes
$\GL(n,\OO_r)\vp^\lambda\GL(n,\OO_r)$ forment une base de $\H^+_r$.
On ne conna\^{\i}t pas d'expression explicite pour
les fonctions $b(c_{r,\lambda})$
mis \`a part le cas trivial 
$\lambda=(0,\ldots,0)$ o\`u $b(c_{r,\lambda})=c_\lambda$ et 
le cas $\lambda=(1,0,\ldots,0)$
o\`u $b(c_{r,\lambda})$ est la fonction de Drinfeld (\cite{Lau1}).

On sait d'apr{\`e}s Lusztig (\cite{Lus}) que
$(\GL(n,F)\cap\gl(n,{\cal O}))/\GL(n,{\cal O})$ s'identifie naturellement {\`a}
l'ensemble des points rationnels d'un sch{\'e}ma $X$ 
qui est une r{\'e}union disjointe 
de $k$-sch{\'e}mas projectifs $X_d$. L'action de $\GL(n,{\cal O})$
sur $(\GL(n,F)\cap\gl(n,{\cal O}))/\GL(n,{\cal O})$ se d{\'e}duit 
d'une action d'un groupe
alg{\'e}brique de dimension infinie $G$ sur $X$, $G$ agissant sur chacune des
composante connexe $X_d$ {\`a} travers un quotient $G_d$ de type fini sur $k$.

On peut param{\'e}trer les orbites de $X_d$ par les $n$-partitions
$\lambda=(\lambda_1\geq\lambda_2\geq\cdots\geq\lambda_n\geq 0)$ de $d$.
La d\'ecomposition en orbites 
$X_d=\coprod_{|\lambda|=d}X_\lambda$
refl\`ete bien entendu la d\'ecomposition de Cartan.
L'adh\'erence ${\bar X}_\lambda$ de l'orbite $X_\lambda$ 
\'etant en g\'en\'eral singuli\`ere, il est naturel de consid\'erer
son complexe d'intersention $\ell$-adique 
$\A_\lambda=\IC({\bar X}_\lambda,\Ql)$.
Le faisceau pervers $\A_\lambda$ est alors d\'efini sur $k$ 
et $G$-\'equivariant.

On d\'efinit \`a la suite de Lusztig les fonctions
$$a_{r,\lambda}:X(k_r)\rta\Ql$$
par
$$a_{r,\lambda}(x)=\Tr(\Fr_{q^r},(\A_\lambda)_x).$$
Choisissons  
\footnote{Ce choix n'est en fait pas important car toutes nos traces sont rationnelles.} une fois pour toutes un isomorphisme $\Ql\simeq\CC$. 
On peut alors consid\'erer ces fonctions $a_{r,\lambda}$ comme 
des \'el\'ements de $\H^+_r$. La matrice de passage des fonctions 
$c_{r,\lambda}$ aux fonctions $a_{r,\lambda}$ \'etant triangulaire
sup\'erieure, les fonctions $a_{r,\lambda}$ forment aussi une base de
$\H^+_r$.

Notre premier objectif consiste \`a interpr\'eter g\'eom\'etriquement
les fonctions $b(a_{r,\lambda})$. 

A la suite de Ginzburg, Mirkovic et Vilonen, on peut d{\'e}finir un produit
de convolution des faisceaux pervers de type $\A_\lambda$,
et donc, pour chaque $\lambda$, la $r$-\`eme puissance convol\'ee
$$\A_\lambda^{*r}=
\underbrace{\A_\lambda*\cdots*\A_\lambda}_{r\,\,{\rm fois}}.$$ 
Ce produit de convolution {\'e}tant commutatif,
on a un automorphisme ${\kappa'}$ de $\A_\lambda^{*r}$ :
$${\kappa'}:
\A_{\lambda,1}*\A_{\lambda,2}*\cdots*\A_{\lambda,r}{\buildrel\kappa\over\rta}
\A_{\lambda,r}*\A_{\lambda,1}*\cdots*\A_{\lambda,r-1}{\buildrel\iota\over\rta}
 \A_{\lambda,1}*\A_{\lambda,2}*\cdots*\A_{\lambda,r}$$
o{\`u} $\A_{\lambda,1},\ldots,\A_{\lambda,r}$ sont $r$ copies de $\A_\lambda$,
o\`u $\kappa$ est l'isomorphisme de commutativit\'e de Ginzburg, Mirkovic
et Vilonen et o\`u $\iota$ se d\'eduit des isomorphismes \'evidents
$\A_{\lambda,i}\rta\A_{\lambda,i+1}$ l'indice $i$ \'etant prise parmi les 
classes modulo $r$.
On d{\'e}finit pour chaque $\lambda$, 
la fonction $\phi_{r,\lambda}\in{\cal H}^+$
comme la trace de $\Fr\circ {\kappa'}$ sur les fibres de $\A_\lambda^{*r}$
au-dessus des points fix{\'e}s par $\Fr$.

\setcounter{theoreme}{2}

\begin{theoreme}
 Pour toute $n$-partition $\lambda$, on a
$b(a_{r,\lambda})=\phi_{r,\lambda}$. 
\end{theoreme}

On utilise ce th\'eor\`eme 3 pour g{\'e}n{\'e}rali\-ser 
le r{\'e}sultat principal de \cite{Ngo1}. 
Pour l'\'enoncer pr\'ecis\'ement, fixons quelques notations.

Notons $A$ le sous-groupe diagonal de $\GL(n)$ et
$N$ son sous-groupe des matrices
triangulaires sup{\'e}rieures unipotentes.
Notons $\theta:N(F)\rta\Ql^\times$ le caract{\`e}re
$$\theta(x)=\Psi\bigl(\sum_{i=1}^{n-1}\res(x_{i,i+1})\bigr)$$
o{\`u} $\Psi:F\rta\Ql^\times$ est un caract{\`e}re additif de conducteur $\OO$.

Pour toute fonction $\phi\in{\cal H}^+$ et pour
tout $a\in A(F)$, posons
$$I(a,\phi)=\int_{N(F)\times N(F)}
        \phi(\,^\t x_1ax_2)\theta(x_1)
        \theta(x_2)\d x_1\d x_2$$
o{\`u} la mesure de Haar normalis{\'e} $\d x$
de $N(F)$ attribue {\`a} $N({\cal O})$ le volume $1$.

Cette int{\'e}grale intervient comme une int{\'e}grale orbitale dans
une formule des traces relative de Jacquet.
Il s'agit d'une int{\'e}grale de Kloosterman
si $\phi$ est la fonction caract{\'e}ristique de $\GL(n,{\cal O})$.

Soient $F_2$ l'extension quadratique non ramifi{\'e}e de $F$ et 
$\OO_{2}$ son anneau des entiers. 
Notons $\theta': N(F_{2})\rta\Ql^\times$ le caract{\`e}re d{\'e}fini par
$$\theta'(x)=\Psi\bigl(\sum_{i=1}^{n-1}
(x_{i,i+1}+{\bar x}_{i,i+1}))\bigr)$$
o{\`u} $x\mapsto{\bar x}$ est l'{\'e}l{\'e}ment non trivial du groupe 
$\Gal(F_{2}/F)$.

Soit $S(F)$ l'ensemble des matrices $g\in\GL(n,F_{2})$
telles que $^\t{\bar g}=g$. 
Le groupe $\GL(n,F_{2})$ agit sur $S(F)$ par $g.s=\,^t{\bar g}sg$.
Notons ${\cal H'}^+$ l'espace des fonctions
{\`a} support compact dans 
$$S(F)^+=S(F)\cap\gl(n,\OO_{2})$$ 
qui sont invariantes sous l'action de $\GL(n,\OO_{2})$. 

Pour toute fonction $\phi'\in{\cal H'}^+$
et pour toute matrice diagonale $a\in A(F)$,
posons
$$J(a,\phi')=\int_{N(F_2)}
\phi'(\,^\t{\bar x}ax)\theta'(x)\d x$$
o{\`u} la mesure de Haar normalis{\'e}e $\d x$ de
$N(F_2)$ attribue {\`a} $N({\cal O}_2)$
le volume $1$.

Soit $b:{\cal H}^+_2\rta{\cal H}^+$ l'homomorphisme de changement de base.
Soit $b':{\cal H}^+_2\rta{\cal H'}^+$ l'application  d{\'e}finie
par $b'(f)=\phi'$ o\`u
$$\phi'(g\,^\t{\bar g})=\int_{H(\OO)}f(gh)\d h$$
o{\`u} $H(\OO)$ est le sous-groupe de $\GL(n,\OO_{2})$ form{\'e} des 
matrices $g$ telles que $^\t{\bar g}=g^{-1}$. L'application $b'$
est bien d\'efinie puisque toute matrice hermitienne $s\in S(F)$
peut s'\'ecrire sous la forme $s=g\,^\t{\bar g}$ avec $g\in\GL(n,F_2)$.

\begin{theoreme}
Pour toute matrice
$$a=\diag(a_1,a_1^{-1}a_2,\ldots,a_{n-1}^{-1}a_n)\in A(F)$$
pour toute fonction $f\in{\cal H}^+_2$, on a
$$I(a,b(f))=(-1)^{\v(a_1a_2\ldots a_{n-1})}J(a,b'(f))$$
o\`u $\v$ est la valuation de $F$.
\end{theoreme}

On d{\'e}montre finalement le lemme fondamental de Jacquet et Ye pour 
l'{\'e}l{\'e}ment long $w_0$ du groupe de Weyl $\S_n$. 
Pour toute $\phi\in{\cal H}^+$, pour toute $\phi'\in{\cal H'}^+$
et pour tout {\'e}l{\'e}ment central $a\in A(F_\vp)$, on introduit,
{\`a} la suite de Jacquet et Ye, les int\'egrales orbitales relatives
$$\displaylines{
I(w_0a,\phi)=\int_{N(F)\times N(F)/
(N(F)\times N(F))^{\scriptstyle w_0a}}
\phi(\,^\t xw_0ax')\theta(xx')\d x\d x'\, ;\cr
J(w_0a,\phi')=\int_{N(F_{2})/N(Fy_{2})^{\scriptstyle w_0a}}
\phi'(\,^\t {\bar x}w_0ax)\theta'(x)\d x\,.}$$

\begin{theoreme}
Pour un {\'e}l{\'e}ment central $a=\diag(\vp^d,\vp^d,\ldots,\vp^d)\in A(F_\vp)$
pour toute fonction $f\in{\cal H}^+_2$, on a
$$I(w_0a,b(f))=(-1)^{d(1+2+\cdots+(n-1))}J(w_0a,b'(f)).$$
\end{theoreme}

Ces {\'e}nonc{\'e}s jouent le r{\^o}le 
d'un lemme fondamental dans une formule
des traces relative. 
Ils ont {\'e}t{\'e} conjectur{\'e}s par Jacquet et Ye dans \cite{JY}
pour un corps local $F$ de caract{\'e}ristique arbitraire. 
Ils les ont d{\'e}montr{\'e}s pour $n=2$ et $n=3$ dans loc.cit. 
Dans \cite{Ngo1}, sous l'hypoth{\`e}se que la caract{\'e}ristique
de $F$ est positive, on a d{\'e}montr{\'e} le th\'eo\-r\`e\-me 4 dans
le cas particulier o{\`u} $f$ est l'unit{\'e} 
de l'alg{\`e}bre ${\cal H}^+_2$.

Gr{\^a}ce au th{\'e}or{\`e}me 3, 
on sait interpr{\'e}ter g{\'e}om{\'e}triquement $b$.
L'applica\-tion $b'$ {\'e}tant d{\'e}finie comme une int{\'e}grale 
le long des fibres,
elle admet aussi une interpr{\'e}tation g{\'e}om{\'e}\-trique.

Ayant une traduction g{\'e}om{\'e}trique des applications $b$ et $b'$, on 
d{\'e}montre le th\'eor\`eme 4 en adaptant les arguments de \cite{Ngo1}
et le th\'eor\`eme 5 en utilisant la preuve d'une conjecture
de Frenkel-Gaitsgory-Kazhdan-Vilonen (\cite{Ngo3}).

\bigskip
\bigskip

L'article est divis\'e en deux parties.

Dans la premi\`ere partie, on propose une nouvelle d\'emonstration
valable en caract\'eristique positive 
des r\'esultats de Ginzburg, Mirkovic et Vilonen (\cite{Gin},\cite{M-V})
sur lesquels s'appuierait notre th\'eor\`eme 3.
L'id\'ee principale est de d\'eformer en consid\'erant une situation
globale. On exploite aussi une analogie avec 
la correspondance de Springer.

On rappelle en 1.1 la d{\'e}finition du 
sch{\'e}ma des r{\'e}seaux de Lusztig (\cite{Lus})
ou plut{\^o}t sa variante globale. 
D{\'e}sormais, on note ${\cal O}=k[\vp]$ l'anneau des polyn{\^o}mes
{\`a} coefficients dans $k$ et {\`a} une variable $\vp$, 
$F$ son corps des fractions, ${\cal O}_\vp$
et $\F_\vp$ leurs compl{\'e}t{\'e}s en $\vp$. Notons $Q_d$
le $k$-sch{\'e}ma affine 
dont l'ensemble des $k$-points est celui des polyn{\^o}mes unitaires 
de degr{\'e} $d$ \`a coefficients dans $k$.
Soit $X_d$ le $k$-sch{\'e}ma dont l'ensemble 
des $k$-points est celui des r{\'e}seaux $\R\subset{\cal O}^n$
tel que $\dim({\cal O}^n:\R)=d$. 
On a un morphisme propre $\phi:X_d\rta Q_d$ d\'efini par le d\'eterminant.
Ce morphisme est en g{\'e}n{\'e}ral singulier en dehors
de l'ouvert $Q_{d,rss}$ des polyn{\^o}mes unitaires s{\'e}parables.
La fibre de $\phi$ en $\vp^d\in Q_d(k)$ 
qu'on notera $X_{d,\vp}$, s'identifie {\`a} la composante connexe
index\'ee par $d$ du sch{\'e}ma  des r{\'e}seaux $\R_\vp\subset{\cal O}_\vp^n$ 
consid{\'e}r{\'e} par Lusztig.
C'est la fibre la plus int{\'e}ressante.

On introduit en 1.2 en imitant \cite{Lau},
une r{\'e}solution simultan{\'e}e
$\pi:{\tilde X}_d\rta X_d$ des fibres de $\phi$. 
Le morphisme $\pi$ est petit,
g{\'e}n{\'e}riquement un $\S_d$-torseur, 
sa restriction {\`a} toutes les fibres de 
$\phi$ est semi-petite : il est compl{\`e}tement analogue {\`a} la
r{\'e}solution simultan{\'e}e de Grothendieck-Springer.

Ainsi, pour chaque repr{\'e}sentation $\ell$-adique de dimension finie
$\rho$ de $\S_d$, on peut construire un faisceau pervers $\A_\rho$ sur $X_d$
dont les restrictions {\`a} toutes les fibres de $\phi$ sont encore perverses,
{\`a} d{\'e}calage pr{\`e}s (1.3). Si $V_\rho$ est l'espace sous-jacent de 
la repr{\'e}sentation $\rho$,
$\A_\rho$ est un facteur direct du faisceau pervers
$$\RR\pi_*\Ql[\dim(X_d)](\dim(X_d)/2)\boxtimes V_\rho.$$
Cette remarque tr{\`e}s simple se r{\'e}v{\`e}le cruciale 
dans la suite de l'article.

Le but de la section 2 est d'{\'e}tablir une correspondance 
{\`a} la Springer pour $X_d$.
Si cette correspondance est certainement 
connue des sp{\'e}cialistes, nous n'avons
pas pu trouver de d\'emonstration d\'etaill\'ee dans la litt\'erature. 
En 3.1 on d{\'e}montre 
que les complexes $\A_\rho$
sont {\'e}quivariants relativement {\`a} l'action d'un 
sch{\'e}ma en groupes lisse
$G_d\rta Q_d$. Gr{\^a}ce {\`a} cette propri{\'e}t{\'e} 
d'{\'e}quivariance, nous montrons en 3.2 que les $\A_\rho$
sont en quelques sortes ind{\'e}pendants de $n$. En 3.3, 
on {\'e}tablit dans le cas $n=d$,
le lien avec la correspondance de Springer habituelle 
due {\`a} Lusztig, Borho et MacPherson
(\cite{Lus},\cite{B-M}).

A la suite de Mirkovic et Vilonen, pour tous $d',d''\in\NN$ 
tels que $d'+d''=d$, on introduit en 3.1 le produit tordu
$X_{d'}{\tilde\times}X_{d''}\rta X_d$ via lequel est d{\'e}fini 
le complexe produit de convolution
$\A_{\rho'}*\A_{\rho''}$ sur $X_d$. 
Du fait que la r{\'e}solution simultan\'ee se factorise {\`a}
travers $X_{d'}{\tilde\times}X_{d''}$, le complexe $\A_{\rho'}*\A_{\rho''}$  
est un faisceau pervers.
On d{\'e}finit en 3.3  un isomorphisme de commutativit{\'e}
$$\kappa:\A_{\rho'}*\A_{\rho''}{\tilde\rta}\A_{\rho''}*\A_{\rho'}$$
en utilisant le prolongement interm{\'e}diaire 
{\`a} partir de l'ouvert $Q_{d,rss}$.
Il est vraisemblable que $\kappa$ co{\"\i}ncide avec 
l'isomorphisme de commutativit{\'e}
d{\'e}fini de mani{\`e}re diff{\'e}rente par Mirkovic et Vilonen (\cite{M-V}).

La section 4 est consacr{\'e}e {\`a} l'{\'e}tude du complexe 
$\RR\phi_*\A_\rho$.
Le but est de d{\'e}montrer que ce complexe est compl{\`e}tement d{\'e}termin{\'e} par
sa restriction {\`a} n'importe quel ouvert dense. En fait, il satisfait {\`a} un ensemble
de propri{\'e}t{\'e}s plus pr{\'e}cises que nous regroupons 
sous le nom de la propri{\'e}t{\'e} $(*)$.
Un complexe born{\'e} ${\cal C}$ de faisceaux constructibles 
$\ell$-adiques sur un sch{\'e}ma $X$ de type fini
sur $k$ a la propri{\'e}t{\'e} $(*)$ 
si et seulement s'il v{\'e}rifie les propri{\'e}t{\'e}s suivantes.
\begin{itemize}
\item ${\cal C}=\bigoplus_i H^i({\cal C})[-i]$ ;
\item pour toute immersion $j:U\hookrightarrow X$ d'un ouvert dense $U$ de $X$,
        le morphisme d'adjonc\-tion $H^i({\cal C})\rta j_*j^*H^i({\cal C})$
        est un isomorphisme ;
\item Pour tout $x\in X(\FF_{q^r})$, $\Fr^r$ agit dans $H^i({\cal C})$
        comme la multiplication par $q^{ir/2}$.
\end{itemize}

Dans tous nos exemples, les faisceaux de cohomologie $H^i({\cal C})$
sont aussi les faisceaux pervers, et on peut remplacer le $j_*$ par le $j_{!*}$.
De plus, le complexe $\C$ est concentr{\'e} en des degr{\'e}s ayant une parit{\'e} fixe.
Toutefois, nous n'avons pas besoin de ces renseignements suppl{\'e}mentaires. 

On d{\'e}montre en 4.1 que la propri{\'e}t{\'e} $(*)$
se conserve par passage aux facteurs directs,
par image directe par un morphisme fini et par le produit tensoriel externe.
On {\'e}voque aussi le th{\'e}or{\`e}me de Jouanolou 
sur la cohomologie des fibr{\'e}s projectifs
qui v{\'e}rifie la propri{\'e}t{\'e} $(*)$.
On d{\'e}montre en 4.2 le r{\'e}sultat principal de la section 4 :

\setcounter{theoreme}{0}
\begin{theoreme}
 Le complexe $\RR\phi_*\A_\rho$ a la propri{\'e}t{\'e} $(*)$.
\end{theoreme}

On d{\'e}montre d'abord cet {\'e}nonc{\'e} lorsque $\A_\rho=\RR\pi_*\Ql[nd](nd/2)$
en utilisant le th{\'e}or{\`e}me de Jouanolou. 
On peut en d{\'e}duire le cas g{\'e}n{\'e}ral parce
que $\A_\rho$ est un facteur direct de
$\RR\pi_*\Ql[nd](nd/2)\boxtimes V_\rho$.

Comme application du th\'eor\`eme 1, on construit en 4.3 un isomorphisme
$$\displaylines{
        \RR\Gamma_c(X_{d,\vp}\otimes_k{\bar k},\A_{\rho'}*\A_{\rho''}){\tilde\rta}\cr
        \RR\Gamma_c(X_{d',\vp}\otimes_k{\bar k},\A_{\rho'})\otimes
        \RR\Gamma_c(X_{d'',\vp}\otimes_k{\bar k},\A_{\rho''})}$$
en prolongeant l'isomorphisme \'evident sur l'ouvert $Q_{d,rss}$.
Ici $X_{d,\vp}$ d{\'e}signe la fibre de $X_d$ au-dessus de $\vp^d\in Q_d(k)$.
Sur $\CC$, cet isomorphisme a {\'e}t{\'e} obtenu de mani{\`e}re 
diff{\'e}rente par Ginzburg, Mirkovic et Vilonen.

La section 5 est consacr{\'e}e {\`a} l'{\'e}tude cohomologique 
des termes constants.
On note $N$ le sous-groupe de $\GL(n)$ des matrices triangulaires sup{\'e}rieures
unipotentes, $A$ son sous-groupe diagonal, $B=AN$ son sous-groupe de Borel
standard.
Soit ${\cal H}^+$  l'alg{\`e}bre des fonctions {\`a} support
compact dans $\GL(n,F_\vp)^+$ qui sont bi-$\GL(n,{\cal O}_\vp)$-invariantes.
Rappelons que le terme constant de $f\in{\cal H}^+$ 
est une fonction $A({\cal O}_\vp)$-invariante
$f^B:A(F_\vp)\rta\Ql$ d{\'e}finie par
$$f^B(\vp^{\und d})=q^{-\<{\und d},\delta\>}
\int_{N(F_\vp)}f(\vp^{\und d}x)\d x$$
o{\`u}
$$\delta={1\over 2}(n-1,n-3,\ldots,1-n),$$
o{\`u} $\vp^{\und d}=\diag(\vp^{d_1},\ldots,\vp^{d_n})\in A(F_\vp)$
et o{\`u} la mesure de Haar normalis{\'e}e $\d x$ de 
$N(F_\vp)$ attribue {\`a} $N({\cal O}_\vp)$
le volume $1$.
Pour tout $\und d\in\NN^n$ avec $|\und d|=d$,
Mirkovic et Vilonen ont introduit un sous-sch{\'e}ma localement ferm{\'e}
$S_{\und d,\vp}$ de $X_{d,\vp}$ tel que, si $a_\rho$ est la fonction
trace de Frobenius d'un complexe $\A_\rho$ sur $X_{d,\vp}$, on a
$$a_\rho^B(\vp^{\und d})=q^{-\<\und d,\delta\>}
\Tr(\Fr,\RR\Gamma_c(S_{\und d,\vp}\otimes_k{\bar k},\A_\rho)).$$

On introduit en 5.2 une variante globale de $S_{\und d,\vp}$
contenue dans le diagramme commutatif suivant :
$$\diagramme{
        S_{\und d} & \hfl{i_{\und d}}{} & X_d \cr
        \vfl{s_{\und d}}{} &       &\vfl{}{\phi}\cr
        Q_{\und d} & \hfl{}{m} & Q_d }$$
o{\`u} $Q_{\und d}=Q_{d_1}\times\cdots\times Q_{d_n}$ et o{\`u} $m:Q_{\und d}\rta Q_d$
est le morphisme d{\'e}fini par $m(P_1,\ldots,P_n)=P_1\ldots P_n$.
La fibre de $S_{\und d}$ au-dessus de $(\vp^{d_1},\ldots,\vp^{d_n})\in Q_{\und d}(k)$
s'identifie {\`a} $S_{{\und d},\vp}$.
On d{\'e}montre en 5.2 le r{\'e}sultat principal de la section 5.

\begin{theoreme}
Le complexe $\RR s_{\und d,!}i_{\und d}^*\A_\rho$
est concentr{\'e} en degr{\'e} $-d+2\<\und d,\delta\>$ 
et a la propri{\'e}t{\'e} $(*)$.
\end{theoreme}

En particulier 
$$\RR\Gamma_c(S_{\und d,\vp}\otimes_k{\bar k},\A_\rho[-d](-d/2))$$
est concentr{\'e} en degr{\'e} $2\<\und d,\delta\>$. Sur $\CC$, cette assertion
a {\'e}t{\'e} d{\'e}montr{\'e} par Mirkovic et Vilonen (\cite{M-V}).
De plus, l'endomorphisme $\Fr$ agit dans
$$H_c^{2\<\und d,\delta\>}(S_{\und d}
(\vp^{\und d}\otimes_k{\bar k},\A_\rho[-d](-d/2))$$
comme la multiplication par $q^{\<\und d,\delta\>}$. Par cons{\'e}quent, si $a_\rho$
est la fonction trace Frobenius sur $\A_{\rho,\vp}=\A_\rho[-d](-d/2)$ 
restreint {\`a} $X_d(\vp^d)$, sa transformation
de Satake est un polyn{\^o}me sym{\'e}trique 
dont les coefficients sont des entiers naturels.
Cette assertion a {\'e}t{\'e} d{\'e}montr{\'e}e par Lusztig de mani{\`e}re 
diff{\'e}rente (\cite{Lus}).

En utilisant le th{\'e}or{\`e}me 2, on d{\'e}montre en 5.3 l'isomorphisme
$$\RR\Gamma(X_{d,\vp}\otimes_k{\bar k},\A_\rho){\tilde\rta}
\bigoplus_{\scriptstyle\und d\in\NN\atop\scriptstyle |\und d|=d}
\RR\Gamma_c(S_{\und d,\vp}\otimes_k{\bar k},\A_\rho)$$
d\^u {\`a} Mirkovic et Vilonen sur $\CC$ (\cite{M-V}).

On d{\'e}montre enfin en 5.4 qu'il existe un isomorphisme
$$\displaylines{\qquad
\RR\Gamma_c(S_{{\und d},\vp}\otimes_k{\bar k},\A_{\rho'}*\A_{\rho''})\,\,
{\tilde\rta}\bigoplus_{\scriptstyle|{\und d}'|=d',\,
|{\und d}''|=d''\atop\scriptstyle{\und d}'+{\und d}''={\und d}}
\RR\Gamma_c(S_{{\und d}',\vp}\otimes_k{\bar k},\A_{\rho'})\hfill\cr\hfill
\otimes \RR\Gamma_c(S_{{\und d}'',\vp}\otimes_k{\bar k},\A_{\rho''}).\qquad}$$
Si cet isomorphisme n'appara{\^\i}t pas explicitement
dans \cite{M-V}, K. Vilonen m'a inform{\'e} qu'il sait le d{\'e}montrer sur $\CC$.
Cela m'a donn{\'e} la confiance n{\'e}cessaire pour chercher 
{\`a} le d{\'e}montrer sur $\FF_q$.

\bigskip
\bigskip

Dans la seconde partie de l'article on propose certaines applications
de la th\'eorie de Ginzburg, Mirkovic et Vilonen dans le probl\`eme 
dit ``lemme fondamental''.

La section 6 est consacr{\'e}e {\`a} la d{\'e}monstration 
du th{\'e}or{\`e}me 3 {\'e}voqu{\'e} plus haut.
Pour cela, on introduit les ``nouvelles'' fontions 
$f_{r,\lambda}\in{\cal H}^+_r$
d\'efinies comme suit.

On peut identifier $X_\vp(k_r)$ {\`a} l'ensemble
des points fix{\'e}s par l'endomorphisme $\Fr\circ\sigma$
de $X_\vp^r$ 
o{\`u} l'endomorphisme 
$$\sigma:\underbrace{X_\vp\times\cdots\times X_\vp}_{r\ {\rm fois}}\rta 
\underbrace{X_\vp\times\cdots\times X_\vp}_{r\ {\rm fois}}$$ 
est d\'efini par
$$\sigma(x_1,\ldots,x_r)=(x_r,x_1,\ldots,x_{r-1}).$$
En effet on peut faire correspondre un \'el\'ement $x\in X_\vp(k_r)$ \`a
l'\'el\'ement
$$(x,\Fr(x),\ldots,\Fr^{r-1}(x))\in\Fix(\Fr\circ\sigma,X_\vp^r).$$
Pour chaque $\lambda$, l'endomor\-phisme $\sigma$ 
se rel{e}vant naturellement sur le faisceau
pervers $\A_{\lambda,_\vp}^{{\boxtimes\, r}}$ au-dessus $X_\vp^r$,
on d{\'e}finit la fonction $f_{r,\lambda}\in{\cal H}^+_r$
comme la trace de $\Fr\circ\sigma$ agissant sur la fibre de 
$\A_{\lambda,_\vp}^{{\boxtimes\, r}}$
au-dessus des points fix{\'e}s par $\Fr\circ\sigma$.

On d\'emontre ensuite l'identit\'e $b(f_{r,\lambda})=\phi_{r,\lambda}$
qui est l'analogue tordu du fait bien connu  suivant.
Si on remplace $F_r$ par le produit de $r$ copies 
de $F$, l'homomorphisme $b$ envoie la fonction
$f_1\boxtimes\cdots\boxtimes f_r$ sur la fonction $f_1*\cdots *f_r$.
L'isomor\-phisme de la section 5.4 intervient de mani{\`e}re 
cruciale dans la d{\'e}monstration de l'identit\'e 
$b(f_{r,\lambda})=\phi_{r,\lambda}$.

On d\'emontre finalement que $f_{r,\lambda}=a_{r,\lambda}$ en comparant
leurs transform\'es de Satake.

A partir de la section 7,
on ne consid\`ere que les extensions de degr\'e $r=2$. On peut ainsi lib\'erer
la lettre $r$ pour d'autres utilisations. On notera $f_\lambda$ et 
$\phi_\lambda$ pour $f_{2,\lambda}$ et $\phi_{2,\lambda}$.

On pr{\'e}pare dans la section 7 les ingr{\'e}dients n{\'e}cessaires 
pour interpr{\'e}ter g{\'e}om{\'e}triquement l'application 
$b':{\cal H}^+_2\rta{\cal H'}^+$.
Notamment, on veut comprendre le comportement 
de la transposition $\tau(g)=\,^\t g$
vis-{\`a}-vis des faisceaux pervers $\A_\rho$.

La transposition $\tau(g)=\,^\t g$ n'{\'e}tant pas d{\'e}finie sur 
le sch{\'e}ma des r{\'e}seaux $X_d$, on doit introduire en 7.1 
le substitut $\g_{d,r}$.
Il s'agit d'un $Q_{d,r}$-sch{\'e}ma o{\`u} $Q_{d,r}(k)$ est l'ensemble des 
couples de polyn{\^o}mes unitaires $(P,R)$ avec $P$ divisant $R$,
qui sont de degr{\'e} respectivement $d<r$. 
L'ensemble des $k$-points de $\g_{d,r}$
au-dessus de $(P,R)\in Q_{d,r}(k)$ est l'ensemble
$$\g_{d,r}(P,R)(k)=\{g\in\gl(n,\OO/R\OO)\mid
 \det(g)\in P(\OO/R\OO)^\times\}.$$
Les variables auxiliaires $r$ et $R$ sont n{\'e}cessaires pour la finitude.
Le morphisme naturel $\g_{d,r}\rta X_d$ qui envoie $g$ sur $g\OO^n$
est lisse, {\`a} fibres g{\'e}om{\'e}triquement
connexes et {\'e}quivariant relativement {\`a} l'action de $G_{r}$. 
On notera $\tA_\rho$ l'image inverse de $\A_\rho$ sur $\g_{d,r}$.

Puis, on d{\'e}finit le rel{\`e}vement ${\tilde\tau}_\rho$ de $\tau$
sur $\tA_\rho$ en imposant la condition que la restriction de  ${\tilde\tau}_\rho$
{\`a} une fibre de ${\tA}_\rho$ au-dessus d'un {\'e}l{\'e}ment r{\'e}gulier,
semi-simple et sym{\'e}trique soit l'identit{\'e}.
Sa d\'efinition apparemment simple cache certains
aspects assez myst\'erieux de ${\tilde\tau}_\rho$. 
Mais on reportera cette discussion
\`a la derni\`ere section de l'article.

On \'etudie en 7.3 le produit de convolution de $\tA_\rho$ sur $\g_{d,r}$
qui est compl\`etement analogue au produit de convolution entre les fonctions.
En fait, on d\'emontre ainsi que le produit de convolution entre les
$\A_\rho$ d\'efini dans la section 3 correspond bien au produit de 
convolution habituel via le dictionnaire faisceaux-fonctions de Grothendieck.

On \'etudie en 7.4 le comportement des ${\tilde\tau}_\rho$
vis-\`a-vis de l'isomorphisme de commutativit\'e $\kappa$.
Il s'agit en fait de traduire cohomologiquement le calcul 
bien connu suivant prouvant la commutativit\'e de l'alg\`ebre 
de Hecke 
\begin{eqnarray*}
(f*g)(x) &=& \int_{G(F_\vp)}f(xy^{-1})g(y)\d y\cr
         &=& \int_{G(F_\vp)}f(\,^\t y^{-1}\,^\t x)g(\,^\t y)\d y\cr
         &=& \int_{G(F_\vp)}g(\,^\t y)f(\,^\t y^{-1}\,^\t x)\d y\cr
         &=& (g*f)(\,^\t x)\cr
         &=& (g*f)(x).\cr
\end{eqnarray*}

On introduit dans la section 8 une autre r\'ealisation g\'eom\'etrique
des fonctions $a_{2,\lambda}$. On peut identifier $\g_{d,r,\vp}(k_2)$
\`a l'ensemble des points fixes de $\Fr\circ\sigma\circ(\tau\times\tau)$ dans 
$\g_{d,r,\vp}\times\g_{d,r,\vp}$ en envoyant un point
$x\in\g_{d,r,\vp}(k_2)$ sur 
$$(x,\,^\t\Fr(x))\in\Fix(\Fr\circ\sigma\circ(\tau\times\tau),
\g_{d,r,\vp}\times\g_{d,r,\vp}).$$
On d\'efinit la fonction $f'_{\lambda}:\g_{d,r,\vp}\rta\Ql$
par
$$f'_{\lambda}(x)=\Tr(\Fr\circ{\tilde\sigma}\circ
({\tilde\tau}_\lambda\times{\tilde\tau}_\lambda),
(\A_{\lambda,\vp}\boxtimes\A_{\lambda,\vp})_{(x,\,^\t\Fr(x))}).$$

On introduit aussi la fonction 
$\phi'_{\lambda}:\Fix(\Fr\circ\tau,\g_{2d,r,\vp})\rta\Ql$
d\'efinie par
$$\phi'_{\lambda}(x)=\Tr(\Tr\circ\kappa'\circ{\tilde\tau}_\rho,(\tA_\rho)_x)$$
o\`u $\rho$ est la repr\'esentation induite de $\S_{2d}$
$$\rho=\Ind_{\S_d\times\S_d}^{\S_{2d}}(\rho_\lambda\times\rho_\lambda)$$
et o\`u $\kappa'$ se d\'eduit de l'automorphisme de commutativit\'e
de $\A_\rho=\A_\lambda\times\A_\lambda$.
Ces fonctions $\phi'_\lambda$ s'identifient naturellement \`a des 
\'el\'ements de ${\H'}^+$.

On d{\'e}montre en utilisant la formule des traces de Grothendieck,
l'identit\'e $b'(f'_{\lambda})=\phi'_\lambda$.
La r\'esultat de la section 7.4 fournit la compatiblit\'e
des divers endomorphismes de Frobenius tordus.

On d\'emontre ensuite que $f_{\lambda}=f'_{\lambda}$ en tant 
qu'\'el\'ements de $\H^+_2$ si bien qu'on a aussi 
$f'_{\lambda}=a_{2,\lambda}$.

Puisque les $a_{2,\lambda}$ forment une base de ${\cal H}^+_2$, 
pour d\'emontrer l'identit\'e
$$I(a,b(f))=(-1)^{\v_\vp(a_1\ldots a_{n-1})}J(a,b'(f))$$
pour toute $f\in\H^+_2$, il suffit de 
d{\'e}montrer
$$I(a,\phi_\lambda)=(-1)^{\v_\vp(a_1\ldots a_{n-1})}J(a,\phi'_\lambda)$$
pour toute $n$-partition $\lambda$.

Le triplet $({\cal X}_\vp(a),h,\tau)$,
o{\`u} ${\cal X}_\vp(a)$ est un sch{\'e}ma de type fini sur $k$ avec
$${\cal X}_\vp(a)(k)=\{(x,x')\in (N(F_\vp)/N(\OO_\vp))^2
\mid \,^\t x_1ax_2\in\gl(n,\OO_\vp)\},$$
o{\`u}
$h$ est le morphisme ${\cal X}_\vp(a)\rta\Ga$
d\'efini par la formule
$$h(x,x')=\sum_{i=1}^{n-1}\res(x_{i,i+1}+x'_{i,i+1})$$ 
et o{\`u} $\tau$ est une involution de ${\cal X}_\vp(a)$ avec $\tau(x,x')=(x',x)$
a {\'e}t{\'e} d{\'e}fini dans \cite{Ngo1}.

Les complexes ${\tilde\A}_\rho$ ne sont pas 
d{\'e}finis {\`a} priori sur ${\cal X}_{\vp}(a)$.
N{\'e}anmoins, pour un entier $r$ assez grand, on peut construire un diagramme
$$\diagramme{
{\cal X}_{r,\vp}(a)  & \hfl{\iota}{} & \g_{d_n,r,\vp} \cr
\vfl{p_r}{} & &\cr
{\cal X}_{\vp}(a) & & }$$
o{\`u} $\iota$ est une immersion ferm{\'e}e et o{\`u} $p_r$ est morphisme
lisse dont les fibres g{\'e}o\-m{\'e}\-triques sont isomorphes 
{\`a} des espaces affines.
Du fait que les restrictions de ${\tilde\A}_\rho$ dans ces fibres 
g{\'e}om{\'e}triques sont constantes, il peut se descendre en un complexe
${\dot\A}_\rho$ sur ${\cal X}_{\vp}(a)$. L'involution ${\tilde\tau}_\rho$
se descend aussi en une involution 
$${\dot\tau}_\rho:\tau^*{\dot\A}_\rho\rta{\dot\A}_\rho.$$

Gr{\^a}ce {\`a} la formule des traces de Grothendieck, on a alors
$$\displaylines{
I(a,\phi_\lambda)=\Tr(\Fr\circ\kappa,\RR\Gamma_c({\cal X}_\vp(a)\otimes_k{\bar k},
{\dot\A}_\rho\otimes h^*\L_\psi))\cr
J(a,\phi'_\lambda)=\Tr(\Fr\circ\kappa\circ{\dot\tau}_\rho,
\RR\Gamma_c({\cal X}_\vp(a)\otimes_k{\bar k},
{\dot\A}_\rho\otimes h^*\L_\psi))}$$
o{\`u} $\rho$ est la repr{\'e}sentation induite
$$\rho=\Ind_{\S_d\times\S_d}^{\S_{2d}}(\rho_\lambda\times\rho_\lambda)$$
avec $d=|\lambda|=d_n/2$. Rappelons que 
pour cette repr\'esentation $\rho$ on a $\A_\rho=\A_\lambda*\A_\lambda$.

Le th{\'e}or{\`e}me 4 r\'esulte alors de 
l'{\'e}nonc{\'e} g{\'e}om{\'e}trique suivant.

\bigskip

\noindent{\sc Th\'eor\`eme  4a}
{\it L'involution ${\dot\tau}_\rho$ agit dans
$\RR\Gamma_c({\cal X}_\vp(a)\otimes_k{\bar k},
{\dot\A}_\rho\otimes h^*\L_\psi))$
comme la multiplication par
$(-1)^{\v_\vp(a_1a_2\ldots a_{n-1})}$.}

\bigskip

La d{\'e}monstration de ce th{\'e}or{\`e}me occupe la fin de la section 9
et la section 10. Il s'agit d'adapter les arguments de \cite{Ngo1}
dans une situation plus g{\'e}n{\'e}rale. On renvoie {\`a} l'introduction
de loc.cit. pour les grandes lignes de cette d{\'e}monstration.

Dans la derni{\`e}re section 11, on propose un probl{\`e}me combinatoire
concernant un signe $\varepsilon_\lambda$ 
d{\'e}pendant de la partition $\lambda$,
inh{\'e}rent {\`a} 
la d{\'e}finition du rel{\`e}vement ${\tilde\tau}_\lambda$.
On utilise le th\'eor\`eme 4 pour montrer que pour $\lambda=(d,\ldots,d)$
on a
$$\varepsilon_\lambda=(-1)^{d(1+2+\cdots+(n-1))}.$$

En combinant cette identit\'e avec la preuve de la conjecture
de Frenkel-Gaitsgory-Kazhdan-Vilonen (\cite{Ngo3}), 
on d\'emontre enfin le th\'eor\`eme 5.

Je remercie K. Vilonen pour une correspondance tr{\`e}s utile.
Je remercie J.-L. Walds\-purger dont une question 
m'a permis d'am\'eliorer sensiblement l'\'enonc\'e du th\'eor\`eme 3.  
Je tiens {\`a} exprimer ma profonde gratitude envers G. Laumon
qui m'a constamment encourag{\'e} durant la pr{\'e}paration de cet article.

\setcounter{theoreme}{0}
\vfill\eject

\part{La th\'eorie de Ginzburg-Mirkovic-Vilonen}
\section{La situation g{\'e}om{\'e}trique}
\subsection{Le sch{\'e}ma des r{\'e}seaux de Lusztig}
Soient $k$ le corps fini ${\Bbb F}_q$, ${\bar k}$ sa cl{\^o}ture
alg{\'e}brique. Notons ${\cal O}$ 
l'anneau des polyn{\^o}mes {\`a} une variable $\vp$
et {\`a} coefficients dans $k$ et $F$ son corps des fractions.

Fixons un entier naturel $n\in\NN$. L'ensemble des r{\'e}seaux 
$\R\subset{\cal O}^n$
s'identifie naturellement {\`a} celui des $k$-points d'un sch{\'e}ma $X$
localement de type fini sur $k$. En fait, $X$ est la r{\'e}union
disjointe infinie des sch{\'e}mas $X_d$ de type fini sur $k$
dont l'ensemble des $k$-points est
\begin{eqnarray*}
X_d(k)&=&\{\R\subset{\cal O}^n\mid \dim_k\,({\cal O}^n:\R)=d\} \cr
      &=&\{g\,\GL(n,{\cal O})\in\GL(n,F)^+/
        \GL(n,{\cal O})\mid\deg(\det\,g)=d\}
\end{eqnarray*}
o{\`u} le polyn{\^o}me $\det\,g$ est bien d{\'e}fini modulo un scalaire inversible.
Si on note $Q_d$ l'espace affine des polyn{\^o}mes unitaires de degr{\'e} $d$,
le d{\'e}terminant d{\'e}finit un morphisme
$$\phi:X_d\rta Q_d.$$

Pour un polyn{\^o}me unitaire $P\in Q_d(k)$ et pour un r{\'e}seau $\R\in X_d(P)$,
d'apr{\`e}s le th{\'e}or{\`e}me de Cramer, on a ${\cal O}^n\supset\R\supset(P)^n$. 
Le r{\'e}seau $\R$
est donc d{\'e}termin{\'e} par son image dans le quotient $({\cal O}/(P))^n$ si bien que
la fibre $X_d(P)$ s'identifie naturellement au sch{\'e}ma des sous-espaces
vectoriels de codimension $d$ de $({\cal O}/(P))^n$ qui sont stables sous la multiplication
par $\vp$. En particulier le morphisme $\phi$ est propre.

Soient ${\cal O}_\vp$ le compl{\'e}t{\'e} de ${\cal O}$ en $\vp$ et $F_\vp$ son corps des
fractions. L'ensem\-ble des r{\'e}seaux $\R_\vp\subset{\cal O}_\vp^n$ tels que
$\dim_k\,({\cal O}_\vp^n:\R_\vp)=d$ s'identifie naturellement {\`a} celui des
$k$-points de la fibre $X_{d,\vp}$ du morphisme $\phi$ au-dessus du $k$-point
$\vp^d\in Q_d(k)$.

\subsection{Une r{\'e}solution simultan{\'e}e, d'apr{\`e}s Laumon}

Le morphisme $\phi_d$ n'est pas lisse. Ses fibres admettent n{\'e}anmoins
une alt{\'e}ration simultan{\'e}ment semi-petite au sens de
Goresky et MacPherson, compl{\`e}tement analogue {\`a} la
r{\'e}solution simultan{\'e}e de Grothendieck-Springer. 
A la suite de Laumon (\cite{Lau}),
consid{\'e}rons le sch{\'e}ma ${\tilde X}_d$ avec
$${\tilde X}_d(k)=\{{\cal O}^n=\R_0\supset\R_1\supset\cdots\supset\R_d\mid
\dim_k(\R_{i-1}:\R_i)=1\}$$
et le morphisme $\pi_d:{\tilde X}_d\rta X_d$ avec
$$\pi_d\,(\R_0\supset\R_1\supset\cdots\supset\R_d)=\R_d.$$
Soit $\tilde\phi_d:{\tilde X}_d\rta\AA^d$ le morphisme qui envoie
le drapeau $(\R_0\supset\R_1\supset\cdots\supset\R_d)$
sur le point de coordonn{\'e}es $(x_1,\ldots,x_d)$
tel que le quotient $\R_{i-1}/\R_i$ est support{\'e} par le point $x_i$.
On renvoie {\`a} \cite{Lau} pour la d{\'e}monstration de 
l'{\'e}nonc{\'e} suivant ;
voir aussi la fin de 2.3 pour une d{\'e}monstration indirecte via la r{\'e}solution
simultan{\'e}e de Grothendieck-Springer.

\begin{proposition}
Soit $\AA^d_{rss}$ l'ouvert dense de $\AA^d$ constitu{\'e} des points ay\-ant
des coordonn{\'e}es $(x_1,\ldots,x_d)$ deux {\`a} deux diff{\'e}rentes. Cet ouvert {\'e}tant
stable sous l'action du groupe sym{\'e}trique ${\frak S}_d$, posons
$Q_{d,rss}=\AA^d_{rss}/{\frak S}_d$
\begin{enumerate}
\item Au-dessus de l'ouvert dense $Q_{d,rss}$, le diagramme
$$\diagramme{
{\tilde X}_{d} & \hfl{\tilde\phi}{} & \AA^d \cr
\vfl{\pi}{} &  & \vfl{}{\pi_Q} \cr
X_{d} &\hfl{}{\phi}& Q_{d}}$$
est cart{\'e}sien.
\item Le morphisme ${\tilde\phi}_d$ est lisse.
\item Le morphisme $\pi$ est petit au sens de Goresky et MacPherson.
\item Pour tout point g{\'e}om{\'e}trique $x$ de $\AA_d$, la fibre ${\tilde\phi}^{-1}(x)$
est une alt{\'e}ration semi-petite de la fibre $\phi^{-1}(\pi_Q(x))$. Si de plus
$x=(0,\ldots,0)$ et $\pi_Q(x)=\vp^d$, c'est une r{\'e}solution semi-petite.
\end{enumerate}
\end{proposition}

\subsection{Les complexes $\A_\rho$}
Il r{\'e}sulte de la proposition pr{\'e}c{\'e}dente que ${\tilde X}_{d,rss}$ est un
${\frak S}_d$-torseur au-dessus de $X_{d,rss}$.
On peut donc associer {\`a} chaque repr{\'e}sentation $\ell$-adique de dimension finie
$\rho:{\frak S}_d\rta\GL(V)$ une repr{\'e}sentation du groupe fondamental de $X_{d,rss}$
et donc un syst{\`e}me local $\L_\rho$ sur $X_{d,rss}$.
Ce syst{\`e}me local peut aussi {\^e}tre d{\'e}fini par
$\L_\rho=(\pi_{rss,*}\Ql\boxtimes V)^{{\frak S}_d}$, les invariants
{\'e}tant pris relativement {\`a} l'action diagonale de ${\frak S}_d$.
Notons ${\cal A}_\rho$ le faisceau pervers qui est le prolongement interm{\'e}diaire
de $\L_\rho[\dim X_d](\dim X_d/2)$.

\begin{corollaire}
\begin{enumerate}
\item Le complexe 
$$\RR\pi_*\Ql[\dim X_d](\dim X_d/2)$$ 
est un faisceau pervers
prolongement interm{\'e}diaire de sa restriction {\`a} l'ouvert $X_{d,rss}$.
De plus, pour tout point g{\'e}om{\'e}trique
$x\in Q_d({\bar k})$, la restriction de
$$\RR\pi_*\Ql[\dim X_d-d](\dim X_d/2-d/2)$$
{\`a} la fibre $\phi^{-1}(x)$ est encore un faisceau pervers. 
\item On a un isomorphisme
$\A_\rho\simeq (\RR\pi_*\Ql\boxtimes V)^{{\frak S}_d}[\dim X_d](\dim X_d/2)$. 
De plus, pour tout point g{\'e}om{\'e}trique $x\in Q_d({\bar k})$,
la restriction de 
$$\A_\rho[-d](-d/2)$$ 
{\`a} la fibre $\phi^{-1}(x)$
est encore un faisceau pervers. 
\end{enumerate}
\end{corollaire}

\DEMONSTRATION
L'assertion 1 r{\'e}sulte de la proposition 1.2.1. On d{\'e}duit les m{\^e}mes
propri{\'e}t{\'e}s pour $(\RR\pi_*\Ql\boxtimes V)^{{\frak S}_d}[\dim X_d](\dim X_d/2)$ ;
celles-ci {\'e}tant conserv{\'e}es par le produit ext{\'e}rieur $\boxtimes V$
et pour les facteurs directs. Par la d{\'e}finition du syst{\`e}me local $\L_\rho$
on a un isomorphisme
$$\A_\rho\simeq (\RR\pi_*\Ql\boxtimes V)^{{\frak S}_d}
[\dim X_d](\dim X_d/2)$$
au-dessus de l'ouvert $X_{d,rrs}$. On en d{\'e}duit l'isomorphisme sur $X_d$
entier par la fonctorialit{\'e} du prolongement interm{\'e}diaire.$\square$

\section{Une correspondance {\`a} la Springer}

\subsection{La propri{\'e}t{\'e} $\GL(n,{\cal O})$-{\'e}quivariante}
Les fonctions traces de Frobenius des $\A_\rho$ sont des fonctions sur
$$(\GL(n,F)\cap\gl(n,{\cal O}))/\GL(n,{\cal O})$$
qui sont $\GL(n,{\cal O})$-invariantes {\`a} gauche.
Pour exprimer g{\'e}om{\'e}triquement cette propri{\'e}t{\'e}, introduisons le sch{\'e}ma en
groupes $G_d\rta Q_d$ dont la fibre au-dessus d'un point $P\in Q_d(k)$
a pour l'ensemble des $k$-points l'ensemble
$G_d(P)(k)=\GL(n,{\cal O}/(P))$.

\begin{lemme}
Le sch{\'e}ma en groupes $G_d\rta Q_d$ est lisse
{\`a} fibres g{\'e}om{\'e}triques connexes de dimension $n^2d$. 
\end{lemme}

\DEMONSTRATION
Le sch{\'e}ma en groupes $G_d$ est naturellement un ouvert du $Q_d$-sch{\'e}ma $M_d$
d{\'e}fini tel que pour chaque $P\in Q_d(k)$ on a $M_{d,n}(P)(k)=\gl(n,O/(P))$.
Il suffit de d{\'e}montrer l'assertion pour $M_{d,n}$. La dimension relative de $M_{d,n}$
est clairement $n^2d$. Pour d{\'e}montrer que $M_{d,n}$ est lisse, il suffit de le faire
pour $n=1$. Par la division euclidienne, chaque classe modulo $(P)$
est repr{\'e}sent{\'e}e par un unique polyn{\^o}me de degr{\'e} 
strictement inf{\'e}rieur {\`a} $d$
si bien que $M_{d,1}$ est isomorphe au fibr{\'e} vectoriel trivial de rang $d$ sur $Q_d$.
En particulier le morphisme $M_{d,1}\rta Q_d$ est lisse. $\carre$

Un r{\'e}seau $\R\in X_d(k)$ avec $\phi(\R)=P$ 
contient d'apr{\`e}s le th{\'e}or{\`e}me de Cramer
le r{\'e}seau $(P)^n$ et donc est d{\'e}termin{\'e} 
par son quotient $\R/(P)^n\subset {\cal O}^n/(P)^n$.
On en d{\'e}duit une action naturelle
$$G_d\times_{Q_d}X_d\rta X_d.$$

\begin{corollaire}
Les faisceaux pervers $\A_\rho$ sont naturellement $G_d$-{\'e}qui\-va\-riants.
\end{corollaire}
\DEMONSTRATION
$G_d$ agit aussi sur ${\tilde X}_d$. Au-dessus de l'ouvert ${\tilde X}_{d,rss}$
cette action commute {\`a} l'action de ${\frak S}_d$ si bien que le syst{\`e}me local
$\L_\rho$ est naturellement $G_d$-{\'e}quivariant. Consid{\'e}rons maintenant
le diagramme
$$\diagramme{
X_d & \hgfl{act}{} & X_d\times_{Q_d} G_d & \hfl{\pr_{X_d}}{} & X_d \cr
\vfl{}{} &         &\vfl{}{}             &                &\vfl{}{}\cr
Q_d & \hgfl{}{}    & G_d                 & \hfl{}{}          &Q_d}$$
o{\`u} les deux carr{\'e}s sont cart{\'e}siens. Le morphisme
$G_d\rta  Q_d$ {\'e}tant lisse et {\`a} fibres g{\'e}om{\'e}tri\-quement connexes,
il en est de m{\^e}me des morphisme $act$ et $\pr_{X_d}$.
Les images inverse $act^*\A_\rho$ et $\pr_{X_d}^*\A_\rho$
sont donc {\`a} d{\'e}calage pr{\`e}s les prolongements interm{\'e}diaires de
$act^*\L_\rho$ et de $\pr_{X_d}^*\L_\rho$. L'isomorphisme
$act^*\L_\rho\rta\pr_{X_d}^*\L_\rho$ induit donc par fonctorialit{\'e}
un isomorphisme $act^*\A_\rho\rta\pr_{X_d}^*\A_\rho$. $\square$

\subsection{La d{\'e}pendance des $\A_\rho$ en $n$}

Ajoutons dans ce paragraphe l'indice $n$ {\`a} 
toutes nos notations pr{\'e}c{\'e}demment pos{\'e}es.
Pour tout $m\in\NN$, l'homomorphisme $\GL(n)\rta\GL(n+m)$
qui envoie $g$ sur $\diag(g,\Id_m)$ induit une immersion ferm{\'e}e
$X_{d,n}\rta X_{d,n+m}$ au-dessus de $Q_d$.

\begin{lemme}
Le ferm{\'e} $X_{d,n}$ est transverse {\`a} l'action du sch{\'e}ma en groupes
$G_{d,m+n}$ sur $X_{d,m+n}$. Plus pr{\'e}cis{\'e}ment le morphisme
$X_{d,n}\times_{Q_d} G_{d,m+n}\rta X_{d,m+n}$ induit par l'action,
est lisse et ses fibres g{\'e}om{\'e}triques sont connexes. Si de plus $n\geq d$,
il est surjectif. 
\end{lemme}

\DEMONSTRATION
On d{\'e}montre d'abord que l'image g{\'e}om{\'e}trique de ce morphisme est ouvert et que
toutes ses fibres g{\'e}om{\'e}triques sont lisses, connexes et ont la m{\^e}me dimension.

Fixons un polyn{\^o}me unitaire 
$P\in Q_d({\bar k})$ avec
$$P(\vp)=\prod (\vp-\gamma_i)^{m_i}$$
o{\`u} les $\lambda_i\in{\bar k}$ sont deux {\`a} deux diff{\'e}rents.
Les orbites de $G_{d,n}(P)({\bar k})$ (resp.  $G_{d,n+m}(P)({\bar k})$
dans $X_{d,n}(P)({\bar k})$ (resp. $X_{d,n+m}(P)({\bar k})$)
sont param{\'e}tr{\'e}es par la donn{\'e}e d'une $n$-partition $\lambda_{i}=(\lambda_{i,1},\ldots,\lambda_{i,n})$
(resp. une $n+m$-partition $\lambda'_{i}=(\lambda_{i,1},\ldots,\lambda_{i,n+m})$) de $m_i$ pour chaque $i$.
L'ordre partiel d{\'e}fini par la relation d'inclusion d'une orbite dans l'adh{\'e}rence d'une autre
correspond {\`a} l'ordre habituel entre les partitions. 
De plus, l'orbite dans $X_{d,n+m}(P)({\bar k})$ coupe $X_{d,n}(P)({\bar k})$
si et seulement si $\lambda_{i,n+1}=\ldots=\lambda_{i,n+m}=0$ pour tout $i$
et dans ce cas l'intersection est l'orbite correspondant {\`a}
$(\lambda_{i,1},\ldots,\lambda_{i,n})$.

On d{\'e}duit notamment que l'image g{\'e}om{\'e}trique du morphisme
$$X_{d,n}\times_{Q_d} G_{d,n+m}\rta X_{d,n+m}$$
est un ouvert dense  lequel est $X_{d,n+m}$ tout entier si $n\geq d$
et que ses fibres g{\'e}om{\'e}triques sont toutes lisses et connexes.
Pour d{\'e}montrer qu'elles ont la m{\^e}me dimension, il faut d{\'e}montrer que
si
$$\lambda'_i=(\lambda_i,\underbrace{0,\ldots,0}_m)$$
pour tout $i$, la diff{\'e}rence 
entre la dimension de l'orbite dans $X_{d,n+m}(P)({\bar k})$ de param{\`e}tre
$\lambda'_i$ et celle de l'orbite dans  $X_{d,n}(P)({\bar k})$
de param{\`e}tre $\lambda_{i}$ est une constante.

La dimension de l'orbite de param{\`e}tre $\lambda_{i}$ est {\'e}gale {\`a}
$$(n-1)d-2\sum_i\<\lambda_i,(0,1,\ldots,n)\>$$
et celle de l'orbite correspondant {\`a} $\lambda'_i$ est {\'e}gale {\`a}
$$(n+m-1)d-2\sum_i\<\lambda'_i,(0,1,\ldots,n+m)\>$$
si bien que la diff{\'e}rence des deux dimensions est {\'e}gale {\`a}
$md$ du fait que $\lambda_{i,n+1}=\ldots=\lambda_{i,n+m}=0$.

Pour terminer la d{\'e}monstration du lemme, il suffit de 
d{\'e}montrer que la source du morphisme
est lisse. Le sch{\'e}ma en groupes $G_{d,n+m}\rta Q_d$ {\'e}tant lisse,
il suffit de d{\'e}montrer le lemme suivant.

\begin{lemme} Le sch{\'e}ma $X_{d,n}$ est lisse et de dimension $nd$.
\end{lemme} 

\DEMONSTRATION
On va couvrir $X_{d,n}$ par des ouverts lisses.

Choisissons une $n$-partition $\lambda=(\lambda_1,\ldots,\lambda_n)$
de $d$.
Soit $e_1,\ldots,e_n$ la ${\cal O}$-base standard de ${\cal O}^n$.
Notons $V$ le sous-$k$-espace vectoriel de ${\cal O}^n$ d{\'e}fini par
$$V=\bigoplus_{i=1}^n\bigoplus_{j=0}^{\lambda_i-1} \vp^j e_i k.$$
Pour tout polyn{\^o}me unitaire $P$ de degr{\'e} $d$, $V$ s'envoie injectivement
dans le quotient $({\cal O}/(P))^n$. Pour tout r{\'e}seau $\R\in\phi^{-1}(P)$,
l'image de $\R$ dans $({\cal O}/(P))^n$ est de codimension $d$.
G{\'e}n{\'e}riquement ces deux sous-espaces vectoriels de $({\cal O}/(P))^n$
se coupent donc en $0$ ; notons $U$ l'ouvert de $X_{d,n}$ des r{\'e}seaux $\R$ qui coupent
$V$ en $0$. On peut associer {\`a} chaque point $\R\in U$ $n$ vecteurs $v_1,v_2,\ldots,v_n\in V$
tels que pour tout $i$, $\vp^{\lambda_i} e_i+v_i\in\R$.

En fait ces $\vp^{\lambda_i} e_i+v_i$ engendrent $\R$. On a
$$\bigoplus_{i=1}^n (\vp^{\lambda_i} e_i+v_i){\cal O}=g{\cal O}^n$$
o{\`u} pour tout $i$ on a $\deg(g_{ii})=\lambda_i$ et pour tout $j\not=i$ on a
$\deg(g_{j,i})<\lambda_j$ (le polyn{\^o}me nul ayant par convention le degr{\'e} $-\infty$).
Un calcul de d{\'e}terminant {\'e}vident montre que $\deg(\det(g))=\sum_{i=1}^n \lambda_i$
donc $\dim({\cal O}^n/g{\cal O}^n)=d$.
Mais $g{\cal O}^n\subset \R$, on a donc $g{\cal O}^n=\R$.

Le morphisme $U\rta V^n$ est donc un isomorphisme. En particulier $U$ est lisse
est de dimension $nd$. Quitte {\`a} changer la ${\cal O}$-base de ${\cal O}^n$,
on couvre $X_{d,n}$ par des ouverts lisses analogues {\`a} $U$. $\square$

\begin{corollaire}
Pour toute repr{\'e}sentation $\rho$ de $\S_d$, notons $\A_{d,n,\rho}$
le faisceau pervers correspondant sur $X_{d,n}$.
Si $n\leq d$, $\A_{d,n,\rho}$ est canononiquement isomorphe {\`a} la restriction
$\A_{d,d,\rho}[(n-d)d]$ {\`a} $ X_{d,n}$.
Si $n\geq d$, $\A_{d,n,\rho}$ est l'unique faisceau pervers d{\'e}fini {\`a} un unique isomorphisme
pr{\`e}s tel que $\pr^*\A_{d,d\rho}[(n-d)d]\cong \act\*\A_{d,n,\rho}$ o{\`u}
$\pr:X_{d,d}\times_{Q_d}G_{d,n}\rta X_{d,d}$ est la projection naturelle
et $\act:X_{d,d}\times_{Q_d}G_{d,n}\rta X_{d,n}$ est la restriction du morphisme action.
\end{corollaire}

\subsection{Une correspondance {\`a} la Springer pour $X_{d,n}$}
Pr{\'e}cisons davantage la construction pr{\'e}c{\'e}dente dans le cas $d=n$ et $\lambda=(1,\ldots,1)$.
Dans ce cas $V=e_1k\oplus\cdots e_nk$. On a une immersion ouverte $\End(V)\rta X_{n,n}$
d{\'e}finie par $g\mapsto (g+\vp\Id){\cal O}^n$. La d{\'e}monstration du lemme suivant est facile
et sera omise. Notons aussi que ce lemme est une variante d'une compactification
de la vari{\'e}t{\'e} des matrices unipotentes due {\`a} Lusztig (\cite{Lus}).

\begin{lemme}
\begin{enumerate}
\item La restriction de $\phi$ {\`a} l'ouvert $\End(V)$ est l'application
polyn{\^o}me caract{\'e}ristique modulo la convention de signe.
\item Pour tout $P\in Q_d({\bar k})$, toute orbite de $G_{n,n}(P)({\bar k})$
dans $\phi^{-1}(P)$ coupe
$\End(V)$ en une orbite adjointe de $\GL(V)$ dans $\End(V)$.
Cette correspondance est compatible {\`a} la param{\'e}trisation
de ces deux sortes d'orbites par la donn{\'e}e de $n$-partitions $\lambda_i$
de la multiplicit{\'e} $m_i$ de chaque racine $\gamma_i$ de $P$.
\item La restriction de ${\tilde X}_{n,n}\rta X_{n,n}$ {\`a} $\End(V)$ est
la r{\'e}solution de Grothen\-dieck-Springer de $\End(V)$.
\end{enumerate}
\end{lemme}

Rappelons que les repr{\'e}sentations 
irr{\'e}ductibles du groupe sym{\'e}trique $\S_d$
sont param{\'e}\-tris{\'e}es par les partitions 
$\lambda=(\lambda_1\geq \lambda_2\geq\cdots)$ de $d$.
Notons $\lg(\lambda)$ le nombre entier maximal tel que $\lambda_n>0$.

\begin{proposition}
Soient $d,n\in{\Bbb N}$ arbitraires, $\lambda$ une partition de $d$, $\rho_\lambda$
la repr{\'e}sentation irr{\'e}ductible correspondant de $\S_d$. Si $\lg(\lambda)\leq n$
alors la restriction de $\A_{\rho_\lambda}[-d](-d/2)$ {\`a} la fibre 
$X_{d,\vp}=X_d(\vp^d)$
est isomorphe au complexe d'intersection de l'adh{\'e}rence de l'orbite
de $X_{\lambda,\vp}$ param{\'e}tr{\'e}e par la $n$-partition obtenue 
de $\lambda$ en tronquant les $\lambda_i$ avec $i\geq n$ (qui sont nuls).
Si $\lg(\lambda)>n$, cette restriction est nulle.
\end{proposition}

\DEMONSTRATION
En combinant le th{\'e}or{\`e}me de Borho-MacPherson (\cite{B-M}) 
avec le lemme pr{\'e}c{\'e}dent,
on obtient le cas $d=n$. Le cas g{\'e}n{\'e}ral 
s'en d{\'e}duit gr{\^a}ce au corollaire 2.2.3. $\carre$

Notons aussi qu'en combinant les lemmes 2.2.1 et 2.3.1 on retrouve la proposition 1.2.1.

\section{Le produit de convolution}
\subsection{D{\'e}finition}
A la suite de Mirkovic et Vilonen (\cite{M-V}),
consid{\'e}rons le sch{\'e}ma $X_{d'}{\tilde\times}X_{d''}$ dont l'ensemble des $k$-points
est
$$X_{d'}{\tilde\times}X_{d''}(k)=
        \{{\cal O}^n\supset\R'\supset\R\mid\dim({\cal O}^n/\R')=d'
        {\rm\ et\ }\dim(\R'/\R)=d''\}$$
o{\`u} $\R'$ et $\R$ sont des sous-r{\'e}seaux de ${\cal O}^n$. 
Si $g,g'\in\GL(n,F)$ tels
que $\R=g{\cal O}^n$ et $\R'=g'{\cal O}^n$ alors $g'$ et 
$g''={g'}^{-1}g$ appartiennent {\`a}
$\gl(n,{\cal O})$. Les d{\'e}terminants $P'=\det(g')$ 
et $P''=\det(g'')$ sont respectivement
des polyn{\^o}mes de degr{\'e} $d'$ et $d''$ 
ind{\'e}pendants des choix de $g'$ et de $g''$.
On en d{\'e}duit un morphisme 
$X_{d'}{\tilde\times}X_{d''}\rta Q_{d'}\times Q_{d''}$.

\begin{lemme}
Au-dessus de l'ouvert $U\subset Q_{d'}\times Q_{d''}$ avec
$$U(k)=\{(P',P'')\in Q_{d'}\times Q_{d''}(k)\mid \pgcd(P',P'')=1\}$$
on a un isomorphisme
$$(X_{d'}\times X_{d''})\times_{Q_{d'}\times Q_{d''}}U
{\tilde\rta}(X_{d'}{\tilde\times} X_{d''})\times_{Q_{d'}\times Q_{d''}}U.$$
\end{lemme}

\DEMONSTRATION
Soient $\R'\in X_{d'}$ et $\R''\in X_{d''}$ tels que $\phi(\R')=P'$
et $\phi(\R'')=P''$. Sous l'hypoth{\`e}se $\pgcd(P',P'')=1$, si $\R=\R'\cap\R''$
on a $\dim({\cal O}^n/\R)=d$ avec $d=d'+d''$. 
On v{\'e}rifie facilement que le morphisme
$$(X_{d'}\times X_{d''})\times_{Q_{d'}\times Q_{d''}}U
\rta(X_{d'}{\tilde\times} X_{d''})\times_{Q_{d'}\times Q_{d''}}U$$
ainsi d{\'e}fini est un isomorphisme. $\carre$

Notons $(Q_{d'}\times Q_{d''})_{rss}$ l'ouvert des $(P',P'')$ tels que
le polyn{\^o}me $P'P''$ n'a pas de racines multiples, $(X_{d'}\times X_{d''})_{rss}$
son image r{\'e}ciproque. D'apr{\`e}s le lemme, on a deux immersions ouvertes
$$j:(X_{d'}\times X_{d''})_{rss}\hookrightarrow X_{d'}\times X_{d''}$$
et
$${\tilde \jmath}:(X_{d'}\times X_{d''})_{rss}
\hookrightarrow X_{d'}{\tilde\times} X_{d''}.$$
Pour toutes repr{\'e}sentations $\rho'$ 
(resp. $\rho''$) de $\S_{d'}$ (resp. de $\S_{d''}$),
la restriction de $\A_{\rho'}\boxtimes\A_{\rho''}$ 
{\`a} $(X_{d'}\times X_{d''})_{rss}$
est {\`a} d{\'e}calage pr{\`e}s la restriction 
du syst{\`e}me local $\L_{\rho'}\boxtimes\L_{\rho''}$.
Notons $\A_{\rho'}{\tilde\boxtimes}\A_{\rho''}$ 
son prolongement interm{\'e}diaire {\`a}
$X_{d'}{\tilde\times}X_{d''}$ 
$$\A_{\rho'}{\tilde\boxtimes}\A_{\rho''}
={\tilde \jmath}_{!*}j^*\A_{\rho'}{\boxtimes}\A_{\rho''}.$$
  
Notons $\mu:X_{d'}{\tilde\times}X_{d''}\rta X_d$ le morphisme
$$\mu({\cal O}^n\supset\R'\supset\R)=\R.$$
Le produit de convolution est d{\'e}fini par
$$\A_{\rho'}*\A_{\rho''}=\RR\mu_*(\A_{\rho'}{\tilde\boxtimes}\A_{\rho''}).$$

\subsection{La perversit{\'e} du produit de convolution}

Le r{\'e}sultat suivant est du {\`a} Ginzburg, Mirkovic et Vilonen (\cite{Gin},\cite{M-V}).
\begin{proposition}
Le convol{\'e} $\A_{\rho'}*\A_{\rho''}$ est un faisceau pervers prolongement
interm{\'e}diaire de sa restriction {\`a} l'ouvert r{\'e}gulier semi-simple.
Plus pr{\'e}cis{\'e}ment, soit
$\rho=\Ind_{\S_{d'}\times\S_{d''}}^{\S_d}(\rho'\otimes\rho'')$
o{\`u} $d=d'+d''$.
Alors on a un isomorphisme $\A_{\rho'}*\A_{\rho''}{\tilde\rta}\A_{\rho}$.
\end{proposition}

\DEMONSTRATION
Le morphisme $\pi:{\tilde X}_d\rta X_d$ se factorise en
$${\tilde X}_d\  {\buildrel{\pi'}\over\longrightarrow}\  X_{d'}{\tilde\times}X_{d''}
        \ {\buildrel{\mu}\over\longrightarrow}\  X_d.$$
Le morphisme $\pi$ {\'e}tant petit, il en est de m{\^e}me de $\pi'$. Au-dessus de l'ouvert
semi-simple r{\'e}gulier,  $\pi'$ est un torseur en groupe $\S_{d'}\times\S_{d''}$
qui est  par ailleurs isomorphe au torseur ${\tilde X}_{d'}\times{\tilde X}_{d''}$
au-dessus de l'ouvert $(X_{d'}\times X_{d''})_{rss}$. On en d{\'e}duit que
$\A_{\rho'}{\tilde\times}{\A_{\rho''}}$ est un facteur direct du faisceau pervers  
$$\RR\pi'_*\Ql\boxtimes(V_{\rho'}\times V_{\rho''})$$
plus pr{\'e}cis{\'e}ment le facteur direct fixe par rapport {\`a} l'action diagonale de
$\S_{d'}\times\S_{d''}$.
Par cons{\'e}quent $\A_{\rho'}*\A_{\rho''}$ est un facteur direct du faisceau pervers
$$\RR\pi_*\Ql\boxtimes(V_{\rho'}\times V_{\rho''})$$
et est donc un faisceau pervers prolongement interm{\'e}diaire de sa restriction
{\`a} l'ouvert r{\'e}gulier semi-simple.

Il suffit donc de v{\'e}rifier la proposition sur l'ouvert r{\'e}gulier semi-simple.
Via la description d'un syst{\`e}me local comme 
une repr{\'e}sentation du groupe fondamental, 
l'op{\'e}ration image directe par $\mu$ correpond bien {\`a} l'op{\'e}ration
induction de $\S_{d'}\times\S_{d''}$ {\`a} $\S_d$. $\carre$

\subsection{La commutativit{\'e} du produit de convolution}

On a vu qu'au-dessus de l'ouvert o{\`u} le polyn{\^o}me $P=P'P''$
est s{\'e}parable, les sch{\'e}mas
$X_{d'}{\tilde\times}X_{d''}$,$X_{d'}\times X_{d''}$ et
$X_{d''}{\tilde\times}X_{d'}$ sont isomorphes.
Cela induit un isomorphisme  entre 
$\A_{\rho'}*\A_{\rho''}$ et $\A_{\rho''}*\A_{\rho''}$
au-dessus de cet ouvert.
On appellera son prolongement interm{\'e}diaire
$$\kappa:\A_{\rho'}*\A_{\rho''}\,{\tilde\rta}\A_{\rho''}*\A_{\rho''}$$
l'isomorphisme de commutativit{\'e}.
Il est vraisemblable que cet isomorphisme de commutativit{\'e} co{\"\i}ncide avec
celui de Mirkovic et Vilonen (\cite{M-V}).

De mani{\`e}re analogue, on obtient un isomorphisme d'associativit{\'e}
$$(\A_\rho*\A_{\rho'})*\A_{\rho''}{\tilde\rta}\A_\rho*(\A_{\rho'}*\A_{\rho''}).$$

\subsection{Restriction}

Soit $P=\prod_{j=1}^r(\vp-\gamma_i)^{d_j}$ avec 
$\gamma_j\in{\bar k}$ et $\sum_j\gamma_j=d$.
L'{\'e}nonc{\'e} suivant est une variante 
du th{\'e}or{\`e}me 3.3.8 de \cite{Lau}.

\begin{proposition}
Via l'isomorphisme
$$X_{d}(P)=\prod_{j=1}^r X_{d_j}((\vp-\gamma_j)^{d_j})$$
on a un isomorphisme
$$\A_\rho|_{X_{d}(P)}=\bigoplus_i\bigotimes_j \A_{\rho_{i,j}}
        |_{X_{d_j}((\vp-\gamma_j)^{d_j})}$$
o\`u 
$$\Res_{\S_{d_1}\times\cdots\times\S_{d_r}}^{\S_d}\rho
=\bigoplus_i\bigotimes_j\rho_{i,j}$$
chaque $\rho_{i,j}$ \'etant une repr\'esentation irr\'eductible
de $\S_{d_j}$.
\end{proposition}
 
\DEMONSTRATION
Soient $Q_{\und d}=\prod Q_{d_i}$
et $m:Q_{\und d}\rta Q_d$ le morphisme 
$$m(P_1,\ldots,P_r)=P_1\ldots P_r.$$
Soit ${\tilde X}_{\und d}$ le $Q_{\und d}$-sch{\'e}ma dont l'ensemble 
des $k$-points au-dessus de 
$$(P_1,\ldots,P_n)\in Q_{\und d}(k)$$
est l'ensemble des drapeaux
$$\OO^n=\R_0\supset\R_1\supset\cdots\supset\R_r=\R$$
tels que pour tout $i=1,\ldots,r-1$, le quotient $\R_{i-1}/\R_i$
est de dimension $d_i$ comme $k$-espace vectoriel et est annul{\'e} par 
$P_i$ comme $\OO$-module.

La r{\'e}solution $\pi:{\tilde X}_d\rta X_d$ se factorise manifestement
{\`a} travers $\pi_{\und d}:{\tilde X}_{\und d}\rta X_d$.

\begin{lemme}
Soit $U$ l'ouvert dense de $Q_{\und d}$ des suites $(P_1,\ldots,P_r)$
tels que les $P_i$ sont deux {\`a} deux premiers entre eux.
La restriction {\`a} $U$ de l'image r{\'e}ciproque $\pi_{\und d}^*\A_\rho$
est un faisceau pervers prolongement interm{\'e}diaire de sa restriction
{\`a} l'ouvert $U\times_{Q_d}Q_{d,rss}$. 
\end{lemme} 

Admettons provisoirement ce lemme.
Soit $X_{\und d}=\prod_{j=1}^r X_{d_i}$. Au-dessus de $U$, on a un
isomorphisme
$$X_{\und d}\times_{Q_{\und d}}U={\tilde X}_{\und d}\times_{Q_{\und d}}U.$$
En admettant le lemme, la restriction {\`a} $U$ de $\pi_{\und d}^*\A_\rho$
et de $\bigoplus_i\boxtimes_j\A_{\rho_{i,j}}$ sont isomorphes
puisqu'ils sont isomorphes au-dessus de l'ouvert r\'egulier semi-simple.
On en d{\'e}duit que
$$(p_{\und d}^*\A_\rho)_{((\vp-\gamma_j)^{d_j})_{j=1}^r}
=\bigoplus_i\bigotimes_j \A_{\rho_{i,j}}
        |_{X_{d_j}((\vp-\gamma_j)^{d_j})}$$
d'o{\`u} la proposition. $\carre$

{\noindent\it D{\'e}monstration du lemme.} 
Formons le produit cart{\'e}sien 
$${\tilde{\!\tilde X}_{\und d}}={\tilde X}_{\und d}\times_{X_d}{\tilde X}_d.$$
Notons ${\tilde{\tilde\pi}}:{\tilde{\!\tilde X}_{\und d}}\rta{\tilde X}_{\und d}$
le morphisme {\'e}vident.
Au-dessus de l'ouvert $U$, ${\tilde{\!\tilde X}}_{\und d}$ est la r{\'e}union
de $|\S_d|/\prod_{j=1}^r|\S_{d_j}|$ composantes connexes dont chacune
est isomorphe {\`a} ${\tilde X}_d\times_{Q_{\und d}}U$. Par cons{\'e}quent,
le morphisme ${\tilde{\tilde\pi}}$ est petit au-dessus de l'ouvert $U$.
On en d{\'e}duit alors par le th{\'e}or{\`e}me de changement de base 
pour un morphisme propre
que $\pi_{\und d}^*\A_\rho|U$ est isomorphe au faisceau pervers
$(\RR{\tilde{\tilde\pi}}_*\Ql\otimes V_\rho)^{\S_d}$. $\square$

\section{L'{\'e}tude de $\RR\phi_*\A_\rho$}
\subsection{La propri{\'e}t{\'e} $(*)$}

Les complexes qu'on rencontre, comme par exemple $\RR\phi_*\A_\rho$, v{\'e}rifient
un ensemble de propri{\'e}t{\'e}s assez fortes qu'on regroupe sous le nom de
la propri{\'e}t{\'e} $(*)$.

\begin{definition}
Un complexe born{\'e} de faisceaux constructibles $\ell$-adiques 
$\C$ sur un sch{\'e}\-ma $X$ de type fini sur $k$
a la propri{\'e}t{\'e} $(*)$ si et seulement si les 
trois conditions suivantes sont satisfaites
\begin{itemize}
\item
on a un isomorphisme
$$\C=\bigoplus_{i\in\ZZ}H^{i}(\C)[-i];$$
\item
pour toute immersion ouverte $j:U\rta X$ d'un ouvert dense $U$, pour tout entier $i$,
le morphisme naturel $H^{i}(\C)\rta j_*j^*H^{i}(\C)$ est un isomorphisme ;
\item
pour tout $x\in X(\FF_{q^r})$, l'endomorphisme $\Fr^r$ agit dans $H^{i}({\cal C})_x$
comme la multiplication par $q^{ri/2}$.
\end{itemize}
Le complexe $\C$ a la propri{\'e}t{\'e} $(*)$ renforc{\'e}e si de plus les faisceaux 
$H^{i}(\C)$ sont des faisceaux constants.
\end{definition}

Dans nos exemples, les faisceaux de cohomologie $H^i({\cal C})$ sont
en fait des faisceaux pervers et on peut remplacer dans la deuxi{\`e}me condition
le $j_*$ par le $j_{!*}$.
De plus, $\C$ est concentr{\'e} en des degr{\'e} ayant une parit{\'e} fixe.
Nous n'utilisons toutefois pas ces renseignements suppl{\'e}mentaires.

\begin{lemme}
\begin{enumerate}
\item Si ${\cal C}$ a la propri{\'e}t{\'e} $(*)$ alors pour tout $d\in\ZZ$,
${\cal C}[d](d/2)$ l'a aussi.
\item Tout facteur direct d'un complexe $\C$ ayant la propri{\'e}t{\'e} $(*)$ l'a aussi.
\item Soit $f:X\rta Y$ un morphisme fini. Si un complexe de faisceaux $\C$ sur $X$
a la propri{\'e}t{\'e} $(*)$ alors le complexe $\RR f_*\C$ l'a aussi.
\item Soient $\C$ et $\C'$ deux complexes de faisceaux respectivement sur $X$ et $X'$
ayant tous les deux la propri{\'e}t{\'e} $(*)$, 
alors le produit tensoriel externe $\C\boxtimes\C'$
l'a aussi.
\item Soient $f:X\rta Y$ et $f':Y\rta Z$ deux morphismes tels que $\RR f_*\Ql$
et $\RR f'_*\Ql$ ont tous les deux la propri{\'e}t{\'e} $(*)$ renforc{\'e}e. 
Alors $\RR(f\circ f')_*\Ql$ l'a aussi.
\end{enumerate}
\end{lemme}

\DEMONSTRATION
La propri{\'e}t{\'e} sur l'action ponctuelle de Frobenius sur les faisceaux de cohomologie
se conserve de mani{\`e}re {\'e}vidente par rapport aux op{\'e}rations passage aux facteurs directs,
l'image directe par un morphisme fini et produit tensoriel externe. Examinons
seulement la conservation des deux autres propri{\'e}t{\'e}s. 

\begin{enumerate}
\item C'est {\'e}vident.
\item Supposons que ${\cal C}={\cal C}'\oplus{\cal C}''$. Alors on a
$H^i({\cal C})=H^i({\cal C}')\oplus H^i({\cal C}'')$ pour tout $i$.
Le morphisme compos{\'e}
$$\bigoplus_{i}H^i({\cal C}')
\rta\bigoplus_{i}H^i({\cal C}){\tilde\rta}{\cal C }\rta{\cal C}'$$
induit un isomorphisme sur les groupes de cohomologie et est donc un
isomorphisme dans la cat{\'e}gorie d{\'e}riv{\'e}e.

Il reste {\`a} d{\'e}montrer que si $j:U\hookrightarrow X$ une immersion ouverte,
pour tout entier $i$, le morphisme d'adjonction
$H^i({\cal C'})\rta j_*j^*H^i({\cal C'})$ est un isomorphisme.
Pour tout faisceau ${\cal F}$ sur $X$, on a la suite exacte habituelle
$$0\rta i_*i^!{\cal F}\rta F\rta j_*j^*{\cal F}\rta 0$$
o{\`u} $i$ est l'immersion du ferm{\'e} compl{\'e}mentaire de $U$ dans $X$.
L'assertion se d{\'e}duit de ce que
$$i_*i^!{\cal C}=i_*i^!{\cal C}'\oplus i_*i^!{\cal C}''.$$

\item Le morphisme $f$ {\'e}tant fini, on a
$$\RR f_*({\cal C}){\tilde\rta}\bigoplus_{i}H^i(\RR f_*{\cal C})$$
avec $H^i(\RR f_*{\cal C})=f_*H^i({\cal C})$.

Soient $j$ l'immersion d'un ouvert dense $U$ dans $Y$, $j_X:U_X\rta X$
son image r{\'e}ciproque. Par le th{\'e}or{\`e}me de changement de base, on a
un isomorphisme 
$$j^*f_* H^i({\cal C}){\tilde\rta} f_{U,*}j_X^*H^i({\cal C}).$$
En composant les isomorphismes suivants
\begin{eqnarray*}
j_*j^*f_*H^i({\cal C}) & {\tilde\rta} & j_*f_{U,*}j_X^*H^i({\cal C}) \cr
                       & {\tilde\rta} & f_*j_{X,*}j_X^*H^i({\cal C}) \cr
                       & {\tilde\rta} & f_*H^i({\cal C})
\end{eqnarray*}
on obtient l'isomorphisme recherch{\'e}.

\item Les op{\'e}rations $j_*$ et $j^*$ se comportent bien par rapport au produit
tensoriel externe. 

\item On utilise la formule de projection.
\end{enumerate}

\begin{lemme}(Jouanolou)
\begin{enumerate}
\item Soit $f:X\rta S$ un fibr{\'e} projectif sur une base $S$
 de type fini sur $k$.
Alors le complexe $\RR f_*\Ql$ a la propri{\'e}t{\'e} $(*)$ renforc{\'e}e.
\item Soit $f:X\rta S$ un fibr{\'e} vectoriel de rang $r$ sur une base $S$ de type fini sur $k$.
Alors le complexe $\RR f_!\Ql$ a la propri{\'e}t{\'e} $(*)$ renforc{\'e}e. Plus pr{\'e}cis{\'e}ment,
la trace induit un isomorphisme $\RR f_!\Ql\rta \Ql[-2r](-r)$.
\end{enumerate}
\end{lemme}

On renvoie {\`a} \cite{Jou} pour la d{\'e}monstration de ce lemme.

La propri{\'e}t{\'e} $(*)$ nous sera utile gr{\^a}ce {\`a} l'{\'e}nonc{\'e} suivant.
Sa d{\'e}monstration est facile et sera omise.

\begin{lemme}
Soient ${\cal C}$ et ${\cal C}'$ deux complexes de faisceaux sur $X$ ayant la propri{\'e}t{\'e} $(*)$.
Tout isomorphisme entre ${\cal C}$ et ${\cal C}'$ au-dessus d'un ouvert dense $U$ de $X$
se prolonge en un isomorphisme sur $X$ tout entier.
\end{lemme} 

\subsection{La propri{\'e}t{\'e} $(*)$ du complexe $\RR\phi_*\A_\rho$}

\begin{theoreme}
Le complexe $\RR\phi_*\A_\rho$ a la propri{\'e}t{\'e} $(*)$.
\end{theoreme}

\DEMONSTRATION
Rappelons qu'on a une r{\'e}solution simultan{\'e}e $\pi:{\tilde X}_d\rta X_d$.
Il suffit de d{\'e}montrer l'assertion en rempla{\c c}ant $\A_\rho$ par $\RR\pi_*\Ql$ 
car $\A_\rho$ est un facteur direct du faisceau pervers
$$\RR\pi_*\Ql[nd](nd/2)\boxtimes V_\rho$$
o{\`u} $V_\rho$ est l'espace sous-jacent de $\rho$.
Du fait que le diagramme
$$\diagramme{
{\tilde X}_{d} & \hfl{\tilde\phi}{} & \AA^d \cr
\vfl{\pi}{} &  & \vfl{}{\pi_Q} \cr
X_{d} &\hfl{}{\phi}& Q_{d}}$$
est commutatif, on a
$$\RR\phi_*\RR\pi_*\Ql{\tilde\rta}\RR\pi_{Q,*}\RR{\tilde\phi}_*\Ql.$$
D{\'e}montrons d'abord un lemme.

\begin{lemme}
Le morphisme $\tilde\phi$ est le compos{\'e} de $d$ fibr{\'e}s projectifs de rang $n-1$.
\end{lemme}

\DEMONSTRATION
Notons ${\tilde X}_{d,i}$ le sch{\'e}ma au-dessus de $\AA^d$ dont la fibre au-dessus de
$$(x_1,\ldots,x_d)\in\AA^d(k)$$ 
est l'ensemble des filtrations de r{\'e}seaux
$${\cal O}^n=\R_0\supset\R_1\supset\cdots\supset\R_i$$
tel que pour tout $j\leq i$, le quotient $\R_{i-1}/\R_i$ est un ${\cal O}$-module
de longueur $1$, support{\'e} par $x_i$.

D{\'e}montrons que ${\tilde X}_{d,i+1}\rta {\tilde X}_{d,i}$ est un fibr{\'e} projectif
de rang $n-1$. En effet, on a
$$\R_i\supset\R_{i+1}\supset(\vp-x_{i+1})\R_i.$$
Le r{\'e}seau $\R_i$ et $x_{i+1}$ {\'e}tant fix{\'e}s, la donn{\'e}e de $\R_{i+1}$ est {\'e}quivalente {\`a} celle
d'un sous-espace vectoriel de codimension $1$ de
$$\R_i/(\vp-x_{i+1})\R_i$$
d'o{\`u} r{\'e}sulte le lemme. $\carre$

{\noindent \it Fin de la d{\'e}monstration.}
Le morphisme ${\tilde\phi}$ {\'e}tant un compos{\'e} de $d$ fibr{\'e}s projectifs de rang $n-1$,
le complexe $\RR{\tilde\phi}_*\Ql$ a la propri{\'e}t{\'e} $(*)$ renforc{\'e}e.
Il s'ensuit que le complexe $\RR\pi_{Q,*}\RR{\tilde\phi}_*\Ql$ a la propri{\'e}t{\'e} $(*)$
du fait que $\pi_Q$ est un morphisme fini. $\square$

\subsection{Cohomologie du produit de convolution}
Soient $\rho'$, $\rho''$ des repr{\'e}sentations de $\S_{d'}$ et de $\S_{d''}$,
$\A_{\rho'}$ et $\A_{\rho''}$ les faisceaux pervers correspondants.
Notons $m_Q:Q_{d'}\times Q_{d''}\rta Q_d$ le morphisme multiplication,
$\phi_d'$ le morphisme $X_{d'}\rta Q_{d'}$ d{\'e}fini en 1.1, $\phi_{d''}$ et $\phi_d$
les morphismes analogues.

\begin{proposition}
On a un isomorphisme
$$m_{Q,*}(\RR\phi_{d',*}\A_{\rho'}\boxtimes \RR\phi_{d'',*}\A_{\rho''}){\tilde\rta}
        \RR\phi_{d,*}(\A_{\rho'}*\A_{\rho''}).$$
\end{proposition}

\DEMONSTRATION
Rappelons qu'on a un diagramme commutatif
$$\diagramme{X_{d'}{\tilde\times}X_{d''} & \hfl{\mu}{} & X_d\cr
              \vfl{\phi_{d'}{\tilde\times}\phi_{d''}}{} & &\vfl{}{\phi_{d}}\cr
                Q_{d'}\times Q_{d''} &\hfl{}{m_Q} & Q_d}$$
qui induit un isomorphisme
$$m_{Q,*}\RR(\phi_{d'}{\tilde\times}\phi_{d''})_*(\A_{\rho'}{\tilde\boxtimes}\A_{\rho''})
        {\tilde\rta}\RR\phi_{d,*}(\A_{\rho'}*\A_{\rho''}).$$
Rappelons qu'au-dessus de l'ouvert r{\'e}gulier semi-simple qu'on a un isomorphisme
$$(X_{d'}\times X_{d''})_{rss}{\tilde\rta}(X_{d'}{\tilde\times}X_{d''})_{rss}$$
de $(Q_{d'}\times Q_{d''})_{rss}$-sch{\'e}mas. Notons ${\tilde\jmath}$ et $j$
les immersions ouvertes
$$\displaylines{{\tilde\jmath}:(Q_{d'}\times Q_{d''})_{rss}\hookrightarrow (Q_{d'}\times Q_{d''})\cr
                j:Q_{d,rrs}\hookrightarrow Q_d.}$$
On a alors un isomorphisme
$${\tilde\jmath}^*(\RR\phi_{d',*}\A_{\rho'}\boxtimes \RR\phi_{d'',*}\A_{\rho''}){\tilde\rta}
{\tilde\jmath}^*\RR(\phi_{d'}{\tilde\times}\phi_{d''})_*(\A_{\rho'}{\tilde\boxtimes}\A_{\rho''})$$
qui induit un isomorphisme
$$j^*m_{Q,*}(\RR\phi_{d',*}\A_{\rho'}\boxtimes \RR\phi_{d'',*}\A_{\rho''}){\tilde\rta}
j^*m_{Q,*}\RR(\phi_{d'}{\tilde\times}\phi_{d''})_*(\A_{\rho'}{\tilde\boxtimes}\A_{\rho''})$$
donc un isomorphisme
$$j^*m_{Q,*}(\RR\phi_{d',*}\A_{\rho'}\boxtimes \RR\phi_{d'',*}\A_{\rho''}){\tilde\rta}
j^*\RR\phi_{d,*}(\A_{\rho'}*\A_{\rho''}).$$
Or, d'apr{\`e}s 3.4, $m_{Q,*}(\RR\phi_{d',*}\A_{\rho'}\boxtimes \RR\phi_{d'',*}\A_{\rho''})$
et $\RR\phi_{d,*}(\A_{\rho'}*\A_{\rho''})$ ont tous les deux la propri{\'e}t{\'e} $(*)$
d'o{\`u} se d{\'e}duit la proposition. $\square$

\begin{corollaire}
\begin{enumerate}
\item Pour toute repr{\'e}sentation $\rho$ de $\S_d$, on a un isomorphisme
$$\RR\Gamma(X_{d,\vp}\otimes_{k}{\bar k},\A_\rho)
=\bigoplus_{i}H^{i}(X_{d,\vp}\otimes_{k}{\bar k},\A_\rho).$$
Ici $X_{d,\vp}$ d{\'e}signe la fibre de $X_d$ au-dessus de $\vp^d\in Q_d(k)$.
\item Pour toutes repr{\'e}sentations $\rho'$ et $\rho''$ de $\S_d'$ et $\S_{d''}$,
pour $d=d'+d''$ et pour un entier $i$ arbitraire, on a
$$\displaylines{H^{i}(X_{d,\vp}\otimes_{k}{\bar k},\A_{\rho'}*\A_{\rho''})\cr
=\bigoplus_{i'+i''=i}
H^{i'}(X_{d',\vp}\otimes_k{\bar k},\A_{\rho'})\otimes
H^{i''}(X_{d'',\vp}\otimes_k{\bar k},\A_{\rho''}).}$$ 
\end{enumerate}
\end{corollaire}

Sur $\CC$, ce corollaire est du {\`a} Ginzburg, Mirkovic et Vilonen (\cite{Gin},\cite{M-V}).

\section{Termes constants}

\subsection{L'isomorphisme de Satake}

Soient $N$ le sous-groupe des matrices triangulaires sup{\'e}rieures unipotentes de $\GL(n)$,
$A$ son tore diagonal et $B=AN$ le sous-groupe de Borel standard.
Soit $\H^+$  l'espace vectoriel des fonctions {\`a} valeurs dans $\Ql$ et
{\`a} support compact dans
$$\GL(n,F_\vp)^+=\GL(n,F_\vp)\cap\gl(n,{\cal O}_\vp)$$ 
qui sont bi-$\GL(n,{\cal O}_\vp)$-invariantes.
Le produit de convolution munit {\`a} $\H^+$ d'une structure d'alg{\`e}bre.
A la suite de Satake (\cite{Sa}) on d{\'e}finit l'application 
$\Phi:\H^+\rta\Ql[Z_1,\ldots,z_n]$ par
$$\Phi(f)=\sum_{{\und d}\in\NN^n} f^B(\vp^{\und d})z^{\und d}$$
o{\`u}
$\vp^{{\und d}}=\diag(\vp^{{d}_1},\ldots,\vp^{{d}_n})\in A(F_\vp)$
et o{\`u} $z^{\und d}=z_1^{d_1}\ldots z_n^{d_n}$. Les termes constants $f^B(\vp^{\und d})$
sont d{\'e}finis par
$$f^B(\vp^{\und d})=q^{-\<\delta,{\und d}\>}
\int_{N(F_\vp)}f(\vp^{\und d} x)\d x,$$
o{\`u} la mesure de Haar $\d x$ normalis{\'e}e attribue {\`a} $N({\cal O}_\vp)$ le volume $1$
et o{\`u}
$$\delta={1\over 2}\<n-1,n-3,\ldots,1-n\>.$$
D'apr{\`e}s Satake, $\Phi$ d{\'e}finit un isomorphisme d'alg{\`e}bres entre $\H^+$
et la sous-alg{\`e}bre de $\Ql[z_1,\ldots,z_n]$ des polyn{\^o}mes sym{\'e}triques.

Pour chaque ${\und d}\in\NN^n$ tel que
$$|\und d|=\d_1+\cdots+d_n=d,$$
{\`a} la suite de Mirkovic et Vilonen (\cite{M-V})
consid{\'e}rons le sous-sch{\'e}ma localement ferm{\'e} $S_{{\und d},\vp}$ de $X_{d,\vp}$
dont l'ensemble des $k$-points est celui des r{\'e}seaux $\R_\vp\subset{\cal O}_\vp^n$
tels que pour tout $i$ on a
$$(\R_\vp\cap \bigoplus_{j=1}^i {\cal O}_\vp e_j)/(\R_\vp\cap \bigoplus_{j=1}^{i-1} {\cal O}_\vp e_j)
=(\bigoplus_{j=1}^{i-1}{\cal O}_\vp e_j\oplus{\cal O}_\vp\vp^{d_i}e_i)
/(\bigoplus_{j=1}^{i-1}{\cal O}_\vp e_j)$$
o{\`u} $e_1,\ldots,e_n$ est la base standard de ${\cal O}_\vp^n$.

\begin{proposition}
Soient $\A_{\rho,\vp}$ un complexe de faisceaux $\ell$-adiques sur $X_{d,\vp}$ qui est
$G_d(\vp^d)$-{\'e}quivariant, $a_\rho$ sa fonction trace de Frobenius qui est
de mani{\`e}re naturelle un {\'e}l{\'e}ment de $\H^+$. Alors, on a
$$a_{\rho}^B(\vp^{\und d})
=q^{-\<{\und d},\delta\>}\Tr(\Fr,
\RR\Gamma_c(S_{{\und d},\vp}\otimes_k{\bar k},\A_{\rho,\vp})).$$
\end{proposition}

\DEMONSTRATION L'ensemble
$$\{x\in N(F_\vp)/N({\cal O}_\vp)\mid\vp^{\und d}x\in\gl(n,{\cal O}_\vp)\}$$
s'identifiant naturellement {\`a} celui des $k$-points de $S(\vp^{\und d})$,
la proposition est une cons{\'e}quence de la formule des traces de
Grothendieck-Lefschetz. $\carre$

Lorsque $\A_{\rho,\vp}$ est la restriction de $\A_\rho[-d](-d/2)$ {\`a} $X_{d,\vp}$,
ce complexe 
$$\RR\Gamma_c(S_{{\und d},\vp}\otimes_k{\bar k},\A_\vp))$$ 
est en fait concentr{\'e} en un seul degr{\'e}. 
La dimension du groupe de cohomologie non trivial peut {\^e}tre 
calcul{\'e}e {\`a} partir de la donn{\'e}e combinatoire associ{\'e}e
{\`a} $\rho$. De plus, l'endomorphisme de Frobenius agit dans ce groupe
comme la multiplication par une puissance de $q$. 
Toutefois, il semble difficile d'obtenir directement 
ces renseignements sans {\'e}tudier au pr{\'e}alable la situation
globale.

\subsection{Termes constants globaux}
Notons $Q_{\und d}=\prod_{i=1}^n Q_{d_i}$ et $m_{\und d}: Q_{\und d}\rta Q_d$ le morphisme
produit
$$m_{\und d}(P_1,\ldots,P_n)=P_1\ldots P_n.$$
Notons $S_{\und d}$ le sous-sch{\'e}ma localement ferm{\'e} de $Q_{\und d}\times_{Q_d} X_d$
dont l'ensemble des $k$-points au-dessus de $(P_1,\ldots,P_n)\in Q_{\und d}(k)$
est celui des r{\'e}seaux $\R\subset{\cal O}^n$ tels que pour tout $i$ on a
$$(\R\cap \bigoplus_{j=1}^i {\cal O} e_j)/(\R\cap \bigoplus_{j=1}^{i-1} {\cal O} e_j)
=(\bigoplus_{j=1}^{i-1}{\cal O} e_j\oplus{\cal O} P_ie_i)
/(\bigoplus_{j=1}^{i-1}{\cal O} e_j).$$
Notons $s_{\und d}:S_{\und d}\rta Q_{\und d}$ le morphisme structurel
et $i_{\und d}:S_{\und d}\rta X_d$ le morphisme {\'e}vident.

\begin{theoreme}
Pour toute repr{\'e}sentation $\rho$ du groupe $\S_d$, le complexe
$$\RR s_{{\und d},!}i_{\und d}^*\A_\rho$$ est concentr{\'e} en degr{\'e}
$-d+2\<{\und d},\delta\>$ et a la propri{\'e}t{\'e} $(*)$.
\end{theoreme}

\DEMONSTRATION
Rappelons qu'on a une r{\'e}solution simultan{\'e}e $\pi:{\tilde X}_d\rta X_d$.
Il suffit de d{\'e}montrer l'assertion en rempla{\c c}ant  $\A_\rho$ par $\RR\pi_*\Ql[nd](nd/2)$ 
car $\A_\rho$ est un facteur direct du faisceau pervers
$$\RR\pi_*\Ql[nd](nd/2)\boxtimes V_\rho$$
o{\`u} $V_\rho$ est l'espace sous-jacent de $\rho$.

Posons ${\tilde S}_{\und d}=S_{\und d}\times_{X_d}{\tilde X}_d$ et d{\'e}signons aussi
par $\pi$ le morphisme ${\tilde S}_{\und d}\rta S_{\und d}$
et par ${\tilde s}_{\und d}$ le compos{\'e} $s_{\und d}\circ\pi$.
On va d{\'e}montrer que le complexe
$$\RR{\tilde s}_{{\und d},!}\Ql=\RR s_{{\und d},!}\RR\pi_*\Ql$$
est concentr{\'e} en degr{\'e} $2\sum_{i=1}^n d_i(n-i)$ et a la propri{\'e}t{\'e} $(*)$.

Soit $(P_1,\ldots,P_n)\in Q_{\und d}(k)$. L'ensemble des $k$-points dans la fibre
de ${\tilde S}_{\und d}$ au-dessus de $(P_1,\ldots,P_d)$ est alors celui des filtrations
compl{\`e}tes
$${\cal O}^n=\R_0\supset\R_1\supset\cdots\supset\R_d=\R$$
d'un r{\'e}seau $\R$ avec $\R\in S_{\und d}(P_1,\ldots,P_n)(k)$.

Pour tout $l$,
le sous-quotient $\R_{l-1}/\R_l$, {\'e}tant de dimension $1$, intervient dans un unique
sous-quotient
$$\bigoplus_{j=1}^i {\cal O} e_j/\bigoplus_{j=1}^{i-1}{\cal O} e_j.$$
On en d{\'e}duit une fonction
$$\tau:\{1,\ldots,d\}\rta\{1,\ldots,n\}$$
telle que pour tout $i\in\{1,\ldots,n\}$, le cardinal de $\tau^{-1}(i)$ est
{\'e}gal {\`a} $d_i$.

Notons ${\tilde S}_{\und d,\tau}$ le sous-sch{\'e}ma de ${\tilde S}_{\und d}$ constitu{\'e}
des points qui induisent la fonction $\tau$. Pour tout
$(x_1,\ldots,x_d)\in \AA^d(k)$, l'ensemble des $k$-points de ${\tilde S}_{{\und d},\tau}$
au-dessus de $(x_1,\ldots,x_d)$ est celui des filtrations de r{\'e}seaux
$${\cal O}^n=\R_0\supset\R_1\supset\cdots\supset\R_d=\R$$
telles que pour tout $l$, le quotient $\R_{l-1}/\R_l$ est support{\'e} par $x_l$
et que l'entier minimal $i$ tel que
$$\R_{l-1}\cap\bigoplus_{j=1}^i{\cal O} e_j\varsupsetneq
\R_l\cap\bigoplus_{j=1}^i{\cal O} e_j$$
est {\'e}gal {\`a} $\tau(l)$.
La restriction du morphisme 
${\tilde s}_{\und d}:{\tilde S}_{{\und d}}\rta Q_{\und d}$ {\`a} ${\tilde S}_{{\und d},\tau}$
se factorise comme suit
$${\tilde S}_{{\und d},\tau}
{\buildrel{{\tilde s}_{{\und d},\tau}}\over\longrightarrow}\AA^d
{\buildrel\pi_\tau\over\longrightarrow}Q_{\und d}$$
o{\`u}
$$\pi_\tau(x_1,\ldots,x_d)=(Q_1,\ldots,Q_n)$$
avec $Q_i=\prod_{l\in\tau^{-1}(i)}(\vp-x_l)$.

\begin{lemme}
Le morphisme
${\tilde S}_{\und d,\tau}{\buildrel{{\tilde s}_{{\und d},\tau}}\over\longrightarrow} \AA^d$
est le compos{\'e} de $d$ fibr{\'e}s vectoriels successivement de rang
$$n-\tau(n),n-\tau(n-1),\ldots,n-\tau(1).$$
En particulier ${\tilde s}_{{\und d},\tau}$ est lisse et {\`a} fibres
g{\'e}om{\'e}triques connexes et de dimension
$$\sum_{l=1}^d\tau(d)=\sum_{i=1}^n d_i(n-i).$$
\end{lemme}

{\noindent\it D{\'e}monstration du lemme.}
Notons ${\tilde S}_{\und d,\tau,l}$ le sch{\'e}ma au-dessus de $\AA^d$ dont l'ensemble
des $k$-points au-dessus de $(x_1,\ldots,x_d)\in\AA^d(k)$ est celui des filtrations
$${\cal O}^n=\R_0\supset\R_1\supset\cdots\supset\R_l$$
telles que pour tout $l'\leq l$, le sous-quotient $\R_{l'-1}/\R_{l'}$ est support{\'e} par $x_{l'}$
et intervient dans 
$\bigoplus_{j=1}^{\tau(l')}{\cal O} e_j/\bigoplus_{j=1}^{\tau(l')-1}{\cal O} e_j$.

Fixons un $k$-point de ${\tilde S}_{\und d,\tau,l-1}$
$$(\R_0\supset\R_1\supset\cdots\supset\R_{l-1} ; x_1,\ldots,x_d)$$
Au-dessus de ce point, la donn{\'e}e d'un $k$-point de ${\tilde S}_{\und d,\tau,l}$
est la donn{\'e}e d'un r{\'e}seau $\R_l$ tel que
le quotient $\R_{l-1}/\R_{l}$ est support{\'e} par $x_{l}$
et que l'entier minimal $i$ tel que
$$\R_{l-1}\cap\bigoplus_{j=1}^i{\cal O} e_j\varsupsetneq
\R_{l}\cap\bigoplus_{j=1}^i{\cal O} e_j$$
est {\'e}gal {\`a} $\tau(l)$.
En particulier, on a
$$R_{l-1}\supset\R_{l}\supset (\vp-x_{l})\R_{l-1}.$$

Notons $V=\R_{l-1}/(\vp-x_{l})\R_{l-1}$.
La filtration
$$0\subset{\cal O} e_1\subset{\cal O} e_1\oplus{\cal O} e_2\subset\cdots\subset\bigoplus_{i=1}^n{\cal O} e_i$$
induit une filtration compl{\`e}te $V_\bullet$
$$0\subset V_1\subset\cdots\subset V_n=V.$$

La donn{\'e}e de $\R_{l}$ est {\'e}quivalente {\`a} la donn{\'e}e d'un sous-espace vectoriel
de codimension $1$ de $V$ qui contient $V_{\tau(l)-1}$ mais ne contient pas
$V_{\tau(l)}$. L'ensemble de ces donn{\'e}e constitue visiblement celui des $k$-points d'un
espace affine de rang $n-\tau(l)$.

La donn{\'e}e de $V$ munie de la filtration $V_\bullet$ d{\'e}finit un fibr{\'e} vectoriel
filtr{\'e} sur ${\tilde S}_{\und d,\tau,l-1}$. Nous aurions termin{\'e} la d{\'e}monstration
de ce lemme si nous savions que ce fibr{\'e} vectoriel filtr{\'e} est trivial.
Cette trivialit{\'e} r{\'e}sulte du lemme suivant.

\begin{lemme}
Soit $\R\in S_{\und d}(P_1,\ldots,P_n)(k)$. Il existe des uniques vecteurs
$v_1,\ldots,v_n$ qui engendrent $\R$ et qui ont la forme
$$v_i =  P_ie_i +\sum_{j\leq i} R_{i,j}e_j$$
o{\`u} $\deg(R_{i,j})<\deg(P_j)$.
\end{lemme}

{\noindent\it D{\'e}monstration du lemme.}
Le vecteur $v_1=P_1e_1$ appartient {\`a} $\R$ par l'hypo\-th{\`e}se
$\R\in S_{\und d}(P_1,\ldots,P_n)(k)$. La m\^eme hypoth{\`e}se implique aussi
qu'il existe un vecteur de la forme $R'_{2,1}e_1+P_2 e_2$ appartenant {\`a} $\R$.
En utilisant la division euclidienne, on d{\'e}montre qu'il existe un unique
vecteur $v_2=R_{2,1}e_1+P_2 e_2\in \R$ tel que $\deg(R_{2,1})<\deg(P_1)$.
En continuant ainsi, on montre qu'il existe pour tout entier $i$ un unique
vecteur
$$v_i=P_ie_i +\sum_{j\leq i} R_{i,j}e_j \in\R$$
avec $\deg(R_{i,j})<\deg(P_j)$. Un calcul de d{\'e}terminant montre que ces
vecteurs engendrent $\R$. $\carre$

{\noindent\it Fin de la d{\'e}monstration du th{\'e}or{\`e}me.}
En combinant les lemmes pr{\'e}c{\'e}dents avec le th{\'e}o\-r{\`e}me de Jouanolou sur la cohomologie
{\`a} support propre des fibr{\'e}s vectoriels, on d{\'e}duit que le morphisme trace
$$\RR {\tilde s}_{{\und d},\tau,!}\Ql
\rta\Ql[-2r](-r)$$
o{\`u} $r=\sum_{i=1}^n d_i(n-i)$, est un isomorphisme.
En particulier, le complexe $\RR {\tilde s}_{{\und d},\tau,!}\Ql$ sur $\AA^d$
est concentr{\'e} en degr{\'e} $2r$ et a la propri{\'e}t{\'e} $(*)$ renforc{\'e}e.
Le morphisme $\pi_\tau:\AA^d\rta Q_{\und d}$ {\'e}tant fini, le complexe
$$\pi_{\tau,*}\RR {\tilde s}_{{\und d},\tau,!}\Ql$$
est donc aussi concentr{\'e} en degr{\'e} $2r$ et conserve la propri{\'e}t{\'e} $(*)$.

Du fait que les complexes
$\pi_{\tau,*}\RR {\tilde s}_{{\und d},\tau,!}\Ql$
sont concentr{\'e}s $2r$, le complexe
$\RR {\tilde s}_{{\und d},!}\Ql$ l'est aussi.
On a donc un isomorphisme dans la cat{\'e}gorie d{\'e}riv{\'e}e
$$\RR {\tilde s}_{{\und d},!}\Ql{\tilde\rta}
\RR^{2r} {\tilde s}_{{\und d},!}\Ql.$$

Or m{\^e}me si les adh{\'e}rences des strates
${\tilde S}_{{\und d},\tau}$ ne sont pas disjointes
dans ${\tilde S}_{\und d}$, on a en degr{\'e} maximal un isomorphisme
$$\RR^{2r}{\tilde s}_{{\und d},!}\Ql
{\tilde\rta}\bigoplus_{\scriptstyle\tau:\{1,\ldots,d\}\rta\{1,\ldots,n\}\atop\scriptstyle |\tau^{-1}(i)|=d_i}
\RR^{2r} {\tilde s}_{{\und d},\tau,!}\Ql.$$
On d{\'e}duit que $\RR {\tilde s}_{{\und d},!}\Ql$
a la propri{\'e}t{\'e} $(*)$.

Le complexe $\RR {\tilde s}_{{\und d},!}\Ql[nd](nd/2)$ est donc concentr{\'e} en
degr{\'e}
$$2\sum_{i=1}^n d_i(n-i)-nd=-d+2\<{\und d},\delta\>$$
et a aussi la propri{\'e}t{\'e} $(*)$. $\square$

\begin{corollaire}
Le complexe $\RR\Gamma_c(S_{{\und d},\vp}\otimes_k{\bar k},\A_\rho[-d](-d/2))$
est concentr{\'e} en degr{\'e} $2\<\und d,\delta\>$ et $\Fr$ agit dans
$H^{2\<\und d,\delta\>}_c(S_{\und d,\vp}\otimes_k{\bar k},\A_\rho[-d](-d/2))$
comme la multiplication par $q^{\<\und d,\delta\>}$.
\end{corollaire}

\begin{corollaire}
Soit $a_\rho$ la fonction trace de Frobenius du faisceau pervers
restriction de $\A_\rho[-d](-d/2)$ {\`a} $X_{d,\vp}$.
Pour tout ${\und d}\in \NN^n$ avec $|{\und d}|=d$, le terme constant
$a_\rho^B(\vp^{\und d})$ est un entier naturel {\'e}gal {\`a} la dimension
du groupe de cohomologie
$$H^{2\<{\und d},\delta\>}_c(S_{\und d,\vp},\A_\rho[-d](-d/2)).$$
\end{corollaire}

Lusztig a d{\'e}montr{\'e} que 
lorsque la repr{\'e}sentation $\rho=\rho_\lambda$ est irr{\'e}ductible,
$a^B_\lambda(\vp^{\und d})$ est {\'e}gale au coefficient de Kostka
$K_{\lambda,\und d}$ (\cite{Lus}).

\subsection{La cohomologie globale est la somme des cohomologie 
des termes constants}

\begin{proposition}
Pour chaque $\und d\in\NN^n$ avec $|\und d|=d$ notons
$m_{\und d}:Q_{\und d}\rta Q_d$ le morphisme
$$m_{\und d}(P_1,\ldots,P_n)=P_1\ldots P_n.$$
Pour toute repr{\'e}sentation $\rho$ de $\S_d$, on a alors un isomorphisme
$$\RR \phi_{*}\A_\rho{\tilde\rta}\bigoplus_{\scriptstyle\und d\in\NN^n\atop\scriptstyle |\und d|=d}
m_{\und d,*}\RR s_{\und d,!}i_{\und d}^*\A_\rho.$$
\end{proposition}

\DEMONSTRATION
Ces complexes ayant tous les deux la propri{\'e}t{\'e} $(*)$, il suffit de construire
l'isomorphisme au-dessus de l'ouvert dense de $Q_{d,rss}$ de $Q_d$.
Au-dessus $X_{d,rss}$, $\A_\rho$ est par d{\'e}finition {\`a} d{\'e}calage et torsion pr{\`e}s
le syst{\`e}me local $\L_\rho$. Puisque $\L_\rho$ provient de $Q_{d,rss}$,
par la formule de projection, il suffit de traiter le cas
o{\`u} $\rho$ est la repr{\'e}sentation triviale et $\L_\rho=\Ql$.
Il est aussi loisible de passer au rev{\^e}tement galoisien
$\AA^d_{rss}$ de $Q_{d,rss}$.

Au-dessus de $\AA^d_{rss}$, on a
$$X_d\times_{Q_d}\AA^d_{rss}={\tilde X}_{d,rss}.$$
La donn{\'e}e d'un $k$-point de ${\tilde X}_d$ au dessus de
$(x_1,\ldots,x_d)\in\AA^d_{rss}(k)$ est {\'e}quivalente {\`a} la donn{\'e}e pour tout $i$
d'un sous-espace vectoriel de codimension $1$ de $V_l={\cal O}^n/(\vp-x_l){\cal O}^n$.
La donn{\'e}e de $V_l$ d{\'e}finit un fibr{\'e} vectoriel ${\cal V}_l$ de rang $n$ sur $\AA^d$.
Les images de $e_1,\ldots e_n$ dans $V_l$ constituent une base de $V_l$ si bien
que le fibr{\'e} vectoriel ${\cal V}_l$ est trivial.
On a donc un isomorphisme
$${\tilde X}_{d,rss}{\tilde\rta}\AA^d_{rss}\times
\underbrace{\PP^{n-1}\times\cdots\times\PP^{n-1}}_{d}.$$

Notons $\PP_l^{n-1}$ le $l$-{\`e}me exemplaire de $\PP^{n-1}$.
Chaque $\PP^{n-1}_l$ admet une d{\'e}composition cellulaire
$$\PP^{n-1}_l=\bigcup_{i=0}^{n-1}\AA^i_l$$
o{\`u} les $k$-points de la cellule $\AA^{n-1-i}_l$ correspondent aux sous-espaces vectoriels 
de codimension $1$ de $V_l=\bigoplus_{j=1}^n ke_j$ qui contiennent
$\bigoplus_{j=1}^{i}ke_j$ mais ne contiennent pas $\bigoplus_{j=1}^{i+1}ke_j$.
On a donc
$$S_{\und d}\times_{Q_d}{\AA^d_{rss}}=
\coprod_{\scriptstyle\tau:\{1,\ldots,d\}\rta\{1,\ldots,n\}\atop\scriptstyle |\tau^{-1}(i)|=d_i}
\prod_{l=1}^n \AA^{n-1-\tau(l)}_l.$$

La proposition se d{\'e}duit donc de l'isomorphisme
$$\RR\Gamma(\PP^{n-1},\Ql){\tilde\rta}\bigoplus_{i=0}^{n-1}\RR\Gamma_c(\AA^i,\Ql)$$
bien connu. $\square$

\begin{corollaire}
On a un isomorphisme
$$\RR\Gamma(X_{d,\vp}\otimes_k{\bar k},\A_\rho){\tilde\rta}
\bigoplus_{\scriptstyle\und d\in\NN^n\atop \scriptstyle|\und d|=d}
\RR\Gamma_c(S_{\und d,\vp}\otimes_k{\bar k},\A_\rho).$$
\end{corollaire}

Sur $\CC$, ce corollaire est du {\`a} Mirkovic et Vilonen (\cite{M-V}).

\subsection{Termes constants du produit de convolution}

\begin{proposition}
Soient $d',d''$ deux entiers naturels, $\rho'$ et $\rho''$
respectivement des repr{\'e}\-sentations de $\S_{d'}$ et $\S_{d''}$.
Soient $d=d'+d''$ et ${\und d}\in\NN^n$ tel que $|{\und d}|=d$.
Pour chaque couple ${\und d}',{\und d}''\in\NN^n$ tel que
$|{\und d}'|=d'$, $|{\und d}''|=d''$ et ${\und d}'+{\und d}''={\und d}$
on a un morphisme
$$m_{{{\und d}'},{{\und d}''}}:Q_{{{\und d}'}}\times Q_{{{\und d}''}}\rta Q_{\und d}$$
d{\'e}fini par
$$m_{{{\und d}'},{{\und d}''}}(P'_1,\ldots,P'_n;P''_1,\ldots,P''_n)=(P'_1P''_1,\ldots,P'_nP''_n).$$
Avec les notations de 5.2, on a un isomorphisme
$$\RR s_{\und d,!}i_{\und d}^*(\A_{\rho'}*\A_{\rho''})
{\tilde\rta}\!\!
\bigoplus_{\scriptstyle|{{\und d}'}|=d',|{{\und d}''}|=d''\atop\scriptstyle{{\und d}'}+{{\und d}''}={\und d}}
m_{{{\und d}'},{{\und d}''},*}(\RR s_{{{\und d}'},!}i_{{{\und d}'}}^*\A_{\rho'}
\boxtimes \RR s_{{{\und d}''},!}i_{{{\und d}''}}^*\A_{\rho''}).$$
\end{proposition}

\DEMONSTRATION
Les morphisme $m_{{{\und d}'},{{\und d}''}}$ {\'e}tant finis,
les deux complexes ont la propri{\'e}t{\'e} $(*)$ gr{\^a}ce au th{\'e}or{\`e}me pr{\'e}c{\'e}dent.
Il suffit donc de construire l'isomorphisme au-dessus d'un ouvert dense.

Soit $Q_{\und d,rss}$ l'ouvert dense de $Q_{\und d}$ constitu{\'e} des points
$(P_1,\ldots,P_n)$ tels que le polyn{\^o}me $P_1\ldots P_n$ n'a pas de racines doubles.
On d{\'e}montre facilement que l'image inverse de
$$S_{\und d}\times_{Q_{\und d}}Q_{\und d,rss}$$
dans
$$(X_{d'}{\tilde\times}X_{d''})\times_{Q_{d}}Q_{\und d,rss}
=(X_{d'}\times X_{d''})\times_{Q_{d}}Q_{\und d,rss}$$
est une r{\'e}union de composantes connexes dont chacune
est isomorphe {\`a}
$$(S_{{{\und d}'}}\times S_{{{\und d}''}})\times_{Q_{\und d}}Q_{\und d,rss}$$
avec $|{{\und d}'}|=d',|{{\und d}''}|=d''$ et ${{\und d}'}+{{\und d}''}={\und d}$,
d'o{\`u} r{\'e}sulte la proposition.

Notons toutefois que les adh{\'e}rences de ces composantes dans
$$(X_{d'}{\tilde\times}X_{d''})\times_{Q_{d}}Q_{\und d}$$
ne sont plus disjointes. $\square$

\begin{corollaire}
On a un isomorphisme
$$\displaylines{\qquad
\RR\Gamma_c(S_{{\und d},\vp}\otimes_k{\bar k},\A_{\rho'}*\A_{\rho''})
{\tilde\rta}\bigoplus_{\scriptstyle|{{\und d}'}|=d',|{{\und d}''}|=d''\atop\scriptstyle{{\und d}'}+{{\und d}''}={\und d}}
\RR\Gamma_c(S_{{\und {d'}},\vp}\otimes_k{\bar k},\A_{\rho'})\hfill\cr\hfill
\otimes \RR\Gamma_c(S_{{{{\und d}''}},\vp}\otimes_k{\bar k},\A_{\rho''}).\qquad}$$
\end{corollaire}

Cet isomorphisme est compatible avec l'isomorphisme de commutativit{\'e}.

\begin{proposition}
Le diagramme
$$\diagramme{
\RR s_{\und d,!}i_{\und d}^*(\A_{\rho'}*\A_{\rho''}) &\rta&
        {\displaystyle\bigoplus_{\scriptstyle|{{\und d}'}|=d',|{{\und d}''}|=d''\atop\scriptstyle{{\und d}'}+{{\und d}''}={\und d}}}
        m_{{{\und d}'},{{\und d}''},*}(\RR s_{{{\und d}'},!}i_{{{\und d}'}}^*\A_{\rho'}
        \boxtimes \RR s_{{{\und d}''},!}i_{{{\und d}''}}^*\A_{\rho''}) \cr
\noalign{\vskip -6mm}\vfl{}{} & &\vfl{}{}\cr
\RR s_{\und d,!}i_{\und d}^*(\A_{\rho''}*\A_{\rho'}) &\rta &
        {\displaystyle\bigoplus_{\scriptstyle|{{\und d}'}|=d',|{{\und d}''}|=d''\atop\scriptstyle{{\und d}'}+{{\und d}''}={\und d}}}
        m_{{{\und d}''},{{\und d}'},*}(\RR s_{{{\und d}''},!}i_{{{\und d}''}}^*\A_{\rho''}
        \boxtimes \RR s_{{{\und d}'},!}i_{{{\und d}'}}^*\A_{\rho'})}$$
o{\`u} la fl{\`e}che verticale {\`a} gauche est induite par l'isomorphisme de commutativit{\'e}
$$\A_{\rho'}*\A_{\rho''}{\tilde\rta}\A_{\rho''}*\A_{\rho'}$$
et o{\`u}
la fl{\`e}che verticale {\`a} droite est la somme des isomorphismes {\'e}vidents
$$m_{{{\und d}'},{{\und d}''},*}(\RR s_{{{\und d}'},!}i_{{{\und d}'}}^*\A_{\rho'}
        \boxtimes \RR s_{{{\und d}''},!}i_{{{\und d}''}}^*\A_{\rho''}){\tilde\rta}
        m_{{{\und d}''},{{\und d}'},*}(\RR s_{{{\und d}''},!}i_{{{\und d}''}}^*\A_{\rho''}
        \boxtimes \RR s_{{{\und d}'},!}i_{{{\und d}'}}^*\A_{\rho'}),$$
est commutatif.
\end{proposition}

\DEMONSTRATION
Il suffit v{\'e}rifier la commutativit{\'e} du diagramme au-dessus d'un ouvert dense.
Mais au-dessus de $Q_{d,rss}$, l'isomorphisme de commutativit{\'e}
$$\A_{\rho'}*\A_{\rho''}{\tilde\rta}\A_{\rho''}*\A_{\rho'}$$
a {\'e}t{\'e} d{\'e}fini via les isomorphismes
\begin{eqnarray*}
(X_{d'}{\tilde\times}X_{d''})_{rss} & {\tilde\rta} &(X_{d'}{\times}X_{d''})_{rss}\cr
                                    & {\tilde\rta} &(X_{d''}{\times}X_{d'})_{rss}\cr
                                    & {\tilde\rta} &(X_{d''}{\tilde\times}X_{d'})_{rss}
\end{eqnarray*}
Il est clair que l'isomorphisme
$$(X_{d'}{\times}X_{d''}) {\tilde\rta} (X_{d'}{\times}X_{d''})$$
est compatible {\`a} l'isomorphisme
$$S_{\und {d'}}\times S_{\und {d''}}{\tilde\rta}S_{\und {d''}}\times S_{\und {d'}}$$
via lequel est d{\'e}fini l'isomorphisme
$$m_{{{\und d}'},{{\und d}''},*}(\RR s_{{{\und d}'},!}i_{{{\und d}'}}^*\A_{\rho'}
        \boxtimes \RR s_{{{\und d}''},!}i_{{{\und d}''}}^*\A_{\rho''}){\tilde\rta}
        m_{{{\und d}''},{{\und d}'},*}(\RR s_{{{\und d}''},!}i_{{{\und d}''}}^*\A_{\rho''}
        \boxtimes \RR s_{{{\und d}'},!}i_{{{\und d}'}}^*\A_{\rho'}),$$
d'o{\`u} r{\'e}sulte la proposition.

\vfill\eject
\part{Applications}

\section{L'homomorphisme $b:\H^+_r\rta\H^+$}
\subsection{Rappel}
 Soient $r$ un entier naturel, $k_r$ l'extension de degr{\'e} $r$ de $k$,
$F_{r,\vp}=k_r((\vp))$ et ${\cal O}_{r,\vp}=k_r[[\vp]]$.
Soit ${\cal H}^+_r$ l'alg{\`e}bre des fonctions {\`a} valeurs dans $\Ql$ et {\`a} support compact 
dans $\GL(n,F_{r,\vp})\cap\gl(n,{\cal O}_{r,\vp})$ qui sont bi-$\GL(n,{\cal O}_{r,\vp})$-invariantes.
Rappelons que l'isomorphisme de Satake correspondant
$$\Phi_r:{\cal H}_r^+\rta\Ql[t_1,\ldots,t_n]^{\frak S_n}$$
est d{\'e}fini par
$$\Phi_r(f)=\sum_{\und d\in\NN^n}f^B(\vp^{\und d})t^{\und d}$$
o{\`u}
$$f^B(\vp^{\und d})=q^{-r\<\und d,\delta\>}\int_{N(F_{r,\vp})}
f(\vp^{\und d} x)\d x$$
la mesure normalis{\'e}e $\d x$ attribuant {\`a} $N({\cal O}_{r,\vp})$ le volume $1$.

On d{\'e}finit l'homomorphisme de changement de base $b:{\cal H}^+_r\rta{\cal H}^+$
de telle mani{\`e}re que le diagramme
$$\diagramme{
{\cal H}^+_r & \hfl{\textstyle\Phi_r}{} & \Ql[t_1,\ldots,t_n]^{\frak S_n}\cr
 \vfl{\textstyle b}{}   &          & \vfl{}{}\cr
 {\cal H}^+  & \hfl{}{\textstyle\Phi} & \Ql[z_1,\ldots,z_n]^{\frak S_n}\cr
}$$
dont la fl{\`e}che verticale {\`a} droite envoie $t_i$ en $z_i^r$, est commutatif.

Pour chaque $n$-partition $\lambda$, on d\'efinit la fonction
$$a_{r,\lambda}:X_\vp(k_r)\rta\Ql$$
par
$$a(x)=\Tr(\Fr_{q^r},(\A_{\lambda,\vp})_x).$$
Puisque les faisceaux pervers $\A_{\lambda,\vp}$ sont $G$-\'equivariants,
les fonctions $a_{r,\lambda}$ s'identifient naturellement \`a des 
\'el\'ements de $\H^+_r$. Le but de cette section est
de donner une interpr\'etation g\'eom\'etrique pour les 
images  $b(a_{r,\lambda})$ des fonctions $a_{r,\lambda}$ par 
l'homomorphisme de changement de base $b:\H^+_r\rta\H^+$.

\subsection{Les fonctions $f_{r,\lambda}$ et $\phi_{r,\lambda}$}

Chaque fonction $X_{d,\vp}(k_r)\rta\Ql$ invariante relativement 
{\`a} l'action {\`a} gauche
de $G_{d,\vp}(k_r)$ s'identifie naturellement 
{\`a} un {\'e}l{\'e}ment de $\H^+_r$.
Toutefois, il vaut mieux de consid{\'e}rer 
les fonctions sur $\Res_{k_r/k}X_{d,\vp}(k)$.

Soit $\sigma$ l'endomorphisme de $X_{d,\vp}^r$ d{\'e}fini par
$$\sigma(x_1,x_2,\ldots,x_r)=(x_r,x_1,\ldots,x_{r-1}).$$
On peut identifier $X_{d,\vp}(k_r)$ {\`a} $\Fix(\Fr\circ\sigma,X_d(\vp^d)^r)$ 
en envoyant
$$x\mapsto (x,\Fr(x),\ldots,\Fr^{r-1}(x)).$$

Pour chaque $n$-partition $\lambda$ de $d$, notons $\A_{\lambda,\vp}$
le faisceau pervers qui est restriction de $\A_\lambda[-d](-d/2)$ 
{\`a} $X_{d,\vp}$.
D'apr{\`e}s 2.3, $\A_{\lambda,\vp}$ est le complexe d'intersection de
l'adh{\'e}rence de l'orbite $G_{d,\vp}\vp^\lambda$ dans $X_{d,\vp}$. 

Prenons $r$ copies $\A_{\lambda,\vp,1},\ldots,\A_{\lambda,\vp,r}$ de
$\A_{\lambda,\vp}$. L'endomorphisme $\sigma$ se rel{\`e}ve naturellement
sur l'isomorphisme
$${\tilde\sigma}:
        \A_{\lambda,\vp,r}\boxtimes\A_{\lambda,\vp,1}
        \boxtimes\cdots\boxtimes\A_{\lambda,\vp,r-1}
        \rta\A_{\lambda,\vp,1}\boxtimes\A_{\lambda,\vp,2}
        \boxtimes\cdots\boxtimes\A_{\lambda,\vp,r}.$$
qui est le produit tensoriel externe des isomorphismes identiques
$\A_{\lambda,\vp,i}\rta\A_{\lambda,\vp,i+1}$ l'indice $i$ \'etant 
prise dans les classes modulo $r$.
On d{\'e}finit la fonction $f_{r,\lambda}:\Fix(\Fr\circ\sigma,X_{d,\vp}^r)\rta\Ql$ par
$$f_{r,\lambda}(x)=\Tr(\Fr\circ{\tilde\sigma},
(\A_{\lambda,\vp,1}\boxtimes\cdots\boxtimes\A_{\lambda,\vp,r})_x).$$
Pour tout $\lambda$, $f_{r,\lambda}$ s'identifie naturellement {\`a} un 
{\'e}l{\'e}ment de ${\cal H}^+_r$.

On d{\'e}finit pour chaque $n$-partition $\lambda$ de $d$ une fonction
$\phi_{r,\lambda}\in{\cal H}^+$ comme suit.
La $r$-\`eme puissance convol\'ee 
$$\A_{\lambda,\vp}^r=\A_{\lambda,\vp,1}*\cdots *\A_{\lambda,\vp,r}$$
admet un automorphisme d'ordre $r$
$$\kappa': \A_{\lambda,\vp,1}*\cdots *\A_{\lambda,\vp,r}
{\buildrel\kappa\over\rta}
\A_{\lambda,\vp,r}*\A_{\lambda,\vp,1}*\cdots *\A_{\lambda,\vp,r-1}
{\buildrel\iota\over\rta}
\A_{\lambda,\vp,1}*\cdots *\A_{\lambda,\vp,r}$$
o\`u $\A_{\lambda,\vp,1},\ldots,\A_{\lambda,\vp,r}$ sont $r$ copies
de $\A_{\lambda,\vp}$, o\`u $\kappa$ est l'isomorphisme de 
commutativit\'e et o\`u $\iota$ se d\'eduit des isomorphismes 
\'evidents $\A_{\lambda,\vp,i}\rta\A_{\lambda,\vp,i+1}$
l'indice $i$ \'etant prise dans les classes modulo $r$. 
Pour tout $x\in X_{rd,\vp}(k)$, on pose
$$\phi_{r,\lambda}=\Tr(\Fr\circ {\kappa'},(\A_{\lambda,\vp}*\cdots*\A_{\lambda,\vp})_x).$$

\begin{proposition}
Quand $\lambda$ parcourt l'ensemble des $n$-partitions, les $f_{r,\lambda}$
forment une base de ${\cal H}^+_r$.
\end{proposition}

\DEMONSTRATION
Soit $c_\lambda\in{\cal H}^+_r$ la fonction caract{\'e}ristique de la double classe
$$\GL(n,{\cal O}_{r,\vp})\vp^{\lambda}\GL(n,{\cal O}_{r,\vp}).$$ 
Les $c_\lambda$ forment clairement
une base de ${\cal H}^+_r$.

Du fait que $\A_{\lambda,\vp}$ est le complexe d'intersection de l'adh{\'e}rence de
l'orbite $G_{d,\vp}\vp^\lambda$ dans $X_{d,\vp}$, on a
$$f_{r,\lambda}=q^{-r\<\lambda,\delta\>}(c_\lambda+\sum_{\mu<\lambda}z_{\lambda,\mu}c_\mu)$$
avec $z_{\lambda,\mu}\in\Ql$. On en d{\'e}duit que les $f_{r,\lambda}$ forment aussi
une base de ${\cal H}^+_r$. $\square$

\begin{proposition}
Pour toute $n$-partition $\lambda$, on a $b(f_{r,\lambda})=\phi_{r,\lambda}$.
\end{proposition}

Notons que si $F_r$ est le produit de $r$ copies de $F$, l'isomorphisme de
changement de base est d{\'e}fini par
$$f_1\boxtimes\cdots\boxtimes f_r\mapsto f_1*\cdots * f_r.$$

\subsection{L'identit\'e $b(f_{r,\lambda})=\phi_{r,\lambda}$}
Il revient au m{\^e}me de d{\'e}montrer que pour tout ${\und d}\in\NN^n$ on a
$$\int_{N(F_{r,\vp})}f_{r,\lambda}(\vp^{\und d}x_r)\d x_r
=\int_{N(F_\vp)}\phi_{r,\lambda}(\vp^{r\und d}x)\d x$$
et que pour tout $\und d\notin r\NN^n$ on a
$$\int_{N(F_\vp)}\phi_{r,\lambda}(\vp^{\und d}x)\d x=0.$$

Suivant une id{\'e}e de Serre, on peut appliquer la formule des traces de Grothendieck
au compos{\'e} de l'endomorphisme de Frobenius avec un endomorphisme d'ordre fini.

Plus pr{\'e}cis{\'e}ment, soit $\sigma$ un endomorphisme d'ordre fini d'un 
sch{\'e}ma $X$ de type fini sur $k=\FF_q$. Il existe alors un sch{\'e}ma 
$X_\sigma$ de type fini
sur $k$ muni d'un isomorphisme $X\otimes_k{\bar k}=X_\sigma\otimes_k{\bar k}$
tel que l'endomorphisme de Frobenius g{\'e}om{\'e}trique 
$\Fr'$ sur $X\otimes_k{\bar k}$
induit par la $k$-structure $X_\sigma$ est {\'e}gale {\`a} 
$\Fr\circ\sigma$ o{\`u} $\Fr$
est l'endomorphisme de Frobenius g{\'e}om{\'e}trique induit par 
la $k$-structure $X$. Si de plus, 
$\sigma$ se rel{\`e}ve sur un complexe de faisceaux $\ell$-adique $\C$
c'est-{\`a}-dire {\'e}tant donn{\'e} ${\tilde\sigma}:\sigma^*\C\rta\C$,
alors il existe aussi $\C_\sigma$ sur $X_\sigma$
muni d'un isomorphisme $\C\otimes_k{\bar k}=\C_\sigma\otimes_k{\bar k}$
tel que l'action ${\Fr'}$ sur $\C\otimes_k{\bar k}$ induit par $\C_\sigma$
est {\'e}gale $\Fr\circ{\tilde\sigma}$.

En appliquant la formule des traces de Grothendieck 
(\cite{Gro}) au pair $(X_\sigma,\C_\sigma)$,
on a
$$\sum_{x\in\Fix(\Fr\circ\sigma, X({\bar k})))}\Tr(\Fr\circ{\tilde\sigma},\C_x)
=\Tr(\Fr\circ{\tilde\sigma},\RR\Gamma_c(X\otimes_k{\bar k},\C)).$$

En appliquant la formule pr{\'e}c{\'e}dente {\`a} 
$(\sigma,{\tilde\sigma})$ agissant sur
$$(\underbrace{S_{\und d,\vp}\times\cdots\times S_{\und d,\vp}}_{r{\rm\ fois}},
\underbrace{\A_{\lambda,\vp}\boxtimes
\cdots\boxtimes\A_{\lambda,\vp}}_{r{\rm\ fois}})$$
o{\`u} par abus de notation $\A_{\lambda,\vp}$ 
d{\'e}signe son image r{\'e}ciproque par
 l'inclusion naturelle
$i_{\und d}:S_{\und d,\vp}\rta X_{d,\vp}$,
on obtient l'{\'e}galit{\'e}
$$\int_{N(F_{r,\vp})}f_{r,\lambda}(\vp^{\und d}x_r)\d x_r
=\Tr(\Fr\circ{\tilde\sigma},
\RR\Gamma_c(S_{\und d,\vp}\otimes_k{\bar k},\A_{\lambda,\vp})^{\otimes r}).$$
L'endomorphisme de Frobenius $\Fr^*$ agit dans
$\RR\Gamma_c(S_{\und d,\vp}\otimes_k{\bar k},\A_{\lambda,\vp})^{\otimes r}$
comme le produit tensoriel de son action dans chacun des facteurs.
L'endomorphisme ${\tilde\sigma}$ agit en permutant circulairement ces facteurs
$$\displaylines{\qquad{\tilde\sigma}
:\RR\Gamma_c(S_{\und d,\vp}\otimes_k{\bar k},\A_{\lambda,\vp,r})
\otimes\RR\Gamma_c(S_{\und d,\vp}\otimes_k{\bar k},\A_{\lambda,\vp,1})
\otimes\cdots\hfill\cr\hfill
\rta\RR\Gamma_c(S_{\und d,\vp}\otimes_k{\bar k},\A_{\lambda,\vp,1})
\otimes\RR\Gamma_c(S_{\und d,\vp}\otimes_k{\bar k},\A_{\lambda,\vp,2})
\otimes\cdots\qquad}$$

Appliquons maintenant la m{\^e}me formule au pair $(\Id,{\kappa'})$ 
agissant sur
$$(S_{\und d,\vp},
\underbrace{\A_{\lambda,\vp}*\cdots *\A_{\lambda,\vp}}_{r{\rm\ fois}}).$$
On obtient l'{\'e}galit{\'e}
$$\int_{N(F_\vp)}\phi_{r,\lambda}(\vp^{\und d}x)\d x
=\Tr(\Fr^*\circ {\kappa'},
\RR\Gamma_c(S_{\und d,\vp}\otimes_k{\bar k},
\underbrace{\A_{\lambda,\vp}*\cdots *\A_{\lambda,\vp}}_{r{\rm\ fois}}))).$$

Or, d'apr{\`e}s 5.4.2, on a
$$\displaylines{\RR\Gamma_c(S_{\und d,\vp}\otimes_k{\bar k},
\underbrace{\A_{\lambda,\vp,1}*\cdots *\A_{\lambda,\vp,r})}_{r{\rm\ fois}})\cr
=\bigoplus_{\scriptstyle \und{d_1},\ldots,\und{d_r}\atop
                \und{d_1}+\cdots+\und{d_r}=\und d}
\bigotimes_{i=1}^r
\RR\Gamma_c(S_{\und{d_i},\vp}\otimes_k{\bar k},\A_{\lambda,\vp,i}).}$$
Examinons l'action de $\Fr$ et de ${\kappa'}$ dans l'expression de droite.
L'action de $\Fr$ laisse stable chacun des termes de la somme directe. 
En revanche, d'apr{\`e}s 5.4.3, on a ${\kappa'}=\iota\circ\kappa$ 
o\`u $\kappa$ envoie le terme
$$\RR\Gamma_c(S_{\und{d_1},\vp}\otimes_k{\bar k},\A_{\lambda,\vp,1})
\otimes\cdots\otimes
\RR\Gamma_c(S_{\und{d_r},\vp}\otimes_k{\bar k},\A_{\lambda,\vp,r})
$$
sur le terme
$$\displaylines{\qquad
\RR\Gamma_c(S_{\und{d_r},\vp}\otimes_k{\bar k},\A_{\lambda,\vp,r})
\otimes \RR\Gamma_c(S_{\und{d_1},\vp}\otimes_k{\bar k},\A_{\lambda,\vp,1})
\hfill\cr\hfill
\otimes\cdots\otimes
\RR\Gamma_c(S_{\und{d_{r-1}},\vp}\otimes_k{\bar k},\A_{\lambda,\vp,{r-1}})
\qquad}$$
lequel s'envoie lui-m\^eme par $\iota$ sur le terme
$$\displaylines{\qquad
\RR\Gamma_c(S_{\und{d_r},\vp}\otimes_k{\bar k},\A_{\lambda,\vp,1})
\otimes \RR\Gamma_c(S_{\und{d_1},\vp}\otimes_k{\bar k},\A_{\lambda,\vp,2})
\hfill\cr\hfill
\otimes\cdots\otimes
\RR\Gamma_c(S_{\und{d_{r-1}},\vp}\otimes_k{\bar k},\A_{\lambda,\vp,{r}}).
\qquad}$$

Si $\und{d_1},\ldots,\und{d_r}$ ne sont pas tous {\'e}gaux, la somme des $r$
termes
$$\bigoplus_{j=0}^{r-1}\bigotimes_{i=1}^r
\RR\Gamma_c(S_{\und{d_{i+j}},\vp},\A_{\lambda,\vp,i})$$
$i+j$ {\'e}tant pris modulo $r$, est stable par ${\kappa'}$.
Puisque l'automorphisme $\kappa'$ y agit par une permutation
circulaire et que l'endomorphisme de Frobenius laisse stable 
chacun de ses termes, le compos\'e $\Fr\circ\kappa'$
a la trace nulle sur cette somme directe.

On en d{\'e}duit que si $\und d\notin r\NN^n$, on a
$$\int_{N(F_\vp)}\phi_{r,\lambda}(\vp^{\und d}x)\d x=0$$
Si maintenant $\und d=r\und d'$, on a
$$\int_{N(F_{r,\vp})}f_{r,\lambda}(\vp^{{\und d}'}x_r)\d x_r=
\int_{N(F_\vp)}\phi_{r,\lambda}(\vp^{\und d}x)\d x$$
du fait que les actions de ${\tilde\sigma}$ et ${\kappa'}$
dans le terme
$$\RR\Gamma_c(S_{{{\und d}'},\vp},\A_{\vp,\lambda})^{\otimes r}$$
sont {\'e}gales. $\square$

\subsection{L'identit\'e $a_{r,\lambda}=f_{r,\lambda}$}

\begin{proposition}
Pour toute $n$-partition $\lambda$, on a $a_{r,\lambda}=f_{r,\lambda}$.
\end{proposition}

En combinant avec l'identit\'e $b(f_{r,\lambda})=\phi_{r,\lambda}$, on 
obtient le r\'esultat principal de cette section.

\begin{theoreme}
Pour toute $n$-partition $\lambda$ on a $b(a_{r,\lambda})=\phi_{r,\lambda}$.
\end{theoreme}

\noindent{\it D\'emonstration de la proposition.} La transformation
de Satake 
$$\Phi_r:\H^+_r\rta\Ql[t_1,\ldots,t_n]^{\S_n}$$
\'etant un isomorphisme, il suffit de d\'emontrer que 
$$\Phi_r(a_{r,\lambda})=\Phi_r(f_{r,\lambda}).$$ 
On va donc comparer les termes constants de $a_{r,\lambda}$ et
de $f_{r,\lambda}$. Il suffit donc de d\'emontrer que pour tout
${\und d}\in\NN^n$ avec $|\und d|=d$, on a
$$\Tr(\Fr^r,\RR\Gamma_c(S_{\und d,\vp}\otimes_k{\bar k},\A_{\lambda,\vp}))
=\Tr(\Fr\circ{\tilde\sigma},\RR\Gamma_c
(S_{\und d,\vp}\otimes_k{\bar k},\A_{\lambda,\vp})^{\otimes r}).$$

D'apr\`es le corollaire 5.2.3, on sait que le complexe 
$\RR\Gamma_c(S_{\und d,\vp}\otimes_k{\bar k},\A_{\lambda,\vp})$
est concentr\'e en degr\'e $2\<{\und d},\delta\>$ et que $\Fr$
agit dans le groupe de cohomogie
$$V={\rm H}^{2\<{\und d},\delta\>}_c
(S_{\und d,\vp}\otimes_k{\bar k},\A_{\lambda,\vp})$$
comme la multiplication par $q^{\<\und d,\delta\>}$.
On est amen\'e \`a d\'emontrer que
$\dim(V)=\Tr({\tilde\sigma},V^{\otimes^r})$.

On peut d\'ecomposer $V^{\otimes r}=V\oplus V'$ o\`u $V$ est le 
sous-espace vectoriel diagonal et o\`u $V'$ est un suppl\'ementaire stable
sous l'action de ${\tilde\sigma}$. On d\'emontre facilement que 
$\Tr({\tilde\sigma},V')=0$.

\section{La transposition}
\subsection{Les sch{\'e}mas $\g_{d,r}$}

D\'esormais, on ne consid\`ere que les entensions de degr\'e
$r=2$. On peut ainsi lib\'erer la lettre $r$ pour d'autres utilisations.

Pour chaque couple d'entiers naturels $(d,r)\in\NN^2$ avec $d<r$,
notons $Q_{d,r}$ l'espace affine dont l'ensemble des $k$-points
est celui des couples de polyn{\^o}mes unitaires $(P,R)$ de degr{\'e}
respectivement $d$ et $r$ tels que $P$ divise $R$. En fait,
$Q_{d,r}$ est isomorphe {\`a} $Q_d\times Q_{r-d}$.

Soit $\g_{d,r}$ le $Q_{d,r}$-sch{\'e}ma dont l'ensemble des $k$-points
au-dessus de $(P,R)\in Q_{d,r}(k)$ est l'ensemble
$$\{g\in\gl(n,\OO/(R))\mid\det(g)\in P(\OO/(R))^\times\}.$$
Pour chaque $g\in\g_{d,r}(P,R)(k)$, si on note $\R$ le r{\'e}seau image
inverse de $g(\OO/(R))^n$ par l'application $\OO^n\rta(\OO/(R))^n$,
on a $\R\in X_d(P)(k)$. Cela d{\'e}finit un morphisme $p:\g_{d,r}\rta X_d$.

Notons $G_{d,r}=G_r\times Q_{d,r}$. Le $Q_{d,r}$-sch{\'e}ma en groupes
$G_{d,r}$ agit de deux c{\^o}t{\'e}s sur $\g_{d,r}$ par
$(h_1,h_2).g=h_1 gh_2^{-1}$.

\begin{lemme}
Le morphisme $p$ est lisse et {\`a} fibres g{\'e}om{\'e}triquement connexes.
Il est invariant relativement {\`a} l'action {\`a} droite de
$G_{d,r}$. Cette action est transitive sur ses fibres
g{\'e}om{\'e}triques. Enfin, le morphisme $p$ est {\'e}qui\-variant
relativement {\`a} l'action {\`a} gauche de $G_{d,r}$ sur $\g_{d,r}$
et l'action de $G_d$ sur $X_d$. 
\end{lemme}

\DEMONSTRATION
Montrons d'abord que pour tout $(P,R)\in Q_{d,r}(k)$, pour tout $\R\in\phi^{-1}(P)$
o{\`u} $\phi:X_d\rta Q_d$ est le morphisme d{\'e}terminant,
le groupe $G_d(R)(k)$ agit transitivement sur $p^{-1}(\R)(k)$.
Soient $g$ et $g'$ dans $\gl(n,{\cal O}/(R))$ tels que
$$g({\cal O}^n/(R)^n)=g'({\cal O}^n/(R)^n)=\R/(R)^n$$
En utilisant le lemme de Nakayama, on peut relever $g$ et $g'$ en ${\tilde g}$
et ${\tilde g}'$ dans $\gl(n,{\cal O}_{(P)})$ o{\`u} ${\cal O}_{(P)}$ 
est le localis{\'e} de ${\cal O}$ en $(P)$, tels que
$$g{\cal O}_{(R)}^n=g'{\cal O}_{(R)}^n=\R\OO_{(R)}.$$
Mais alors on a
$${\tilde h}={\tilde g}^{-1}{\tilde g}'\in\GL(n,{\cal O}_{(P)}).$$
Il est alors clair que $gh=g'$
o{\`u} $h\in \GL(n,{\cal O}/(R))$ est la r{\'e}duction de $h$ modulo $(R)$. 

L'action de $G_d$ {\'e}tant transitive sur les fibres de $\pi$, ces fibres sont
automatiquement lisses. Calculons leur dimension. Un r{\'e}seau $\R\in\phi^{-1}(P)({\bar k})$
{\'e}tant fix{\'e}, la donn{\'e}e de $g\in\gl(n,{\cal O}/(R))$ telle que
$g({\cal O}^n/(R)^n)=\R/(R)^n$ est {\'e}quivalente {\`a} celle de $n$ vecteurs dans $\R/(P)^n$
qui engendrent ce module.  La derni{\`e}re condition {\'e}tant ouverte, la fibre $\pi^{-1}(\R)$
est de dimension $n(nr-d)$.

Pour d{\'e}montrer que $p$ est lisse, il suffit maintenant de 
d{\'e}montrer que la source $\g_{d,r}$
est lisse. Dans la d{\'e}monstration du lemme 2.2.2, on a construit un
recouvrement d'ouverts lisses $U$ de $X_d$ munis des sections $U\rta\g_{d,r}$.
En utilisant l'action transitive de $G_{d,r}$, on a un morphisme surjectif
$U\times_{Q_{d,r}} G_{d,r}\rta\g_{d,r}$.  
On v{\'e}rifie facilement que ses fibres g{\'e}om{\'e}triques 
sont lisses et ont la m{\^e}me dimension. Du fait que sa source est lisse, ce morphisme est
lisse. Il en est de m{\^e}me de son but. $\carre$

Pour toute repr{\'e}sentation $\ell$-adique $\rho$ de ${\frak S}_d$, on notera
$\tA_\lambda=p^*\A_\lambda$. On d{\'e}duit du lemme pr{\'e}c{\'e}dent que le
complexe $\tA_\lambda$ est bi-$G_{d,r}$-{\'e}quivariant et est un faisceau
pervers d{\'e}cal{\'e}.

\subsection{Les rel{\`e}vements
        ${\tilde\tau}_\rho$ de $\tau$}

Notons $\tau:\g_{d,r}\rta\g_{d,r}$ l'involution d{\'e}finie par
$$(g,P,R)\mapsto(\,^\t g,P,R).$$
Pour chaque repr{\'e}sentation $\rho$ du groupe sym{\'e}trique $\S_d$,
on va relever $\tau$ sur $\tA_\rho$ comme suit.

Soit $Q_{d,r,rss}$ l'ouvert de $Q_{d,r}$ dont l'ensemble des 
$k$-points est celui des couples de polyn{\^o}mes $(P,R)$  
avec $P$ s{\'e}parable divisant $R$.
Au-dessus de $\g_{d,r,rss}=\g_{d,r}\times_{Q_{d,r}}Q_{d,r,rss}$,
le complexe $\tA_\rho$ provient d'un syst{\`e}me local $\L_\rho$
sur $Q_{d,r,rss}$ par l'image r{\'e}ciproque du morphisme lisse
$$p_{rss}:\g_{d,rss}\rta Q_{d,rss}.$$
Or puisque $p\circ\tau=p$, on a isomorphisme
$$\tau^* p_{rss}^*\L_\rho\rta p_{rss}^*\L_\rho$$
qu'on {\'e}tend par le prolongement interm{\'e}diaire pour obtenir
$${\tilde\tau}_\rho:\tau^*\tA_\rho\rta\tA_\rho.$$

Leur d\'efinition apparemment simple cache certains aspects assez
myst\'erieux de ces ${\tilde\tau}_\rho$. On reporte cette discussion
\`a la derni\`ere section de l'article.

\subsection{Convolution}

Soient $d',d'',r$ des entiers naturels tels que $d'+d''=d<r$.
Soit $Q_{d',d'',r}$ l'espace affine dont l'ensemble des
$k$-points est celui des triplets de polyn{\^o}mes unitaires
$(P',P'',R)$ de degr{\'e} respectivement $d',d''$ et $r$ tels que
le produit $P'P''=P$ divise $R$. 
En fait $Q_{d',d'',r}$ est isomorphe {\`a}
$Q_{d'}\times Q_{d''}\times  Q_{r-d'-d''}$. Notons
$m_Q:Q_{d',d'',r}\rta Q_{d,r}$ le morphisme d{\'e}fini par
$$(P',P'',R)\mapsto(P,R) $$
o{\`u} $P=P'P''$.

Soit $\g_{d',d'',r}$ le $Q_{d',d'',r}$-sch{\'e}ma dont l'ensemble des
$k$-points au-dessus de $(P',P'',R)\in Q_{d',d'',r}(k)$ est
l'ensemble
$$\{g',g''\in\gl(n,\OO/(R))\mid\det(g')\in P'(\OO/(P))^\times,
                \det(g'')\in P''(\OO/(R))^\times\}.$$
On a aussi un morphisme $m:\g_{d',d'',r}\rta\g_{d,r}$
au-dessus de $m_Q$ d{\'e}fini par
$$(g',g'';P',P'',R)\rta(g;P,R)$$
o{\`u} $g=g'g''$ et $P=P'P''$.

Notons $G_{d',d'',r}=G_r\times_{Q_r}Q_{d',d'',r}$. Le morphisme
$m$ est invariant relativement {\`a} l'action $\alpha$ de
$G_{d',d'',r}$ sur $Q_{d',d'',r}$ d{\'e}fini par
$$\alpha(h,(g',g''))=(g'h^{-1},hg'').$$
Le quotient de $\g_{d',d'',r}$ par cette action peut {\^e}tre d{\'e}finie par 
$$\g_{d',d'',r}^\alpha(k)
=\{(\R',g)\in X_{d'}\times\g_{d,r}(k)\mid \R'\supset g\OO^n\}$$
o{\`u} $g\OO^n$ d{\'e}signe le r{\'e}seau image inverse de $g(\OO/(R))^n$
par l'application $\OO^n\rta(\OO/(R))^n$.

\begin{lemme}
Le morphisme $\pi_\alpha:\g_{d',d'',r}\rta\g_{d',d'',r}^\alpha$ d{\'e}fini par
$$(g',g'',P',P'',R)\mapsto((g'\OO^n,P'),(g'g'',P'P'',R))$$
est lisse et invariant relativement {\`a} l'action $\alpha(G_{d',d'',r})$.
De plus, cette action est transitive sur ses fibres g{\'e}om{\'e}triques.
\end{lemme}

Notons $q:\g_{d',d'',r}\rta \g_{d'}\times \g_{d''}$ le morphisme
d{\'e}fini par
$$(g',g'',P',P'',R)\rta ((g',P',R),(g'',P'',R)).$$
Notons $p'$ et $p''$ les morphismes $\g_{d'}\rta X_{d'}$
et $\g_{d''}\rta X_{d''}$.

\begin{lemme}
Le morphisme compos{\'e} $(p'\times p'')\circ q$ 
est lisse {\`a} fibres g{\'e}o\-m{\'e}\-tri\-quement connexes.
\end{lemme}

Ces deux derniers lemmes se d{\'e}montrent de mani{\`e}re tr{\`e}s analogue
au lemme 7.1.1. On omet leurs d{\'e}monstrations.

Pour toutes $\rho'$ et $\rho''$ repr{\'e}sentations $\ell$-adiques
respectivement de $\S_{d'}$ et de $\S_{d''}$, 
le complexe image inverse 
$q^*(p'\times p'')^*(\A_{\rho'}\boxtimes\A_{\rho''})$
est donc un faisceau pervers {\`a} d{\'e}calage pr{\`e}s.
Ce complexe est clairement $\alpha(G_{d',d'',r})$-{\'e}quivariant
si bien qu'il existe un faisceau pervers d{\'e}cal{\'e}
$\tA_{\rho'}\boxtimes^{\,\alpha}\tA_{\rho''}$
sur $\g_{d',d'',r}^\alpha$ d{\'e}fini {\`a} un unique isomorphisme pr{\`e}s
tel que
$$q^*(p'\times p'')^*(\A_{\rho'}\boxtimes\A_{\rho''})
{\tilde\rta}\pi_\alpha^*(\tA_{\rho'}
        \boxtimes^{\,\alpha}\tA_{\rho''}).$$

\begin{proposition}
\begin{enumerate}
\item Soit $p^\alpha:\g_{d',d'',r}^\alpha
\rta X_{d'}{\tilde\times}X_{d''}$ le morphisme d{\'e}fini par
$$(\R',g,P,R)\mapsto (\R',\R)$$
avec $\R=g\OO^n$. On a un isomorphisme
$$p^{\alpha,*}(\A_{\rho'}{\tilde\boxtimes}\A_{\rho''})
{\tilde\rta} \tA_{\rho'}\boxtimes^{\,\alpha}\tA_{\rho''}.$$
\item Le diagramme
$$\diagramme{
        \g_{d',d'',r}^\alpha & \hfl{m^\alpha}{} & \g_{d,r} \cr
        \vfl{p^\alpha}{}             &          & \vfl{}{p} \cr
        X_{d'}{\tilde\times}X_{d''} &\hfl{}{\mu}& X_d \cr     
        }$$
o{\`u} $m^\alpha$ est le morphisme induit 
par le morphisme $\alpha$-invariant 
$$m:\g_{d',d'',r}\rta\g_{d,r},$$
est cart{\'e}sien.
\item
On a un isomorphisme
$$\RR m^\alpha_*(\tA_{\rho'}\boxtimes^{\,\alpha}\tA_{\rho''}) 
{\tilde\rta} \tA_\rho$$
o{\`u} $\rho$ est la repr{\'e}sentation induite 
$$\rho=\Ind_{\S_{d'}\times\S_{d''}}^{\S_d}(\rho'\times\rho'').$$
\end{enumerate}
\end{proposition}

\DEMONSTRATION
\begin{enumerate}
\item 
Le faisceau pervers d{\'e}cal{\'e}  
$p^{\alpha,*}(\A_{\rho'}{\tilde\boxtimes}\A_{\rho''})$
v{\'e}rifie la propri{\'e}t{\'e} suivante :
au-dessus de $\g_{d',d'',r}$ on a un isomorphisme
entre les faisceaux pervers d{\'e}cal{\'e}s
$$q^*(p'\times p'')^*(\A_{\rho'}\boxtimes\A_{\rho''}){\tilde\rta}
\pi^{\alpha,*} p^{\alpha,*}(\A_{\rho'}{\tilde\boxtimes}\A_{\rho''}).$$
En effet, au-dessus de l'ouvert dense o{\`u}
le polyn{\^o}me $P'P''$ est s{\'e}parable, on a un isomorphisme entre
${X_{d'}\times X_{d''}}$ et $X_{d'}{\tilde\times}X_{d''}$, 
d'o{\`u} on d{\'e}duit l'isomorphisme  
$$q^*(p'\times p'')^*(\A_{\rho'}\boxtimes\A_{\rho''}){\tilde\rta}
\pi^{\alpha,*} p^{\alpha,*}(\A_{\rho'}{\tilde\boxtimes}\A_{\rho''})$$
au-dessus de cet ouvert. Puis on l'{\'e}tend {\`a} $\g_{d',d'',r}$
par le prolongement interm{\'e}diaire.

Or, $\tA_{\rho'}\boxtimes^{\,\alpha}\tA_{\rho''}$ est d{\'e}fini
par cette propri{\'e}t{\'e} {\`a} l'unique isomorphisme pr{\`e}s, d'o{\`u} 
l'assertion.

\item
C'est trivial sur la d{\'e}finition des fl{\`e}ches du diagramme. 

\item
Rappelons qu'on a
$$\tA_\rho=p^*(\A_{\rho'}*\A_{\rho''})
=p^*\RR\mu_*(\A_{\rho'}{\tilde\boxtimes}\A_{\rho''})$$
d'apr{\`e}s la proposition 3.2.1.
Il suffit d'appliquer le th{\'e}or{\`e}me de changement de base 
pour un morphisme propre au diagramme cart{\'e}sien pr{\'e}c{\'e}dent. 
\end{enumerate}

Gr\^ace \`a cette nouvelle description du produit de convolution,
on voit facilement que cette notion correspond bien
au produit de convolution habituel entre fonctions de Hecke,
via le dictionnaire faisceaux-fonctions de Grothendieck.

\subsection{Commutativit{\'e} revue}
A l'aide des ${\tilde\tau}_\rho$, on peut traduire g{\'e}om{\'e}triquement
le calcul habituel suivant
\begin{eqnarray*}
(f*g)(x) &=& \int_{G(F_\vp)}f(xy^{-1})g(y)\d y\cr
         &=& \int_{G(F_\vp)}f(\,^\t y^{-1}\,^\t x)g(\,^\t y)\d y\cr
         &=& \int_{G(F_\vp)}g(\,^\t y)f(\,^\t y^{-1}\,^\t x)\d y\cr
         &=& (g*f)(\,^\t x)\cr
         &=& (g*f)(x)\cr
\end{eqnarray*}
qui d{\'e}montre la commutativit{\'e} du produit de convolution entre
les fonctions de Hecke.

Consid{\'e}rons le diagramme cart{\'e}sien suivant
$$\diagramme{
\g_{d',r}\times \g_{d'',r} & \hgfl{q}{} & \g_{d',d'',r} &
        \hfl{\pi^\alpha}{} & \g^\alpha_{d',d'',r}
        &\hfl{m^\alpha}{} & \g_{d,r} \cr
\vfl{\sigma\circ(\tau'\times\tau'')}{} & &
        \vfl{\sigma\tau}{}
& & \vfl{}{(\sigma\tau)^\alpha} & &\vfl{}{\tau}\cr
\g_{d'',r}\times \g_{d',r} & \hgfl{}{q_1} & \g_{d'',d',r} &
        \hfl{}{\pi^\alpha_1} & \g^\alpha_{d'',d',r}
        &\hfl{}{m_1^\alpha} & \g_{d,r} \cr
}$$
o{\`u} $\sigma$ est d{\'e}fini par $(g',g'')\mapsto(g'',g')$,
o{\`u} $\tau'$ (resp. $\tau''$) est d{\'e}fini par 
$g'\mapsto\,^\t g'$ (resp. $g''\mapsto\,^\t g''$),
o{\`u} 
$\sigma\tau(g',g'')=(\,^\t g'',\,^\t g')$,
o{\`u}
$(\sigma\tau)^\alpha$ est le morphisme induit par le
morphisme $\alpha$-invariant $\sigma\tau$
et o{\`u} les indices $1$ accol{\'e}es aux morphismes de ligne inf{\'e}rieure
est destin{\'e}e {\`a} distinguer ces morphismes de ceux de la ligne 
sup{\'e}rieure.   

Soient maintenant $\rho'$ et $\rho''$ deux repr{\'e}sentations
$\ell$-adiques respectivement de ${\frak S}_{d'}$ et de 
${\frak S}_{d''}$.
On a un isomorphisme canonique
$${\tilde\sigma}:\sigma^*(\tA_{\rho''}\boxtimes\tA_{\rho'})
{\tilde\rta}\tA_{\rho'}\boxtimes\tA_{\rho''}.$$
Notons
$${\widetilde{\sigma\tau}}_{\rho'\times\rho''}:
      (\sigma\tau)^*q_1^*(\tA_{\rho''}\boxtimes\tA_{\rho'})
      {\tilde\rta} q^*(\tA_{\rho'}\boxtimes\tA_{\rho''})$$
l'isomorphisme qui se d{\'e}duit de 
$({\tilde\tau}_{\rho'}\times{\tilde\tau}_{\rho''})\circ{\tilde\sigma}$
par $q^*$.
Cet isomorphisme {\'e}tant $\alpha$-invariante, il induit un isomorphisme
$${\widetilde{\sigma\tau}}^\alpha_{\rho'\times\rho''}
   : ({\sigma}{\tau})^{\alpha,*}
   (\tA_{\rho''}\boxtimes^\alpha \tA_{\rho'})
   {\tilde\rta}(\tA_{\rho'}\boxtimes^\alpha \tA_{\rho''}).$$

\begin{proposition}
On a un diagramme commutatif 
$$\diagramme{
\tau^*\RR m_{1,*}^\alpha(\tA_{\rho''}\boxtimes^\alpha\tA_{\rho'})
    &\hfl{\sim}{} 
    &\tau^* \tA_{\rho_1} \cr
\vfl{{\rm c-d-b}}{} & & \vfl{}{\tau^* (\kappa)}\cr
\RR m_{1,*}^\alpha (\sigma\tau)^{\alpha,*}
        (\tA_{\rho''}\boxtimes^\alpha \tA_{\rho'})
    &  
    &\tau^* \tA_{\rho}\cr
\vfl{\RR m_*^\alpha({\widetilde{\sigma\tau}}^\alpha_{\rho'\times\rho''})}{} 
    & 
    & \vfl{}
        {\tilde\tau}_{\,\,\rho}\cr
\RR m_{*}^\alpha (\tA_{\rho'}\boxtimes^\alpha \tA_{\rho''})
    &\hfl{}{\sim}
    &\tA_\rho\cr}$$
o{\`u} les deux fl{\`e}ches horizontales sont des isomorphismes de la
proposition 7.3.3,
o{\`u} ${\rm c-d-b}$ est l'isomorphisme de changement de base pour
le morphisme propre $m^\alpha$,
o{\`u} $\rho$ et $\rho_1$ sont les repr{\'e}sentations induites 
$$\displaylines{
\rho=\Ind_{\S_{d'}\times\S_{d''}}^{\S_d}(\rho'\times\rho'')\cr
\rho_1=\Ind_{\S_{d''}\times\S_{d'}}^{\S_d}(\rho''\times\rho')}
$$
pour lesquelles on a
$$\displaylines{
\tA_\rho=p^*(\A_{\rho'}*\A_{\rho''})\cr
\tA_{\rho_1}=p^*(\A_{\rho''}*\A_{\rho'})}$$
et o{\`u} $\kappa$ se d{\'e}duit de l'isomorphisme de commutativit{\'e}
$$\kappa:\A_{\rho''}*\A_{\rho'}\,{\tilde\rta}\,\A_{\rho'}*\A_{\rho''}.$$
\end{proposition}

\DEMONSTRATION
Les fl{\`e}ches du diagramme {\'e}tant tous des isomorphismes entre faisceaux
pervers d{\'e}cal{\'e}s, il suffit de d{\'e}montrer que le diagramme commute
au-dessus de l'ouvert $\g_{d,r,rss}$ o{\`u} le polyn{\^o}me $P$ est s{\'e}parable.

Au-dessus de cet ouvert, les $\tA_\rho$, $\tA_\rho'$, $\tA_\rho''$ proviennent
respectivement de $Q_{d,rss}$, $Q_{d',rss}$, $Q_{d'',rrs}$.
Par construction, les ${\tilde\tau}_\rho$ sont triviaux au niveau 
de $Q_{d,rss}$, $Q_{d',rss}$, $Q_{d'',rrs}$. Il s'agit de v\'erifier
la commutativit\'e du diagramme
$$\diagramme{
\RR m_{1,Q,*}(\L_{\rho''}\boxtimes\L_{\rho'})
& {\tilde\rta}& \L_{\rho''}*\L_{\rho'}\cr
\vfl{\RR m_{Q,*}({\tilde\sigma})}{} & & \vfl{}{\kappa}\cr
\RR m_{Q,*}\sigma^*(\L_{\rho'}\boxtimes\L_{\rho''})
& {\tilde\rta}& \L_{\rho'}*\L_{\rho''}\cr
}$$
Mais c'est pr\'ecis\'ement la d\'efinition de 
l'isomorphisme de commutativit\'e $\kappa$.

\subsection{Compatibilit{\'e}}
Soient $(P,R)\in Q_{d,r}({\bar k})$ avec $P=\prod_{i=1}^r(\vp-\gamma_j)^{d_j}$
et $R=\prod_{i=1}^r(\vp-\gamma_j)^{r_j}$ o{\`u} $\gamma_j\in{\bar k}$
et o{\`u} $d_j,r_j\in\NN$ avec $d_j<r_j$.

Soit $\rho$ une repr{\'e}sentation $\ell$-adique de $\S_d$. Sa restriction {\`a}
$\S_{d_1}\times\cdots\S_{d_r}$ admet une d{\'e}composition
$$\Res_{\S_{d_1}\times\cdots\S_{d_n}}^{\S_d}\rho
=\bigoplus_i\bigotimes_j \rho_{i,j}$$
o{\`u} chaque $\rho_{i,j}$ est une repr{\'e}sentation irr{\'e}ductible de $\S_{d_j}$.

\begin{proposition}
On a un isomorphisme
$$\g_{d,r}(P,R)=\prod_{j=1}^r\g_{d_j,r_j}((\vp-\gamma_j)^{d_j},(\vp-\gamma_j)^{r_j})$$
via lequel on a un isomorphisme
$$\tA_\rho|_{\g_{d,r}(P,R)}=\bigoplus_i\bigotimes_{j=1}^r
        \tA_{\rho_{i,j}}|_{\g_{d_j,r_j}((\vp-\gamma_j)^{d_j},(\vp-\gamma_j)^{r_j})}.$$

De plus l'action de ${\tilde\tau}_\rho$ sur le membre de gauche s'identifie
{\`a} l'action de 
$\bigoplus_i \bigotimes_{j=1}^r{\tilde\tau}_{\rho_{i,j}}$
sur le membre de droite.
\end{proposition}

La d{\'e}monstration est analogue 
{\`a} celle de la proposition 3.4.1. 

\section{L'application $b':\H^+_2\rta{\H'}^+$}
\subsection{L'identit{\'e} $b(f_{\lambda})=\phi_{\lambda}$ : rappel}

Soit $k_2$ l'extension quadratique de $k$ contenue dans ${\bar k}$.
Notons $F_{\vp,2}=k_2((\vp))$ et $\OO_{\vp,2}=k_2[[\vp]]$.
Soit $\H_2^+$ l'alg{\`e}bre des fonctions {\`a} support compact
dans $\gl(n,\OO_{\vp,2})\cap\GL(n,F_{\vp,2})$ qui sont
bi-$\GL(n,\OO_{\vp,2})$-invariantes. Rappelons que $b:\H^+_2\rta\H^+$
d{\'e}signe l'homomorphisme de changement de base.

Le degr\'e d'extension $r=2$ sera d\'esormais fixe, on notera
$f_\lambda$ et $\phi_\lambda$ pour $f_{2,\lambda}$ et $\phi_{2,\lambda}$.

On a construit dans la section 6, pour toute $n$-partition $\lambda$,
une nouvelle r\'ealisation g\'eom\'etrique $f_\lambda $de la fonction 
$a_{2,\lambda}\in\H^+_2$ de Lusztig.
La fonction $f_{\lambda}$ est d\'efinie comme la trace de
${\tilde\sigma}\circ\Fr$ sur les fibres de
$\A_{\lambda,\vp}\boxtimes\A_{\lambda,\vp}$
au-dessus des points fixes de $\Fr\circ\sigma$. 
La fonction $\phi_{\lambda}=b(f_{\lambda})$ peut alors {\^e}tre 
d{\'e}finie comme la fonction
trace de $\Fr\circ\kappa'$ dans les fibres de
$\A_{\lambda,\vp}*\A_{\lambda,\vp}$ au-dessus
des points fixes de $\Fr$ o\`u $\kappa'$ est l'automorphisme
de $ \A_{\lambda,\vp}*\A_{\lambda,\vp}$ qui se d\'eduit 
de l'isomorphisme de commutativit\'e $\kappa$.

Cette construction reste valable avec les $\tA_\lambda$ ; 
d{\'e}crivons-la bri{\`e}vement.

Soit $\lambda$ une $n$-partition de $d'$. Choisissons un entier $r>d=2d'$.
Notons $\g_{d',r,\vp}$ la fibre de $\g_{d',r}$ au-dessus de
$(\vp^{d'},\vp^r)$. Notons $\tA_{\lambda,\vp}$ la restriction de
$\tA_\lambda[-d](-d/2)$ {\`a} $\g_{d',r,\vp}$. Soient $\A_{\lambda,\vp,1}$
et $\A_{\lambda,\vp,2}$ deux copies de $\A_{\lambda,\vp}$.

L'involution $\sigma:\g_{d',r,\vp}^2\rta\g_{d',r,\vp}^2$ d{\'e}finie par
$\sigma(g,g')=(g',g)$ se rel{\`e}ve en un isomorphisme
$${\tilde\sigma}: \sigma^*(\tA_{\lambda,\vp,1}\boxtimes\tA_{\lambda,\vp,2})
        \rta (\tA_{\lambda,\vp,1}\boxtimes\tA_{\lambda,\vp,2})$$
qui se d{\'e}duit de l'isomorphisme canonique
$\A_{\lambda,\vp,1}{\tilde\rta}\A_{\lambda,\vp,2}$.

Les points fixes de $\sigma\circ\Fr$ dans $\g_{d',r,\vp}^2$
sont de la forme $(g,{\bar g})$ o{\`u} $g\in\g_{d',r,\vp}(k_2)$
et o{\`u} $x\mapsto{\bar x}$ est l'{\'e}l{\'e}ment non trivial du groupe
$\Gal(k_2/k)$.
La trace de ${\tilde\sigma}\circ\Fr$ sur les fibres de
$\tA_{\lambda,\vp,1}\boxtimes\tA_{\lambda,\vp,2}$ au-dessus des points fixes
de $\sigma\circ\Fr$ d{\'e}finit donc une fonction sur $\g_{d',r,\vp}(k_2)$.
Le complexe $\tA_{\lambda,\vp}$ {\'e}tant bi-$G_{r,\vp}$-{\'e}quivariant,
cette fonction est bi-$G_{r,\vp}(k_2)$-invariante.
Elle d{\'e}finit donc une fonction dans $\H^+_2$ ; c'est la fonction 
$f_\lambda$ de la section 6.

La fonction $\phi_\lambda=b(f_\lambda)$ peut {\^e}tre d{\'e}finie par la fonction
trace de $\kappa'\circ\Fr$ sur les fibres de 
$\tA_{\lambda,\vp,1}*\tA_{\lambda,\vp,2}$
au-dessus de $\g_{d',r,\vp}(k)$ o{\`u} $\kappa'$ d{\'e}signe
l'automorphisme
$$\tA_{\lambda,1}*\tA_{\lambda,2}\,{\buildrel \kappa\over\rta}\,
        \tA_{\lambda,2}*\tA_{\lambda,1}
{\buildrel\iota\over\rta}\tA_{\lambda,1}*\tA_{\lambda,2}.$$
o\`u $\kappa$ d\'esigne l'isomorphisme de commutativit\'e
et o\`u $\iota$ sa d\'eduit des isomorphismes identiques entre
$\tA_{\lambda,1}$ et $\tA_{\lambda,2}$.

\subsection{Les fonctions $f'_\lambda,\phi'_\lambda$ et l'identit{\'e} 
            $b'(f'_\lambda)=\phi'_\lambda$}

Soit $S(F_\vp)$ l'ensemble des matrices $s\in\GL(n,F_{2,\vp})$ telles
que $^\t\Fr(s)=s$. Le groupe $\GL(n,F_{2,\vp})$ agit sur $S(F_{\vp})$
par $g.s=\,^\t{\bar g}sg$ o{\`u} $x\mapsto {\bar x}$ est l'{\'e}l{\'e}ment non 
trivial du groupe de Galois $\Gal(F_{2,\vp}/F_\vp)$. 

Notons ${\H'}^+$ l'espace des fonctions {\`a}
support compact dans 
$$S(F_\vp)^+=\gl(n,\OO_{2,\vp})\cap S(F_\vp)$$ 
qui sont
$\GL(n,\OO_{2,\vp})$-invariantes. 
Soit $b':{\H_2}^+\rta{\H'}^+$ l'application lin{\'e}aire qui associe {\`a}
chaque fonction $f\in\H_2^+$ la fonction $\phi'\in{\H'}^+$ d{\'e}finie
par 
$$\phi'(g\,^\t{\bar g})=\int_{H(F_\vp)} f(gh)dh$$
o{\`u} 
$$H(F_\vp)=\{h\in\GL(n,F_{2,\vp})\mid \,^\t\Fr(h)=h^{-1}\}$$
et o{\`u} la mesure de Haar normalis{\'e}e du sous-groupe unitaire $H(F_\vp)$
attribue  {\`a} $H(\OO_\vp)=H(F_\vp)\cap\GL(n,\OO_{2,\vp})$ le volume $1$.
La fonction $\phi'$ est bien d\'efinie puisque toute matrice hermitienne 
$s\in S(F_\vp)$ s'\'ecrit sous la forme $s= g\,^\t{\bar g}$.

En utilisant la proposition 7.4.1 on peut interpr{\'e}ter 
g{\'e}om{\'e}triquement l'application $b'$.

Les points fixes de
$$\displaylines{
(\tau'\times\tau')\circ\sigma\circ\Fr :
\g_{d',r.\vp}^2\rta\g_{d',r.\vp}^2\cr
(g_1,g_2)\rta(\,^\t\Fr(g_2),\,^\t\Fr(g_1))}$$
sont de la forme $(g,\,^\t\Fr(g))$ avec $g\in\g_{d',r,\vp}(k_2)$,
si bien qu'on peut identifier
cet ensemble de points fixes {\`a} $\g_{d',r,\vp}(k_2)$.

Pour toute $n$-partition $\lambda$ de $d'$, la trace de 
$({\tilde\tau}_{\lambda}\times{\tilde\tau}_{\lambda}) 
\circ{\tilde\sigma}\circ\Fr$ dans les fibres de 
$\tA_{\lambda}\boxtimes\tA_{\lambda}$ au-dessus des points fixes de 
$(\tau'\times\tau')\circ\sigma\circ\Fr$ d{\'e}finit donc une fonction sur
$\g_{d',r,\vp}(k_2)$ :
$$f'_\lambda(g)=\Tr(({\tilde\tau}_\lambda\times{\tilde\tau}_\lambda)
\circ{\tilde\sigma}\circ\Fr,
(\tA_\lambda\boxtimes\tA_\lambda)_{(g,\,^\t\Fr(g))}).$$
Le complexe $\tA_\lambda$ {\'e}tant
bi-$G_{r,\vp}$-{\'e}quivariant, cette fonction est 
bi-$\GL(n,\OO_{2,\vp})$-invariante. Notons $f'_\lambda$ l'{\'e}l{\'e}ment de
$\H^+_2$ ainsi d{\'e}fini. On verra dans la section suivante que $f'_\lambda$
n'est qu'une nouvelle r\'ealisation g\'eom\'etrique de la m\^eme
fonction $a_{2,\lambda}$ de Lusztig.

L'ensemble des points fixes de $\tau\circ\Fr$ agissant sur
$\g_{d,r,\vp}$ est l'ensemble
$$\{s\in\g_{d,r,\vp}(k_2)\mid \,^\t\Fr(s)=s\}$$
qui s'identifie {\`a} un sous-ensemble de l'ensemble des classes modulo
$\vp^r$ des {\'e}l{\'e}ments de $S(F_\vp)^+$. 

Rappelons que pour la repr\'esentation induite de $\S_d$
$$\rho=\Ind_{\S_{d'}\times\S_{d'}}^{\S_d}
(\rho_\lambda\boxtimes\rho_\lambda)$$
on a $\A_\rho=\A_\lambda*\A_\lambda$.
La trace de ${\tilde\tau}_\rho\circ\Fr\circ \kappa$ sur les fibres de 
$\tA_\rho$ au-dessus des points fixes de $\tau\circ\Fr$
d{\'e}finit donc une fonction sur $S(F_\vp)^+$ :
$$\phi'_\lambda(s)=\Tr({\tilde\tau}_\rho\circ\kappa'\circ\Fr,
(\tA_\rho)_s).$$
Le complexe $\tA_\rho$ {\'e}tant bi-$G_r$-{\'e}quivariant, cette fonction
appartient {\`a} ${\H'}^+$. Notons-la $\phi'_\lambda$.

\begin{proposition}
Pour toute $n$-partition $\lambda$ de $d'$,
o a $b'(f'_\lambda)=\phi'_\lambda$.
\end{proposition}

\DEMONSTRATION
Soit $s$ un {\'e}l{\'e}ment quelconque $\Fix(\tau\circ\Fr,\g_{d,r,\vp})$.
Il faut d{\'e}montrer que 
$$\displaylines{
\sum_{\scriptstyle g\in\g_{d',r,\vp}(k_2)/H_{r,\vp}(k)\atop
       \scriptstyle g\,^\t\Fr(g)=s}
\Tr(({\tilde\tau}_\lambda\times{\tilde\tau}_\lambda)
     \circ{\tilde\sigma}\circ\Fr, 
     (\tA_\lambda\boxtimes\tA_\lambda)_{(g,\,^\t\Fr(g))})\cr
=\Tr({\tilde\tau}_{\rho}\circ\Fr\circ\kappa,
     (\tA_\rho)_s)}$$
o{\`u} 
$$H_{r,\vp}(k)=\{h\in\GL(\OO_{2,\vp}/\vp^r\OO_{2,\vp})
      \mid\,^\t{\bar h}=h^{-1}\}.$$

L'ensemble o{\`u} s'{\'e}tend la sommation est celui des points fixes de 
$(\tau'\times\tau')\circ\sigma\circ\Fr$ dans la fibre du morphisme
$$m^\alpha:\g_{d',d',r}^{\alpha}\rta\g_{d,r}$$
au-dessus de $s$. Gr{\^a}ce {\`a} la proposition 7.4.1 on sait que
l'endomorphisme de $(\tA_\rho)_s$ induit par
l'endomorphisme de Frobenius tordu
$$({\tilde\tau}_\lambda\times{\tilde\tau}_\lambda)
\circ{\tilde\sigma}\circ\Fr :
((\tau'\times\tau')\times\sigma\times\Fr)^*
(\tA_\lambda\boxtimes^\alpha\tA_\lambda)\rta
\tA_\lambda\boxtimes^\alpha\tA_\lambda$$
est pr{\'e}cis{\'e}ment la restriction de l'endomorphisme
de Frobenius tordu
$${\tilde\tau}_{\rho}\circ\Fr\circ\kappa :
(\tau\circ\Fr)^*(\tA_\rho)\rta\tA_\rho$$
{\`a} $s$.

La proposition r{\'e}sulte donc de la formule des traces de
Grothendieck. 

\subsection{L'identit{\'e} $f_\lambda=f'_\lambda$}

Il revient au m{\^e}me de d{\'e}montrer l'{\'e}nonc{\'e} suivant.

\begin{proposition}
Pour tout $g\in\g_{d,r,\vp}(k_2)$, on a
$$\Tr({\tilde\sigma}\circ\Fr,(\tA_\lambda\boxtimes\tA_\lambda)_{g,\Fr(g)})
=\Tr(({\tilde\tau}_\lambda\times{\tilde\tau}_\lambda)\circ
      {\tilde\sigma}\circ\Fr,
      (\tA_\lambda\boxtimes\tA_\lambda)_{g,\,^\t\Fr(g)})$$
\end{proposition}

\DEMONSTRATION
Puisque les actions de ${\tilde\tau}_\lambda$, de ${\tilde\sigma}$ et
de $\Fr$ commutent dans un sens {\'e}vident, le diagramme
$$\diagramme{
\tA_{\lambda,g}\otimes\tA_{\lambda,\Fr(g)} 
    &\hfl{{\tilde\sigma}\circ\Fr}{}
    &\tA_{\lambda,g}\otimes\tA_{\lambda,\Fr(g)} \cr
\vfl{1\otimes{\tilde\tau}_\lambda}{}
    &
    &\vhfl{}{1\otimes{\tilde\tau}_\lambda^{-1}}\cr
\tA_{\lambda,g}\otimes\tA_{\lambda,\,^\t\Fr(g)}
    &\hfl{}{({\tilde\tau}_\lambda\times{\tilde\tau}_\lambda)\circ
      \,{\tilde\sigma}\circ\,\Fr}
    &\tA_{\lambda,g}\otimes\tA_{\lambda,\,^\t\Fr(g)}\cr
}$$
est commutatif. La trace des deux fl{\`e}ches horizontales sont 
donc {\'e}gales.

\begin{corollaire}
Pour toute $n$-partition $\lambda$, on a 
$a_{2,\lambda}=f_\lambda=f'_\lambda$
et $b'(a_{2,\lambda})=\phi'_\lambda$.
\end{corollaire}

\section{L'interpr{\'e}tation g{\'e}om{\'e}trique 
d'une conjecture de Jacquet et Ye}

\subsection{Enonc{\'e}}

Dans la section pr{\'e}c{\'e}dente, les applications lin{\'e}aires
$$\displaylines{b:{\cal H}^+_2\rta{\cal H}\cr
                b':{\cal H}^+_2\rta{\cal H}'}$$
{\'e}tant interpr{\'e}t{\'e}s g{\'e}om{\'e}triquement, nous sommes maintenant en mesure
de donner une traduction g{\'e}om{\'e}trique du lemme fondamental de Jacquet
et Ye et de le d{\'e}montrer dans les m{\^e}me lignes que \cite{Ngo1}.

Notons $A$ le sous-groupe diagonal de $\GL(n)$ et
$N$ son sous-groupe des matrices
triangulaires sup{\'e}rieures unipotentes. Pour chaque
$$\alpha=(\alpha_1,\ldots,\alpha_{n-1})\in(k^\times)^{n-1}$$
notons $\theta_\alpha:N(F_\vp)\rta\Ql^\times$ le caract{\`e}re
$$\theta_\alpha(x)=\psi\bigl(\sum_{i=1}^{n-1}
\alpha_{i,i+1}\res(x_{i,i+1})\bigr)$$
o{\`u} $\psi:k\rta\Ql^\times$ est un caract{\`e}re additif non trivial
qu'on fixe une fois pour toutes.

Pour toute fonction $\phi\in{\cal H}^+$,
toute matrice diagonale $a\in A(F_\vp)$, pour
tout $\alpha\in(k^\times)^{n-1}$, posons
$$I_\vp(a,\alpha,\phi)=\int_{N(F_\vp)\times N(F_\vp)}
        \phi(\,^\t x_1ax_2)\theta_\alpha(x_1)
        \theta_\alpha(x_2)\d x_1\d x_2$$
o{\`u} la mesure de Haar normalis{\'e} $\d x$
de $N(F_\vp)$ attribue {\`a}
$N({\cal O}_\vp)$ le volume $1$.

Cette int{\'e}grale intervient comme une int{\'e}grale orbitale dans
une formule des traces relative de Jacquet.
Il s'agit d'une int{\'e}grale de Kloosterman
si $\phi$ est la fonction caract{\'e}ristique de $GL(n,{\cal O}_\vp)$,

Pour tout $\alpha:(k^\times)^{n-1}$, notons
$\theta'_\alpha: N(F_{2,\vp})\rta\Ql^\times$ le caract{\`e}re d{\'e}fini par
$$\theta'_\alpha(x)=\psi\bigl(\sum_{i=1}^{n-1}
\alpha_i\res(x_{i,i+1}+{\bar x}_{i,i+1})\bigr).$$

Pour toute fonction $\phi'\in{\cal H'}^+$
pour toute matrice diagonale $a\in A(F_\vp)$
et pour tout $\alpha\in(k^\times)^{n-1}$, posons
$$J_\vp(a,\alpha,\phi')=\int_{N(F_{2,\vp})}
\phi'(\,^\t{\bar x}ax)\theta'_\alpha(x)\d x$$
o{\`u} la mesure de Haar normalis{\'e}e $\d x$ de
$N(F_{2,\vp})$ attribue {\`a} $N({\cal O}_{2,\vp})$
le volume $1$.

\begin{theoreme}
Pour toute fonction $f\in{\cal H}^+_2$
pour 
$$a=\diag(a_1,a_1^{-1}a_2,\ldots,a_{n-1}^{-1}a_n)\in A(F_\vp)$$
pour tout $\alpha\in(k^\times)^{n-1}$, on a
$$I_\vp(a,\alpha,b(f))=(-1)^{\v_\vp(a_1\ldots  a_{n-1})}
J_\vp(a,\alpha,b'(f)).$$
\end{theoreme}

Cet {\'e}nonc{\'e} joue le r{\^o}le d'un lemme fondamental dans
une formule des traces relative. Jacquet et Ye l'ont conjectur{\'e} dans
\cite{JY} sans hypoth{\`e}se sur la caract{\'e}ristique et l'ont
d{\'e}montr{\'e} pour $n=2$ et $n=3$.
Le cas o{\`u} $f$ est l'unit{\'e} de l'alg{\`e}bre de Hecke a
{\'e}t{\'e}  d{\'e}montr{\'e} dans \cite{Ngo1}.

\subsection{Sommes locales}

Pour tout $a\in A(F_\vp)$, on a d{\'e}fini dans \cite{Ngo1} 
un triplet $({\cal X}_\vp(a),h,\tau)$ o{\`u} 
${\cal X}_\vp(a)$ est un sch{\'e}ma de type fini sur $k$ tel que 
$${\cal X}_\vp(a)(k)=\{(x,x')\in(N(F_\vp)/N(\OO_\vp))^2
   \mid \,^\t xax'\in\gl(n,\OO_\vp)\}$$
o{\`u} le morphisme
$$h:{\cal X}_\vp(a)\times\Gm^{n-1}\rta\Ga$$
induit sur les $k$-points l'application
$$h(x,x',\alpha)=\sum_{i=1}^{n-1}\alpha_i\res(x_{i,i+1}+x'_{i,i+1})\d\vp$$
et o{\`u} l'involution $\tau:{\cal X}_\vp(a)\rta{\cal X}_\vp(a)$
est d{\'e}finie par $\tau(x,x')=(x',x)$.

On peut {\'e}crire la matrice $a$ sous la forme
$$a=\diag(a_1,a_1^{-1}a_2,\ldots,a_{n-1}^{-1}a_n).$$
Pour que le sch{\'e}ma ${\cal X}_{\vp}(a)$ ne soit pas vide, il est
n{\'e}cessaire que $a_1,\ldots,a_n\in\OO_\vp$ (\cite{Ngo1}).

Pour chaque entier $r\in\NN$ tel que
$$r>\v_\vp(a_1)+\cdots+\v_\vp(a_n)$$
soit ${\cal X}_{\vp,r}(a)$ le sch{\'e}ma de type fini sur $k$ dont
l'ensemble des $k$-points est l'ensemble
$$\{g\in\gl(n,\OO_\vp/\vp^r\OO_\vp)\mid \Delta_i(g)\cong a_i\ \mod \vp^r
    {\rm\  pour\ tout\ }i=1,\ldots,n \}$$
o{\`u} $\Delta_i(g)$ est le d{\'e}terminant de la sous-matrice de $g$ faite
des $i$ premi{\`e}res lignes et des $i$ premi{\`e}res colonnes. 

La limite projective 
$${\cal X}_{\vp,\infty}(a)(k)=\limpro{\cal X}_{\vp,r}(a)(k)$$
s'identifie {\`a} l'ensemble  
$$\{g\in\gl(n,\OO_\vp)\mid \Delta_i(g)=a_i
    {\rm\  pour\ tout\ }i=1,\ldots,n \}.$$
Tout {\'e}l{\'e}ment $g\in{\cal X}_{\vp,\infty}(a)(k)$ s'{\'e}crit
alors  de mani{\`e}re
unique sous la forme $g=\,^\t nan'$ avec $n,n'\in N(F_\vp)$ si bien 
qu'on a une application 
$$p_{\infty}(k):{\cal X}_{\vp,\infty}(a)(k)\rta{\cal X}_{\vp}(a)(k).$$

\begin{proposition}
\begin{enumerate}
\item
Sous la condition 
$$r> \v_\vp(a_1)+\cdots+\v_\vp(a_n)$$
on a un morphisme
$$p_r:{\cal X}_{\vp,r}(a)\rta{\cal X}_{\vp}(a)$$
tel que l'application $p_\infty(k)$ se factorise {\`a} travers $p_r(k)$.
\item
De plus, $p_r$ est un compos{\'e} de fibrations vectorielles.
\end{enumerate}
\end{proposition}

\DEMONSTRATION
\begin{enumerate}
\item
Soit $g\in\gl(n,\OO_\vp)$ tel que $\Delta_i(g)=a_i$ pour tout $i=1,\ldots,n$. 
Notons $g_i$ la sous-matrice de $g$  faite des $i$ premi{\`e}res lignes et 
des $i$ premi{\`e}res colonnes. On a pour tout $i=1,\ldots,n-1$
$$g_{i+1}=\pmatrix{g_{i} & y'_i \cr ^\t y_i & z_i }$$
o{\`u} $y_i$ et $y'_i$ sont des vecteurs colonnes de taille $i\times 1$
et o{\`u} $z_i\in k$. On a donc 
$$g_{i+1}=\,^\t x'_i\pmatrix{g_i & 0\cr 0& a_i^{-1}z_i} x_i$$
o{\`u} 
$$x_i=\pmatrix{\Id_{i} & g_i^{-1}y_i\cr 0 & 1}$$
et o{\`u} 
$$x'_i=\pmatrix{\Id_{i} & ^\t g_i^{-1}y'_i\cr 0  & 1}.$$
Si on pose 
$$\displaylines{x=x_1 x_2\ldots x_{n-1}\cr 
x'=x'_1 x'_2\ldots x'_{n-1}}$$
alors on a 
$$g=\,^\t x\,\diag(a_1,a_1^{-1}a_2,\ldots,a_{n-1}^{-1}a_n)\,x'.$$

Maintenant, si les $y_i,y'_i$ sont d{\'e}termin{\'e}s  modulo 
$\vp^r$ avec 
$$r>\v_\vp(a_1\ldots a_n),$$ 
en utilisant les relations de 
commutation, on montre que $x$ et $x'$ sont d{\'e}termin{\'e}s modulo 
$N(\OO_\vp)$. 

L'application $p_r(k):{\cal X}_{\vp,r}(a)(k)\rta {\cal X}_{\vp}(a)$
ainsi d{\'e}finie provient d'un mor\-phisme
$p_r:{\cal X}_{\vp,r}(a)\rta{\cal X}_{\vp}(a)$.

\item
On proc{\`e}de par r{\'e}currence sur $n$.
L'assertion est triviale pour $n=1$.
Suppo\-sons qu'elle est vraie pour $n-1$. 

Notons $a'$ la matrice
$$a'=\diag(a_1,a_1^{-1}a_2,\ldots,a_{n-2}^{-1}a_{n-1})\in\GL(n-1,F_\vp).$$
Par r{\'e}currence, on peut supposer que 
$${\cal X}_{\vp,r}(a')\times_{{\cal X}_{\vp}(a')}{\cal X}_{\vp}(a)
\rta {\cal X}_\vp(a)$$
peut se factoriser en fibrations vectorielles.

La donn{\'e}e d'un $k$-point de 
${\cal X}_{\vp,r}(a')\times_{{\cal X}_{\vp}(a')}{\cal X}_{\vp}(a)$
est {\'e}quivalente {\`a} la donn{\'e}e de 
$$\pmatrix{g' & y'_{n-1}\cr y_{n-1} & z_{n-1}}$$
o{\`u} 
$$g'\in{\cal X}_{\vp,r}(a')\subset\gl(n,\OO_\vp/\vp^r)$$
o{\`u} $y'_{n-1}$ (resp. $y_{n-1}$) est dans $\OO_{\vp}^{n-1}$ et d{\'e}termin{\'e}
modulo $g'\OO_\vp^{n-1}$ (resp. $^\t g'\OO_\vp^{n-1})$)
et o{\`u} $z_{n-1}\in\OO_\vp$ est d{\'e}termin{\'e} modulo $a_{n-1}\OO_\vp$.

Il est clair que 
$${\cal X}_{\vp,r}(a)\rta
{\cal X}_{\vp,r}(a')\times_{{\cal X}_{\vp}(a')}{\cal X}_{\vp}(a)$$
est un fibr{\'e} vectoriel, d'o{\`u} l'assertion.
\end{enumerate}

Le sch{\'e}ma ${\cal X}_{\vp,r}$ est naturellement une partie localement
ferm{\'e}e de $\g_{d_n,r,\vp}$ o{\`u} $d_n=\v_\vp(a_n)$.
Pour toute repr{\'e}sentation $\rho$ de $\S_{d_n}$, on d{\'e}signe encore par
${\tA}_\rho$ sa restriction {\`a} ${\cal X}_{\vp,r}$.

\begin{proposition}
Posons
$${\dot\A}_{\rho,\vp}=\RR p_{r,!}\tA_{\rho,\vp}[2d_r](d_r)$$
o{\`u} $d_r$ est la dimension relative de $p_r$.
Ce morphisme {\'e}tant lisse on a un isomorphisme
$$\RR p_r^!{\dot\A}_{\rho,\vp}{\simeq} p_r^*{\dot\A}_{\rho,\vp}[2d_r](d_r).$$
La fl{\`e}che d'adjonction
$$p_r^*{\dot\A}_\rho{\simeq}
\RR p_r^!{\dot\A}_\rho[-2d_r](-d_r)\rta\tA_\rho$$
est alors un isomorphisme. 
\end{proposition}

\DEMONSTRATION
Point par point, l'isomorphisme r{\'e}sulte de ce que 
les fibres g{\'e}om{\'e}triques de $p_r$ sont isomorphes {\`a} l'espace affine de
dimension de $d_r$ et qu'elles sont incluses dans les orbites de 
$G_{r,\vp}\times G_{r,\vp}$ sur $\g_{d,r}$ au-dessus desquelles
le complexe $\tA_\rho$ est constant. $\carre$

La fl{\`e}che ${\tilde\tau}_\rho:\tau_{d_n}^*\tA_\rho\rta\tA_\rho$
induit une fl{\`e}che 
${\dot\tau}_\rho:\tau^*{\dot\A}_\rho\rta{\dot\A}_\rho$.
On retrouve ${\tilde\tau}_\rho$ comme l'image r{\'e}ciproque de 
${\dot\tau}_\rho$.

\begin{corollaire}
Pour toute $n$-partition $\lambda$ de $d'=d_n/2$, on a 
$$\displaylines{
I_\vp(a,\alpha,\phi_\lambda)=\Tr(\Fr\circ\kappa,
        \RR\Gamma_c({\cal X}_\vp(a)\otimes_k{\bar k},
        {\dot\A}_\rho\otimes h_\alpha^*\L_\psi));\cr
J_\vp(a,\alpha,\phi'_\lambda)=\Tr(\Fr\circ\kappa\circ{\dot\tau}_\rho,
        \RR\Gamma_c({\cal X}_\vp(a)\otimes_k{\bar k},
        {\dot\A}_\rho\otimes h_\alpha^*\L_\psi)).
}$$
o{\`u} $\rho$ est la repr{\'e}sentation induite 
$$\Ind_{\S_{d'}\times\S_{d'}}^{\S_{2d'}}(\rho_\lambda\times\rho_\lambda)$$
et o{\`u} $\L_\psi$ est le faisceau d'Artin-Schreier sur $\Ga$ associ{\'e}
{\`a} $\psi:k\rta\Ql^\times$.

Dans le cas o{\`u} $|\lambda|\not=d'$, les int{\'e}grales 
pr{\'e}c{\'e}dentes sont nulles.
\end{corollaire}

\DEMONSTRATION
Compte tenu de l'interpr{\'e}tation g{\'e}om{\'e}trique des fonctions 
$\phi_\lambda$ et $\phi'_\lambda$, c'est la formule des traces 
de Grothendieck. $\square$

\bigskip

\noindent{\sc Th\'eor\`eme 4a}
{\it Pour toute repr{\'e}sentation $\rho$ de $\S_d$, pour tout 
$\alpha\in{\bar k}^{n-1}$, l'involution ${\dot\tau}_\rho$ agit dans 
$\RR\Gamma_c({\cal X}_\vp(a)\otimes_k{\bar k},
{\dot\A}_\rho\otimes h_\alpha^*\L_\psi)$ 
comme la multiplication par
$(-1)^{\v_\vp(a_1\ldots a_{n-1})}$.}

\subsection{Sommes globales}

Soit ${\und d}=(d_i)_{i=1}^n\in{\Bbb N}^n$. Notons 
$Q_{\und d}=\prod_{i=1}^n Q_{d_i}$.
On a d{\'e}fini dans \cite{Ngo1} un quadruplet 
$({\cal X}_{\und d},f_{\und d},h_{\und d},\tau_{\und d})$
o{\`u} le morphisme de type fini
$$f_{\und d}:{\cal X}_{\und d}\rta Q_{\und d}$$ 
tel que pour tout $a=(a_i)_{i=1}^n\in Q_{\und d}(k)$, on a
$${\cal X}(a)(k)=\{(x,x')\in(N(F)/N(\OO))^2\mid \,^\t xax'\in\gl(n,\OO)\}$$
o{\`u} le morphisme 
$$h:{\cal X}_{\und d}\times\Gm^{n-1}\rta\Ga$$
induit au niveau des $k$-points l'application
$$h(x,x',\alpha)=\sum_{i=1}^{n-1}\sum_{v\not=\infty}
        \res_v(x_{i,i+1}+x'_{i,i+1})\d\vp$$
et o{\`u} $\tau_{\und d}:{\cal X}_{\und d}\rta{\cal X}_{\und d}$ est 
l'involution $\tau_{\und d}(x,x')=(x',x)$.

Soit ${\tilde{\cal X}}_{\und d}$ le $Q_{\und d}$-sch{\'e}ma dont l'ensemble 
des $k$-points au-dessus de $a=(a_i)_{i=1}^n\in Q_{\und d}(k)$ est 
l'ensemble
$$\{g\in\gl(n,\OO/(R))\mid\Delta_i(g)=a_i\,\mod R,\ i=1,\ldots,n\}$$
o{\`u} $R=(a_1\ldots a_n)^2$. Il nous faut choisir $R$ strictement 
divisible par $a_1\ldots a_n$.

Les assertions suivantes se d{\'e}montrent exactement que leurs analogues 
locaux.

\begin{proposition}
On a un morphisme $p_{\und d}:{\tilde{\cal X}}_{\und d}\rta{\cal X}_{\und d}$
qui, au niveau des $k$-points, envoie la r{\'e}duction modulo $R$ de
$$g=\,^\t xax'\in\gl(n,\OO)$$
sur $(xN(\OO),x'N(\OO))\in(N(F)/N(\OO))^2$.
De plus, $p_{\und d}$ peut se factoriser en fibrations vectorielles.
\end{proposition}

\begin{proposition}
Posons
$${\dot\A}_{\rho}=\RR p_{r,!}\tA_{\rho}[2d_r](d_r)$$
o{\`u} $d_r$ est la dimension relative de $p_r$.
Ce morphisme {\'e}tant lisse on a un isomorphisme
$$\RR p_r^!{\dot\A}_{\rho}{\simeq} p_r^*{\dot\A}_{\rho}[2d_r](d_r).$$
La fl{\`e}che d'adjonction
$$p_r^*{\dot\A}_\rho{\simeq}
\RR p_r^!{\dot\A}_\rho[-2d_r](-d_r)\rta\tA_\rho$$
est alors un isomorphisme. 

L'involution ${\tilde\tau}_\rho:\tau^*\tA_\rho\rta\tA_\rho$
descend aussi en une involution 
$${\dot\tau}_\rho:\tau^*_{\und d}{\dot\A}_\rho\rta{\dot\A}_\rho.$$
\end{proposition}

Voici la variante globale du th{\'e}or{\`e}me {\sc 4a}.

\bigskip
\noindent{\sc Th\'eor\`eme 4b}
{\it L'involution ${\dot\tau}_{\rho}$ agit sur le complexe de faisceaux
$$\RR(f_{\und d}\times\Id_{\Gm^{n-1}})_!
({\dot\A}_\rho\otimes h_{\und d}^*\L_\psi)$$
comme la multiplication par $(-1)^{d_1+\cdots+d_{n-1}}$.}

\bigskip 

Cet {\'e}nonc{\'e} se d{\'e}duit de son analogue local en utilisant 
la formule de multiplicativit{\'e} {\'e}nonc{\'e}e dans la section 
qui suit. 

Toutefois, comme dans \cite{Ngo1}, on commencera par d{\'e}montrer
l'{\'e}nonc{\'e} global dans le cas tr{\`e}s particulier 
${\und d}=(1,2,\ldots,n)$ puis en d{\'e}duire l'{\'e}nonc{\'e} local 
en utilisant la formule de multiplicativit{\'e}.

\subsection{Compatibilit{\'e}}

Pour tout id{\'e}al maximal $v$ de ${\bar\OO}=\OO\otimes_k{\bar k}$, 
pour tout $a\in Q_{\und d}({\bar k})$
et $\alpha\in({\bar k}^{\times})^{n-1}$,
on a d{\'e}fini dans \cite{Ngo1} un triplet
$({\cal X}_{{v}}(a),h_{\alpha,{v}},\tau_{{v}})$
o{\`u} ${\cal X}_v(a)$ est un sch{\'e}ma de type fini sur ${\bar k}$
dont l'ensemble des ${\bar k}$-points est l'ensemble
$${\cal X}_v(a)({\bar k})=
\{x,x'\in N({\bar F}_v)/N({\bar\OO}_v)
\mid \,^\t xax'\in\gl(n,\OO_v)\}$$ 
o{\`u} $h_{\alpha,v}:{\cal X}_v(a)\rta{\Ga}_{\bar k}$ est le morphisme
$$h_{\alpha,v}(x,x')=\sum_{i=1}^{n-1}\res_v(x_{i,i+1}+x'_{i,i+1})$$
o{\`u} $\tau_v$ est le morphisme $\tau(x,x')=(x',x)$.
On renvoie {\`a} \cite{Ngo1} pour la d{\'e}monstration de l'{\'e}nonc{\'e} suivant.

\begin{proposition}
Pour tout $a\in Q_{\und d}({\bar k})$, on a un isomorphisme
$${\cal X}_{\und d}(a)=\prod_{v|a_1\ldots a_{n-1}}
                {\cal X}_{{v}}(a)$$
via lequel on a
$$h_{\und d}(a)=\sum_{v|a_1\ldots a_{n-1}}
                h_{{v}}$$
et
$$\tau_{\und d}(a)=\prod_{v|a_1\ldots a_{n-1}}
                \tau_{{v}}.$$
Pour $v$ ne divisant pas $a_1\ldots a_{n-1}$, ${\cal X}_v(a)$
est r{\'e}duit {\`a} un point.
\end{proposition}

Soit maintenant $\rho$ une repr{\'e}sentation de $\S_{d_n}$ 
o{\`u} $d_n=\v_\vp(a_n)$. Notons $d_{n,v}=\v_{{v}}(a_n)$.
On a $d_n=\sum_{v|a_n}d_{n,\lambda}$.

On d{\'e}compose la restriction de $\rho$ {\`a} 
$\prod_{v|a_n}\S_{d_{n,v}}$ en somme de repr{\'e}sentations 
irr{\'e}ductibles. Celles-ci sont tous  de la forme 
$\bigotimes_{v|a_n}\rho_{v}$ o{\`u} $\rho_v$ est une repr{\'e}sentation 
irr{\'e}ductible de $\S_{d_{n,v}}$.

\begin{proposition}
Supposons que
$$\Res_{\prod_{v|a_n}\S_{d_{n,v}}}^{\S_{d_n}}\rho=
\bigoplus_i\bigotimes_{v|a_n}\rho_{i,v}.$$
Pour tout $a\in Q_{\und d}({\bar k})$, on a
$$\displaylines{
\RR\Gamma_c({\cal X}_{\und d}(a),
{\dot\A}_\rho\otimes h_{\und d}^*\L_\psi)\cr
=\bigoplus_i\bigotimes_{v|a_1\ldots a_{n}}
\RR\Gamma_c({\cal X}_v(a),
{\dot\A}_{\rho_{i,v}}\otimes h_{v}^*\L_\psi)}$$
o{\`u} pour $v$ ne divisant pas $a_n$, on prend pour
$\rho_{i,v}$ la repr{\'e}sentation triviale.
De plus, l'action ${\dot\tau}_\rho$ sur le premier membre est {\'e}gale
{\`a} celle de 
$\bigoplus_i\bigotimes_{v|a_1\ldots a_{n}}{\dot\tau}_{\rho_{i,v}}$
sur le second membre.
\end{proposition}

Dans le cas o\`u la repr\'esentation $\rho$ est triviale, on retrouve 
la variante cohomologique de la formule de multiplicativit\'e
pour les int\'egrales de Kloosterman (corollaire 3.2.3 \cite{Ngo1}).

\DEMONSTRATION
Compte tenu de la proposition pr{\'e}c{\'e}dente, il suffit de d{\'e}montrer
que via l'isomorphisme
$${\cal X}_{\und d}(a)=\prod_{v|a_1\ldots a_{n}}{\cal X}_v(a)$$
on a un isomorphisme
$${\dot\A}_\rho=\bigoplus_{i}\bigotimes_{v|a_1\ldots a_{n}}\tA_{\rho_{i,v}}.$$
Via 9.3.2, on se ram{\`e}ne {\`a} la proposition 7.5.1. $\square$

\begin{proposition}
Si $v$ divise $a_1\ldots a_{n-1}$ avec la multiplicit{\'e} $1$, et si $v$
ne divise pas $a_n$ alors dans la d{\'e}composition pr{\'e}c{\'e}dente
$$\RR\Gamma_c({\cal X}_v(a)\otimes_k{\bar k},
{\dot\A}_{\rho_{i,v}}\otimes h_{v}^*\L_\psi)$$
est un $\Ql$-espace vectoriel de rang $2$ plac{\'e} en degr{\'e} $1$
dans lequel $\tau_{v}$ agit comme $-1$.
\end{proposition}

\DEMONSTRATION
La repr{\'e}sentation $\rho_{i,v}$ est triviale du fait que $v$ ne divise
pas $a_n$. On se ram{\`e}ne donc {\`a} l'{\'e}nonc{\'e} 2.4 dans \cite{Ngo1}
et donc finalement {\`a} un th{\'e}or{\`e}me de Deligne sur les sommes de 
Kloosterman classiques (\cite{Del}). $\square$

\section{D{\'e}monstration du th{\'e}or{\`e}me 5}
\subsection{L'ouvert $U_{\und d}$}

Soit $U_{\und d}$ l'ouvert dense de $Q_{\und d}$ form{\'e} des suites
$a=(a_i)_{i=1}^n$ telles que le polyn{\^o}me produit $\prod_{i=1}^n a_i$
est s{\'e}parable.

\begin{proposition} 
On a un isomorphisme
$$\displaylines{\RR(f_{\und d}\times\Id_{\Gm^{n-1}})_!
({\dot\A}_\rho\otimes h_{\und d}^*\L_\psi)|_{U_n\times\Gm^{n-1}}\cr
=\RR(f_{\und d}\times\Id_{\Gm^{n-1}})_!\, 
h_{\und d}^*\L_\psi|_{U_n\times\Gm^{n-1}}
\otimes\pr_{U_{d_n}}^*\L_\rho}$$
o{\`u} $\pr_{U_{d_n}}:U_{\und d}\rta U_{d_n}$ est la projection de $U_{\und d}$
sur l'ouvert $U_{d_n}\subset Q_{d_n}$ des polyn{\^o}mes unitaires s{\'e}parables
$a_n$ de degr{\'e} $d_n$ et o{\`u} $\L_\rho$ est le syst{\`e}me local associ{\'e}
{\`a} la repr{\'e}sentation $\rho$ de groupe de Galois $\S_{d_n}$ du 
rev{\^e}tement {\'e}tale galoisien $\AA^{d_n}_{rss}\rta U_{d_n}$.
De plus, ${\tilde\tau}_{\rho}$ agit dans 
$\RR(f_{\und d}\times\Id_{\Gm^{n-1}})_! h_{\und d}^*\L_\psi
|_{U_n\times\Gm^{n-1}}
\otimes\pr_{U_{d_n}}^*\L_\rho$ comme ${\tilde\tau}\otimes\Id_{\L_\rho}$.
\end{proposition}

\DEMONSTRATION
Par d{\'e}finition de $\A_\rho$, on sait que la restriction de ${\dot\A}_\rho$
{\`a} l'ouvert ${\cal X}_{\und d}\times_{Q_{\und d}}U_{\und d}$
est l'image r{\'e}ciproque de $\L_\rho$. 
La proposition r{\'e}sulte donc de la formule de projection. $\carre$

L'{\'e}nonc{\'e} suivant, extrait de \cite{Ngo1} se d{\'e}duit de la formule
de multiplicativit{\'e}.

\begin{proposition} 
$\RR(f_{\und d}\times\Id_{\Gm^{n-1}})_! h_{\und d}^*\L_\psi
|_{U_{\und d}\times\Gm^{n-1}}$
est un syst{\`e}me local de rang $2^{d_1+\cdots +d_{n-1}}$ plac{\'e} en degr{\'e}
$d_1+\cdots+\d_{n-1}$ dans lequel $\tilde\tau$ agit comme 
$(-1)^{d_1+\cdots+d_{n-1}}$.
\end{proposition}

\begin{corollaire}
L'involution  ${\tilde\tau}_{\rho}$ agit dans 
$$\RR(f_{\und d}\times\Id_{\Gm^{n-1}})_!
({\dot\A}_\rho\otimes h_{\und d}^*\L_\psi)|_{U_n\times\Gm^{n-1}}$$
comme la multiplication par $(-1)^{d_1+\cdots+d_{n-1}}$.
\end{corollaire}

\subsection{Le cas ${\und d}=(1,2,\ldots,n)$}

On a d{\'e}montr{\'e} dans \cite{Ngo1} que pour toute suite $a=(a_i)_{i=1}^n$
dont chaque membre $a_i$ est un polyn{\^o}me unitaire de degr{\'e} $i$,
pour tous $x,x'\in N(F)$ tels que $\,^\t xax'\in\gl(n,\OO)$,
il existe une unique matrice de la forme $\gamma+\Id_n\vp$ avec 
$\gamma\in\gl(n,k)$ telle que 
$$\gamma+\Id_n\vp\in\,^\t N(\OO)\,^\t xax' N(\OO).$$

L'application $(x,x')\mapsto\gamma$ d{\'e}finit une section 
$$\iota:\gl(n)={\cal X}_{\und d}\rta{\tilde{\cal X}}_{\und d}$$
via laquelle on a ${\dot\A}_\rho=\iota^*{\tilde A}_\rho$.
En utilisant le lemme 2.3.1, on sait alors que ${\dot A}_\rho$
est un faisceau pervers {\'e}quivariant pour l'action adjointe de 
$\GL(n)$ et qui est isomorphe au prolongement interm{\'e}diaire 
de sa restriction {\`a} l'ouvert $\gl(n)_{rss}$ form{\'e} des {\'e}l{\'e}ments
r{\'e}guliers semi-simples.

Identifions $\GL(n-1)$ au sous-groupe 
$$\diag(\GL(n-1),1)\subset\GL(n).$$
L'{\'e}nonc{\'e} suivant se d{\'e}duit imm{\'e}diatement de la proposition
5.2.2 de \cite{Ngo1}.

\begin{proposition}
Si $K$ est un faisceau pervers sur $\gl(n)$ qui est 
$\GL(n-1)$-{\'e}quivariant et est isomorphe {\`a} son prolongement 
interm{\'e}diaire de sa restriction {\`a} l'ouvert $\gl(n)_{rss}$,
le complexe de faisceaux 
$$\RR(f_{\und d}\times\Id_{\Gm^{n-1}})_!
(K\otimes h_{\und d}^*\L_\psi)$$
est {\`a} d{\'e}calage pr{\`e}s un faisceau pervers, prolongement interm{\'e}diaire
de sa restriction {\`a} l'ouvert $U_{\und d}\times\Gm^{n-1}$.
\end{proposition} 

Compte tenu de ce r\'esultat et du corollaire 10.1, on obtient
l'{\'e}nonc{\'e} suivant.

\begin{corollaire}
Lorsque ${\und d}=(1,2,\ldots,n)$, ${\dot\tau}_\rho$ agit dans 
$$\RR(f_{\und d}\times\Id_{\Gm^{n-1}})_!
({\dot\A}_\rho\otimes h_{\und d}^*\L_\psi)$$
comme la multiplication par $(-1)^{1+2+\cdots+(n-1)}$.
\end{corollaire}

\subsection{Augmenter $n$}

Pour d{\'e}duire du corollaire 10.2.2 le th{\'e}or{\`e}me {\sc 4a}, 
et donc aussi son analogue global {\sc 4b}
l'astuce consiste {\`a} remplacer
$n$ par un entier assez grand.

\begin{lemme}
\begin{enumerate}
\item Pour tout $a_\vp\in A(F_\vp)$, pour $\alpha=1$ 
pour tout $m\in\NN$, les donn{\'e}es
$$\displaylines{({\cal X}_\vp(a_\vp),h_\alpha,\tau,{\dot\A}_\rho,{\dot\tau}_\rho)\cr 
({\cal X}_\vp(\diag(\Id_{m},a_\vp)),h_\alpha,\tau,{\dot\A}_\rho,{\dot\tau}_\rho)}$$ 
sont isomorphes.
\item Pour tous $a_\vp,a'_\vp\in A(F_\vp)$ tels que pour tout 
$i=1,\ldots,n$, 
$$a_i\cong a'_i\,\mod\,\vp^r$$ 
pour 
$r=\v_\vp(a_1\ldots a_{n})$
les donn{\'e}es
$$\displaylines{({\cal X}_\vp(a_\vp),h_\alpha,\tau,{\dot\A}_\rho,{\dot\tau}_\rho)\cr
({\cal X}_\vp(a'_\vp),h_\alpha,\tau,{\dot\A}_\rho,{\dot\tau}_\rho)}$$ 
sont isomorphes.
\end{enumerate}
\end{lemme}

\DEMONSTRATION
\begin{enumerate}
\item On a construit dans \cite{Ngo1}, un isomorphisme 
$$({\cal X}_\vp(a_\vp),h_\alpha,\tau){\tilde\rta}
({\cal X}_\vp(\diag(\Id_{m},a_\vp)),h_\alpha,\tau).$$ 
L'isomorphisme entre les ${\dot\A}_\rho$ r{\'e}sulte du corollaire 2.2.3.
\item La construction de l'isomorphisme
$$({\cal X}_\vp(a_\vp),h_\alpha,\tau){\tilde\rta}
({\cal X}_\vp(a'_\vp),h_\alpha,\tau)$$
dans \cite{Ngo1} fournit {\'e}galement un isomorphisme pour 
les ${\dot\A}_\rho$. $\carre$ 
\end{enumerate}

Le lemme suivant est extrait de \cite{Ngo1}.

\begin{lemme}
Pour tout $b_\vp\in A(F_\vp)$, pour un entier $m\in\NN$ assez grand,
il existe une suite $a=(a_i)_{i=1}^{m+n}$ dont chaque membre $a_i$
est un polyn{\^o}me unitaire de degr{\'e} $i$ {\`a} coefficients dans ${\bar k}$,
qui satisfait deux conditions suivantes
\begin{itemize}
\item
si on {\'e}crit $b=\diag(\Id_m,b_\vp)$ sous la forme
$$b=(b_{1},{b}_{1}^{-1}b_{2},\ldots,
{b}_{m+n-1}^{-1}b_{m+n})$$
alors on a $a_i\cong b_{i}\,\mod \vp^r$ pour tout $i=1,2\ldots,m+n$
et pour $r=\v_\vp(b_1\ldots b_{m+n})$ ;
\item
pour tout id{\'e}al maximal $v$ de $\OO\otimes_k{\bar k}$ 
diff{\'e}rente de de $\vp$,
$v$ divise $a_1\ldots a_{m+n}$ avec au plus une multiplicit{\'e} $1$. 
\end{itemize}
\end{lemme}

\noindent{\it Fin de la d{\'e}monstration du th{\'e}or{\`e}me 5.}
Soit $b$ un {\'e}l{\'e}ment quelconque de $A(F_\vp)$. Choisissons une
suite $a=(a_i)_{i=1}^{m+n}$ comme dans le lemme pr{\'e}c{\'e}dent afin que 
les donn{\'e}es  
$$\displaylines{
({\cal X}_\vp(b_\vp),h_\alpha,\tau,{\dot\A}_\rho,{\dot\tau}_\rho)\cr
({\cal X}_\vp(a),h_\alpha,\tau,{\dot\A}_\rho,{\dot\tau}_\rho)}$$ 
soient isomorphes.

Soit $\rho$ une repr{\'e}sentation de $\S_{\v_\vp(a_{m+n})}$ 
qu'on peut supposer irr{\'e}ductible. Notons 
$$\rho'=\Ind_{\S_{\v_\vp(a_{m+n})}}
^{\S_{m+n}}\rho$$
o{\`u} on a identifi{\'e} $\S_{\v_\vp(a_{m+n})}$ au sous-groupe
$$\S_{\v_\vp(a_{m+n})}\times\S_1^{m+n-\v_\vp(a_{m+n})}$$
de $\S_{m+n}$.

Dans la formule 10.4.2,  
$$\displaylines{
\RR\Gamma_c({\cal X}_{\und d}(a)\otimes_k{\bar k},
{\dot\A}_{\rho'}\otimes h_{\und d}^*\L_\psi)\cr
=\bigoplus_i\bigotimes_{v|a_1\ldots a_n}
\RR\Gamma_c({\cal X}_v(a)\otimes_k{\bar k},
{\dot\A}_{\rho_{i,v}}\otimes h_{v}^*\L_\psi)}$$
$\rho$ est une des repre\'sentations irr{\'e}ductibles $\rho_{i,\vp}$.
Sachant que ${\dot\tau}_{\rho'}$ agit dans le membre de gauche
comme la multiplication par $(-1)^{1+2+\cdots+(m+n-1)}$, il agit 
de la m{\^e}me mani{\`e}re dans tous les facteurs directs du membre de 
droite. Or dans chacun de ces facteurs directs, d'apr{\`e}s 10.4.3,
les complexes locaux en $v\not=\vp$ sont des espaces vectoriels de 
rang $2$ plac{\'e}s en degr{\'e} $1$ dans lesquels ${\dot\tau}_{\rho_{i,v}}$ 
agit comme $-1$. En utilisant la formule de Kunneth
on en d{\'e}duit que ${\dot\tau}_{\rho}$ agit dans 
$\RR\Gamma_c({\cal X}_\vp(a),
{\dot\A}_{\rho}\otimes h_{\vp}^*\L_\psi)$ 
comme 
$$(-1)^{\v_\vp(a_1\ldots a_{m+n-1})}=(-1)^{\v_\vp(b_1\cdots b_{n-1})}.$$

\section{D'autres remarques}

\subsection{Le signe $\varepsilon_\lambda$}

On a construit en 7.2 les rel{\`e}vements ${\tilde\tau}_\rho$ de $\tau$
sur $\tA_\rho$. Rappelons cette construction dans la situation
plus habituelle de la correspondance de Springer. D'apr{\`e}s 
2.2.3 et 2.3.1, la d{\'e}finition de ${\tilde\tau}_\rho$ dans 7.2 
et celle qui suit co{\"\i}ncident.

Soient $\g=\gl(n)$ et $\tau:\g\rta\g$ la transposition $g\mapsto\,^\t g$.  
On a $\phi\circ\tau=\phi$ o{\`u} $\phi:\g\rta Q_d$ est le morphisme
polyn{\^o}me caract{\'e}ristique.

Pour toute repr{\'e}sentation $\rho$ du groupe $\S_n$, la restriction de 
$\A_\rho$ {\`a} l'ouvert des {\'e}l{\'e}ments r{\'e}guliers semi-simples $\g_{rss}$
provient d'un syst{\`e}me local $\L_\rho$ sur $Q_{d,rss}$ par l'image inverse
de $\phi$. On {\'e}tend alors le morphisme {\'e}vidente 
$$\tau^*\phi^*_{rss}\L_\rho\rta\phi^*_{rss}\L_\rho$$
en un morphisme
$${\tilde\tau}_\rho:\tau^*\A_\rho\rta\A_\rho$$
par le prolongement interm{\'e}diaire.

Maintenant, si $\rho$ est {\'e}gale {\`a} $\rho_\lambda$ 
la repr{\'e}sentation irr{\'e}ductible
correspondant {\`a} une partition $\lambda$ de $n$, la restriction de $\A_\lambda$
au c{\^o}ne $\Nil$ des {\'e}l{\'e}ments nilpotents est 
{\`a} d{\'e}calage pr{\`e}s, le complexe d'intersection de 
${\overline{\Nil}}_\lambda$ l'adh{\'e}rence de l'orbite adjointe ${\Nil}_\lambda$ 
des {\'e}l{\'e}ments nilpotents de bloc de Jordan de taille $\lambda$ 
(\cite{Lus},\cite{B-M}).

En restreignant alors ${\tilde\tau}_\lambda$ {\`a} la fibre d'un {\'e}l{\'e}ment 
$x\in{\Nil}_\lambda({\bar k})$ qui est sym{\'e}trique, on obtient un automorphisme
d'ordre $2$ de $\Ql$. Le signe $\varepsilon_\lambda$ ainsi obtenu ne d{\'e}pend 
clairement pas du choix de la matrice $x$.

Il est naturel de se demander comment calculer $\varepsilon_\lambda$ de fa\c con
combinatoire. Nous avons obtenu le r{\'e}sultat partiel suivant.

\begin{proposition} 
Soit $V_\lambda$ la repr{\'e}sentation irr{\'e}ductible de $\S_n$ correspondant
{\`a} la partition $\lambda$. Lorsque la trace de la permutation longue
$w_0\in\S_n$ dans $V_\lambda$ est non nulle, 
son signe est {\'e}gal {\`a} celui de $\varepsilon_\lambda$ . 
\end{proposition}

\DEMONSTRATION
Soit ${\cal B}$ la vari{\'e}t{\'e} des drapeaux de $\GL(n)$. 
Rappelons que la r{\'e}solution simultan{\'e}e de
Grothendieck-Springer est d{\'e}finie par
$${\tilde\g}=\{(x,B)\in\g\times{\cal B}\mid x\in\Lie(B)\}
\ {\buildrel\pi\over\rta}\ \g$$
La transposition $\tau:\g\rta\g$
se rel{\`e}ve en une involution ${\tilde\tau}:{\tilde\g}\rta{\tilde\g}$
d{\'e}finie par ${\tilde\tau}(x,B)=(\,^\t x,\,^\t B)$ si bien qu'on a un morphisme
$${\tilde\tau}:\tau^*\RR\pi_*\Ql\rta\RR\pi_*\Ql.$$

D'apr{\`e}s le th{\'e}or{\`e}me de d{\'e}composition (\cite{BBD}), on a 
$$\RR\pi_*\Ql=\bigoplus_{\lambda}V_\lambda\otimes\A_\lambda.$$

\begin{lemme}
L'action induite par l'involution ${\tilde\tau}:{\tilde\g}\rta{\tilde\g}$ sur
$\RR\pi_*\Ql$ est {\'e}gale {\`a} celle de 
$\bigoplus_\lambda \rho_\lambda(w_0)\otimes{\tilde\tau}_\lambda$
agissant sur $\bigoplus_{\lambda}V_\lambda\otimes\A_\lambda$.
\end{lemme}

\DEMONSTRATION
Au-dessus de l'ouvert $\g_{rss}$, on a le morphisme entre
$${\tilde g}'=\{(x,gA)\in \g\times G/A\mid g^{-1}xg\in\Lie(A)\}$$
et $\tilde g$ d{\'e}fini par 
$$(x,gA)\mapsto(x,gBg^{-1})$$
est un isomorphisme.
Ici $G$ d{\'e}signe $\GL(n)$ et $A$ son sous-groupe diagonal. 
L'action de $W$ sur ${\tilde\g}_{rss}$ est d{\'e}duite de son action sur
${\tilde g}'$ d{\'e}finie par $(x,gA)\mapsto(x,gwA)$.
 
Soit maintenant $x\in A({\bar k})$ r{\'e}gulier semi-simple. 
On sait que la restriction de 
${\tilde\tau}_\lambda$ {\`a} la fibre de $\L_\lambda$ au-dessus de $x$
est l'identit{\'e}. Il suffit donc d'examiner l'action de ${\tilde\tau}$
dans la fibre de ${\tilde\g}$ au-dessus de $x$.
Dans cette fibre, $\tilde\tau$ agit par 
$$(x,wBw^{-1})  \mapsto  (\,^\t x,\,^\t (wBw^{-1}))$$
o{\`u} 
$$\,^\t wBw^{-1}=ww_0 Bw_0w^{-1},$$             
si bien qu'il agit dans $(\RR\pi_*\Ql)_x$ comme l'action de $w_0$ 
dans la repr{\'e}sentation r{\'e}guli{\`e}re d'o{\`u} de d{\'e}duite le lemme. $\carre$

{\noindent\it Fin de la d{\'e}monstration.}
On consid{\`e}re maintenant la fibre g{\'e}om{\'e}trique 
de $\pi$ au-dessus d'une matrice $x\in{\Nil}_\lambda({\bar k})$ qui est sym{\'e}trique.
Ses composantes irr{\'e}ductibles ayant la m{\^e}me dimension $d_\lambda$, 
la contribution de
$(V_\lambda\otimes\A_\lambda)_x$ dans $(\RR\pi_*\Ql)_x$ est pr{\'e}cis{\'e}ment 
le groupe de cohomologie de degr{\'e} maximal $H^{2d_\lambda}(\pi^{-1}(x))$.
Par l'application trace, on a un isomorphisme
$$ H^{2d_\lambda}(\pi^{-1}(x)){\tilde\rta}\bigoplus_{c\in C}\Ql(-d_\lambda)$$
o{\`u} $C$ est l'ensemble des composantes irr{\'e}ductibles de $\pi^{-1}(x)$.
L'action de $\tilde\tau$ sur le groupe de cohomologie
de degr\'e maximal $H^{2d_\lambda}(\pi^{-1}(x))$
se d\'eduit de son action sur cet ensemble $C$.
Une composante fix{\'e}e par $\tilde\tau$ contribue une valeur propre $1$ ;
deux composantes diff{\'e}rentes permut{\'e}es par $\tilde\tau$ contribuent une valeur propre $1$
et une valeur propre $-1$. Ainsi la trace de $\tilde\tau$ dans 
$H^{2d_\lambda}(\pi^{-1}(x))$
est toujours un nombre entier positif ou nul. 

En comparant cette assertion au lemme pr{\'e}c{\'e}dent,
on d{\'e}duit la proposition. $\square$

\subsection{Les fonctions $a'_\lambda$}

Les applications $b:{\cal H}_2^+\rta{\cal H}^+$
et $b':{\cal H}_2^+\rta{\cal H'}^+$ n'{\'e}tant pas surjectives,
les fonctions $\phi_\lambda$ (resp. $\phi'_\lambda$) n'engendrent
pas ${\cal H}^+$ (resp. ${\cal H'}^+$).

Pour toute $n$-partition $\lambda$ de $d$, la trace de l'endomorphisme
de Frobenius sur les fibres de $\tA_{\lambda,\vp}$ au-dessus des
$k$-points de $\g_{d,r,\vp}$ d{\'e}finit un {\'e}l{\'e}ment 
$a_\lambda\in{\cal H}^+$. La trace de $\Fr\circ{\tilde\tau}_\lambda$
dans les fibres de $\tA_{\lambda,\vp}$ au-dessus des points fixes
de $\Fr\circ\tau$ dans $\g_{d,r,\vp}$ d{\'e}finit un {\'e}l{\'e}ment
$a'_\lambda\in{\cal H'}^+$.

Les fonctions caract{\'e}ristiques $c_\lambda\in{\cal H}^+$ de la
double classe 
$$GL(n,\OO_\vp)\vp^{\lambda}\GL(n,\OO_\vp)
\subset\gl(n,\OO_\vp)\cap\GL(n,F_\vp)$$
forment une base de ${\cal H}^+$.
D'apr{\`e}s Lusztig (\cite{Lus}), on a
$$a_\lambda=q^{-2\<\lambda,\delta\>}(c_\lambda+
\sum_{\mu<\lambda}K_{\lambda,\mu}(q)c_\mu)$$
o{\`u} $K_{\lambda,\mu}$ sont des polyn{\^o}mes {\`a} coefficients entiers 
naturels. En particulier, les $a_\lambda$ forment une base de 
${\cal H}^+$.

Notons $c'_\lambda$ la fonction caract{\'e}ristique de l'orbite de 
$\vp^\lambda$ sous l'action de $\GL(n,\OO'_\vp)$ dans 
$S(F_\vp)\cap\gl(n,\OO_\vp)$. Ces fonctions $c'_\lambda$ forment
une base de ${\cal H'}^+$. On peut encore {\'e}crire
$$a'_\lambda=q^{-2\<\lambda,\delta\>}(\varepsilon_\lambda c_\lambda
 +\sum_{\mu<\lambda}K'_{\lambda,\mu}(q)c'_\mu)$$
o{\`u} $\varepsilon_\lambda$ est le signe d{\'e}finie 
dans la section pr\'ec\'edente.

L'{\'e}nonc{\'e} suivant se d{\'e}duit du th{\'e}or{\`e}me {\sc 4a} 
en utilisant la formule des traces de Gro\-then\-dieck.

\begin{proposition}
Pour tout $a\in A(F_\vp)$, on a
$$I(a,\alpha,a_\lambda)=(-1)^{\v_\vp(a_1\ldots a_{n-1})}J(a,\alpha,a'_\lambda).$$
\end{proposition} 

En appliquant ceet \'enonc\'e {\`a} 
$$a=\diag(\vp^d,\vp^d,\ldots,\vp^d)$$ 
on obtient de mani{\`e}re sans doute tr{\`e}s d{\'e}tourn{\'e}e 
le corollaire suivant.

\begin{corollaire}
Lorsque $\lambda=(d,d,\ldots,d)$, on a
$$\varepsilon_\lambda=(-1)^{d(1+2+\cdots+(n-1))}.$$  
\end{corollaire}

Lorsque $d=1$, la repr{\'e}sentation associ{\'e}e {\`a} $\lambda=(1,\ldots,1)$
est la repr{\'e}sentation signe. On a bien 
$${\rm Sgn}(w_0)=(-1)^{(1+2+\cdots+(n-1))}.$$
Lorsque $d=2$, $n=2$, on v{\'e}rifie que $\Tr(w_0,V_\lambda)=0$.
Ainsi le corollaire 11.2.1 n'est pas strictement contenu dans la
proposition 11.1.1.

\begin{proposition}
Notons $t:{\cal H}^+\rta{\cal H'}^+$ l'application lin\'eaire d{\'e}finie par 
$$t(a_\lambda)=a'_\lambda.$$ 
On a alors $t\circ b=b'$.
\end{proposition} 

\DEMONSTRATION
Soit $\lambda$ une $n$-partition de $d$. Notons $\rho$ la repr{\'e}sentation
induite 
$$\rho=\Ind_{\S_d\times\S_d}^{\S_{2d}}(\rho_\lambda\times\rho_\lambda)$$
On peut d{\'e}composer $\rho$ en somme de repr{\'e}sentations irr{\'e}ductibles
$$\rho=\bigoplus_{|\mu|=2d} \rho_\mu\otimes M_\mu$$
o{\`u} les multiplicit{\'e}s $M_\mu$ sont des $\Ql$-espaces vectoriels de dimension
finie. On en d{\'e}duit la d{\'e}composition
$$\tA_\rho=\bigoplus_{|\mu|=2d}\tA_\mu\otimes M_\mu.$$

L'endomorphisme de commutativit{\'e} $\kappa$ de $\tA_\rho$ pr{\'e}serve les
composantes isotypiques $\tA_\mu\otimes M_\mu$ et est donc de la forme
$$\kappa=\bigoplus_{|\mu|=2d}\Id_\tA\otimes \kappa_\mu$$
o{\`u} $\kappa_\mu$ est un endomorphisme de $M_\mu$. On a alors 
$$b_\lambda=\sum_{|\mu|=2d}\Tr(\kappa_\mu,M_\mu) A_\mu.$$

Du fait que 
$${\tilde\tau}_\rho=\bigoplus_{|\mu|=2d}{\tilde\tau}_\mu\otimes\Id_{M_\mu}$$
on a
$$b'_\lambda=\sum_{|\mu|=2d}\Tr(\kappa_\mu,M_\mu) A'_\mu$$
d'o{\`u} l'assertion.

\subsection{L'int{\'e}grale orbitale relative associ{\'e}e {\`a} $w_0$}

Identifions la permutation longue $w_0\in\S_n$ {\`a} la matrice de permutation
correspondant dans $\GL(n)$. 
Pour un {\'e}l{\'e}ment central 
$$a=\diag(\vp^d,\ldots,\vp^d)$$
pour toute fonction $\phi\in{\cal H}^+$, {\`a} la suite 
de Jacquet et Ye, posons
$$I(w_0a,\phi)=\int_{(N(F_\vp)\times N(F_\vp))/
(N(F_\vp)\times N(F_\vp))^{\scriptstyle w_0a}}
\phi(\,^\t xw_0ax')\theta(xx')\d x\d x'$$
et pour toute $\phi'\in{\cal H'}^+$, posons
$$J(w_0a,\phi')=\int_{N(F_{2,\vp})/N(F_{2,\vp})^{\scriptstyle w_0a}}
\phi'(\,^\t {\bar x}w_0a x)\theta'(x)\d x$$
o{\`u} $(N(F_\vp)\times N(F_\vp))^{w_0a}$ (resp. $N(F_{2,\vp})^{w_0a}$)
est le stabilisateur de $N(F_\vp)\times N(F_\vp)$ 
(resp. $N(F_{2,\vp})$) en $w_0a$.
Avec la normalisation habituelle attribuant la mesure $1$ {\`a} 
$N(\OO_\vp)\times N(\OO_\vp)$ (resp. $N(\OO_{2,\vp})$), on a
$$\displaylines{
I(w_0a,\phi)
=\sum_{\scriptstyle x\in N(F_\vp)/N(\OO_\vp)
\atop\scriptstyle w_0ax\in\GL(n,F_\vp)^+}
\psi(w_0ax)\theta(x)\cr
J(w_0a,\phi')
=\sum_{\scriptstyle x\in N(F_\vp)/N(\OO_\vp)\atop
\scriptstyle w_0ax\in S(F_\vp)^+}
\phi'(w_0a)\theta'(x)}$$

\begin{proposition}
Pour la matrice $a$ comme ci-dessus,
pour toute $n$-partition $\lambda$, on a
$$I(w_0a,a_\lambda)=(-1)^{d(1+2+\ldots+(n-1))} J(w_0a,a'_\lambda).$$
\end{proposition}

\DEMONSTRATION
Choisissons un entier $r>dn$. 
Consid{\'e}rons le sous-sch{\'e}ma ferm{\'e} ${\dot S}_{(d,\ldots,d),\vp}$ de 
$\g_{dn,r,\vp}$ dont l'ensemble des
$k$-points est celui des matrices de la forme $w_0x\in \g_{d,r,\vp}(k)$ avec 
$$x=\pmatrix{\vp^d & x_{1,2} & \cdots & x_{1,n}\cr 
                0  & \vp^d   & \cdots & x_{2,n}\cr
            \vdots & \ddots  & \ddots & \vdots \cr
               0   & \cdots  & 0      & \vp^d }$$
o{\`u} $x_{i,j}\in\OO_\vp/\vp^r\OO_\vp$.
Soit ${\dot h}:{\dot S}_{(d,\ldots,d),\vp}\rta\Ga$ le morphisme d{\'e}fini par
$${\dot h}(w_0x)=\sum_{i=1}^{n-1}\res_\vp(\vp^{-d}x_{i,i+1}).$$

Soit $S_{\lambda,\vp}$ le sch{\'e}ma de type fini sur $k$ d{\'e}fini en 5.1.
Le morphisme $p:{\dot S}_{(d,\ldots,d),\vp}\rta S_{(d,\ldots,d),\vp}$ 
d{\'e}fini par
$$p(w_0x)=x\OO_\vp^n$$
est lisse et {\`a} fibres g{\'e}om{\'e}triques isomorphismes {\`a} l'espace affine
de dimension fixe qu'on note $d_r$. La fonction ${\dot h}$ provient en fait d'une
fonction sur $h:S_{(d,\ldots,d),\vp}\rta\Ga$.

Alors, en utilisant la formule des traces de Grothendieck, on a
\begin{eqnarray*}
I(w_0a,a_\lambda) & =&
q^{-d_r}\Tr(\Fr,\RR\Gamma_c({\dot S}_{(d,\ldots,d)}\otimes_k{\bar k},
{\tA}_\lambda\otimes h^*\L_\psi))\cr
& = & \Tr(\Fr,\RR\Gamma_c({S}_{(d,\ldots,d)}\otimes_k{\bar k},
{\A}_\lambda\otimes h^*\L_\psi))
\end{eqnarray*}

La transposition $\tau$ laisse stable ${\dot S}_{(d,\ldots,d),\vp}$ ainsi
que la fonction $h$ si bien qu'on a
$$J(w_0a,a'_\lambda)=q^{-d_r}\Tr(\Fr\circ{\tilde\tau}_\lambda,
\RR\Gamma_c({\dot S}_{(d,\ldots,d)}\otimes_k{\bar k},
{\tA}_\lambda\otimes h^*\L_\psi)).$$

Or, on a d{\'e}montr{\'e} dans \cite{Ngo3} que pour
$\lambda\not=(d,\ldots,d)$, on a
$$ \RR\Gamma_c({S}_{(d,\ldots,d)}\otimes_k{\bar k},
{\A}_\lambda\otimes h^*\L_\psi))=0$$
d'o{\`u}
$$\RR\Gamma_c({\dot S}_{(d,\ldots,d)}\otimes_k{\bar k},
{\tA}_\lambda\otimes {\dot h}^*\L_\psi)=0.$$
Dans ce cas, on a
$$I(w_0a,a_\lambda)=J(w_0a,a'_\lambda)=0.$$

Si maintenant $\lambda=(d,\ldots,d)$, on d{\'e}duit du lemme 
2.3 de \cite{Ngo3} que le support de $\tA_\lambda$ coupe
${\dot S}_\lambda$ en un espace affine ${\dot S}_0$ dont les $k$-points sont
de la forme $w_0x$ avec
$$x=\pmatrix{\vp^d & x_{1,2} & \cdots & x_{1,n}\cr 
                0  & \vp^d   & \cdots & x_{2,n}\cr
            \vdots & \ddots  & \ddots & \vdots \cr
               0   & \cdots  & 0      & \vp^d }$$
o{\`u} les $x_{i,j}\OO_\vp/\vp^r\OO_\vp$ sont tous divisibles par
$\vp^d$. De plus, la restriction de $\tA_\lambda$ {\`a} cet espace affine 
est $\Ql$ {\`a} d{\'e}calage pr{\`e}s. 

Sur ${\dot S}_0$, on a un rel{\`e}vement {\'e}vident 
$${\tilde\tau}:\tau^*\Ql\rta \Ql$$
dont la restriction {\`a} une fibre d'un point fixe de $\tau$
est l'identit{\'e}. L'action de $\tau$ sur ${\dot S}_0$ {\'e}tant 
homotope {\`a} l'identit{\'e}, ${\tilde\tau}$ agit 
sur 
$\RR\Gamma_c({\dot S}_0\otimes_k{\bar k},\Ql)$
comme l'identit{\'e}.
Comme ${\tilde\tau}_\lambda=\varepsilon_\lambda{\tilde\tau}$
par d{\'e}finition de $\varepsilon_\lambda$, ${\tilde\tau}_\lambda$
agit sur ce complexe de cohomologie comme $\varepsilon_\lambda$.

Or, d'apr{\`e}s le corollaire 11.2.1, on a
$$\varepsilon_\lambda=(-1)^{d(1+2+\cdots+(n-1))}$$
o{\`u} 
$$I(w_0a,a_\lambda)=(-1)^{d(1+2+\cdots+(n-1))}J(w_0a,a'_\lambda).$$

En combinant avec la proposition 11.2.3, on obtient le lemme fondamental de
Jacquet et Ye pour l'{\'e}l{\'e}ment $w_0$ du groupe de Weyl.

\begin{theoreme}
Pour toute fonction $f\in{\cal H}^+_2$, on a
$$I(w_0a,b(f))=(-1)^{d(1+2+\ldots+(n-1))} J(w_0a,b'(f)).$$
\end{theoreme}

Le lemme fondamental de Jacquet et Ye dans toute sa g{\'e}n{\'e}ralit{\'e} 
concerne un {\'e}l{\'e}ment $w_M$ du groupe de Weyl qui est l'{\'e}l{\'e}ment 
le plus long du groupe de Weyl
d'un sous-groupe de Levi standard $M\supset A$ (\cite{JY}).
Nous l'avons donc d{\'e}montr{\'e} dans les deux cas extr{\^e}mes $w=1$
et $w=w_0$. Nous esp{\'e}rons d{\'e}montrer le cas g{\'e}n{\'e}ral par une sorte
de descente en combinant les id{\'e}es de d{\'e}monstration de ces deux cas.

\def\refname{{R{\'e}f{\'e}rence}}

\bigskip
Ng\^o Bao Ch\^au\\
D\'epartement de math\'ematique\\
Universit\'e de Paris-Nord\\
av. J.-B. Cl{\'e}ment\\
93430 Villetaneuse\\
FRANCE\\
ngo@math.univ-paris13.fr\\

\end{document}